\documentclass[12pt]{article}
\usepackage{amsfonts}
\usepackage{amsmath}
\usepackage{epsfig}

\title{The Farthest Point Map on the Regular Dodecahedron}
\author{Richard Evan Schwartz \thanks{Supported by N.S.F. Grant DMS-1807320}}

\newtheorem{theorem}{Theorem}[section]

\newtheorem{lemma}[theorem]{Lemma}

\newtheorem{corollary}[theorem]{Corollary}

\def\startproof{{\bf {\medskip}{\noindent}Proof: }}

\def\endproof{$\spadesuit$  \newline}

\def\C{\mbox{\boldmath{$C$}}}%
\def\Q{\mbox{\boldmath{$Q$}}}%
\def\R{\mbox{\boldmath{$R$}}}%
\def\Z{\mbox{\boldmath{$Z$}}}%

\begin{document}
\maketitle

\begin{abstract}
  Let $X$ be the regular dodecahedron, equipped with
  its intrinsic path metric.  Given $p \in X$ let
  $G(p)=-q$ where $q$ is the point on $X$ which
  maximizes the distance to $p$. (Generically, $G$ is
  single-valued.)  We give a complete
  description of the map $G$ and as a consequence
  show that the $\omega$-limit set of $G$ is the $1$-skeleton
  of a subdivision of $X$ into $180$ convex quadrilaterals.
  $G$ is a piecewise
  bi-quadratic map, and each algebraic piece is defined by a
  straight line construction involving a rhombus.
  The rhombi involved have the same shapes as the
  ones in the Penrose tiling.  Our proof is 
  computer-assisted but rigorous.
\end{abstract}

\section{Introduction}

Let $(X,d_X)$ be a compact metric space.
The {\it farthest point map\/}, or
{\it farpoint map\/} for short, associates
to each point $p \in X$ the set ${\cal F\/}_p \subset X$ of
points $q \in X$ which maximize the distance function
$q \to d_X(p,q)$.  When $X$ is a the surface of a convex
polyhedron we always take $d_X$ to be the intrinsic metric
measured in terms of paths on $X$ rather than the
chordal metric coming from $\R^3$.  The farpoint map
is pretty boring with respect to the chordal metric.

J. Rouyer's paper [{\bf R1\/}] gives a complete description of
the farthest point map on the regular tetrahedron.
My recent paper [{\bf S\/}] gives a complete description
for the regular octahedron.   The
paper [{\bf W\/}] has some results for the
case of centrally symmetric octahedra
having all equal cone angles.
The papers [{\bf R2\/}], [{\bf R3\/}] study
the farthest point map for general convex polyhedra.
The papers
[{\bf V1\/}], [{\bf V2\/}], [{\bf VZ\/}], and [{\bf Z\/}]
study the map on general convex surfaces.

Given the work in [{\bf R1\/}] and [{\bf S\/}], it is natural
to wonder about what happens for the other platonic solids.
The case of the cube and the icosahedron seem quite
similar to that of the octahedron.  The case of the
dodecahedron is the most intricate and beautiful.
I had originally planned to write about
all cases at the same time, but the
dodecahedron case already makes for a
long story.
This paper is a companion to a
Java program I wrote, which shows all
the structure.  One can get this program
on GitHub:
\newline
\newline
    {\bf http://www.github.com/RichardEvanSchwartz/Dodecahedron\/}
\newline

Henceforth $X$ denotes the regular dodehahedron equipped with
its intrinsic path metric. It is nicer to think about
the set
\begin{equation}
  {\cal G\/}_p=A({\cal F\/}_p),
\end{equation}
where $A$ is the antipodal map.
When ${\cal G\/}_p$ is a singleton, we define
$G(p)$ to be this point. Whenever we write
$G(p)$ we mean implicitly that ${\cal G\/}_p$ is a singleton.

\begin{center}
\resizebox{!}{3in}{\includegraphics{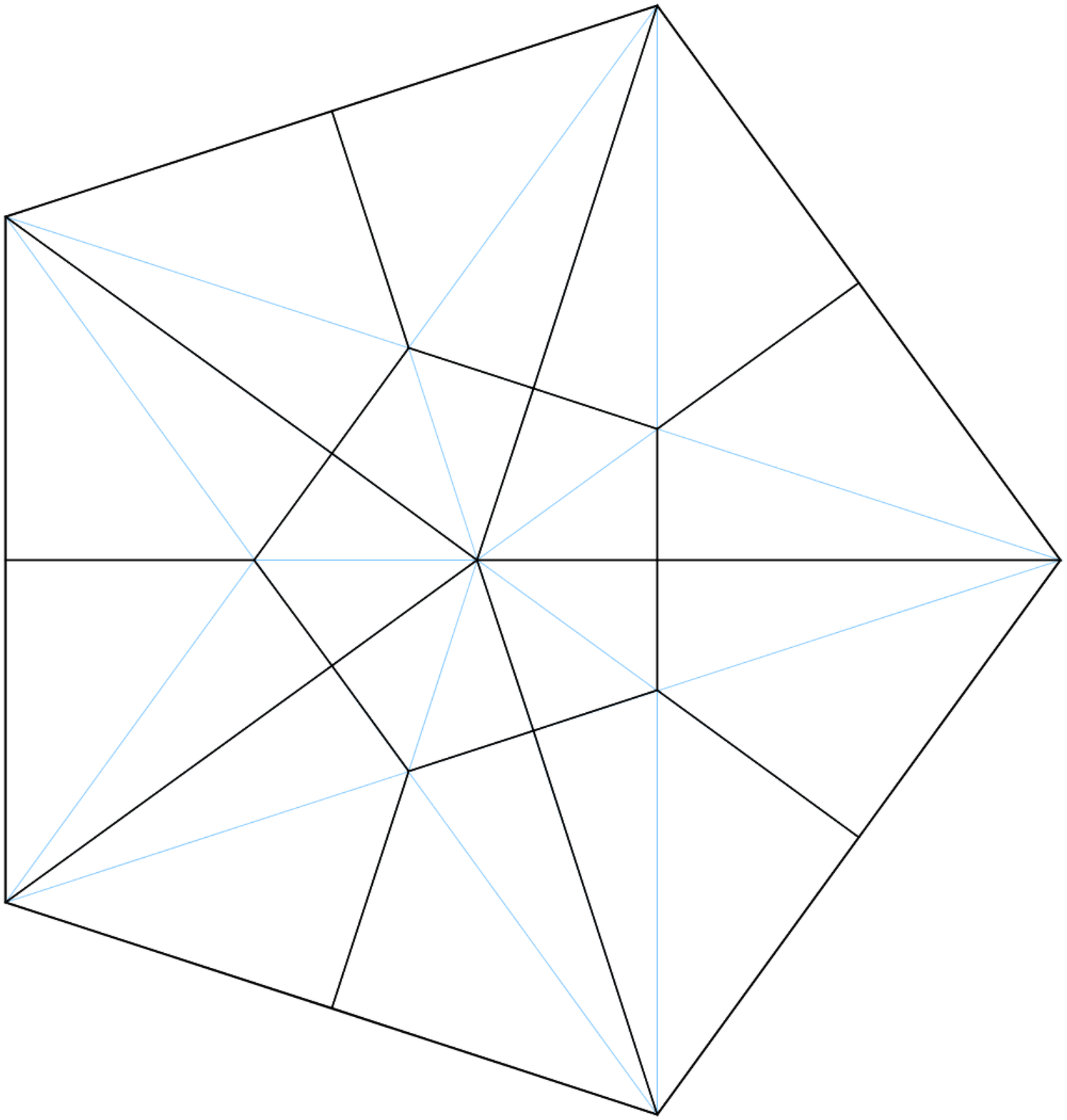}}
\newline
    {\bf Figure 1.1:\/} The decomposition of $\Pi$ into states
\end{center}

Since $G$ commutes with every isometry of $X$,
it suffices to describe the action of $G$ on
a single pentagonal face $\Pi$ of $X$.
Henceforth we take our points in $\Pi$.
We identify $\Pi$ with the planar pentagon whose
vertices are the $5$th roots of unity. Figure 1.1
shows a subdivision of $\Pi$ into $15$ quadrilaterals,
which we call {\it states\/}.  The blue segments are
drawn just as guides.  The edges of the states are the black segments.
Let ${\cal E\/}$ denote the union of the state edges.

Recall that the
$\omega$-{\it limit set\/} $L_{\omega}(G)$ of $G$ is the accumulation set of
the well-defined $G$-orbits.  Here is a corollary of our main result.

\begin{theorem}
  \label{omega}
  $G(p)=p$ if and only if $p \in {\cal E\/}$, and
  $L_{\omega}(G) \cap \Pi={\cal E\/}$.
\end{theorem}

Now we turn towards describing our main result.
We introduce a map which we call a
{\it rhombus map\/}.  A very similar map turned up in [{\bf S\/}] though
we did not study it as formally.

\begin{center}
\resizebox{!}{2.5in}{\includegraphics{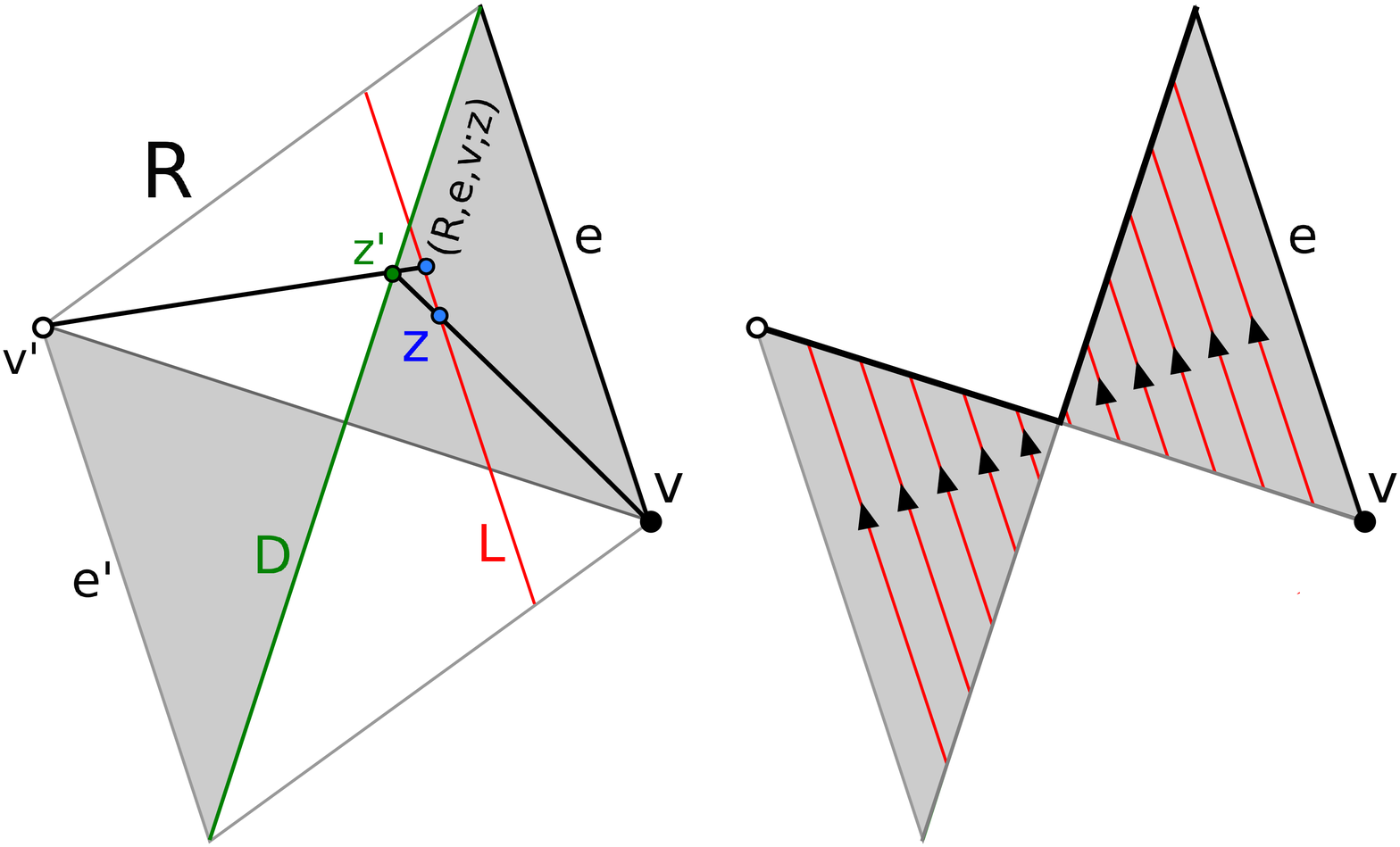}}
\newline
    {\bf Figure 1.2:\/} A rhombus map $z \to (R,e,v;z)$ and its
    invariant foliation.
\end{center}

A rhombus map is defined by a
triple $(R,e,v)$ where $R$ is a rhombus, $e$ is an edge of $R$, and
$v$ is a vertex of $R$ incident to $e$. That is, $(v,e)$ is a
flag of $R$.  Let $(v',e')$ denote the opposite flag.
Let $D$ be the diagonal of $R$ that does not contain $v$.
Let $R^o$ be the interior of $R$.
Given $z \in R^o$ let $L_z$ be the line parallel to $e$ through $z$.
We define

\begin{equation}
  (R,e,v;z)=\overline{v'z'} \cap L_z, \hskip 30 pt
  z'=\overline{vz} \cap D.
\end{equation}
Our map carries $z$ to $(R,e,v;z)$.
Figure 1.2 shows the construction.

We take the domain of $(R,e,v)$ to be the union $R^*$ of the
open shaded triangles in Figure 1.2.
The map $(R,e,v)$ is a bi-quadratic self-diffeomorphism of $R^*$
which fixes $\partial R^*-(e \cup e')$ pointwise.
The restriction of $(R,e,V)$ to
each segment of $R^*$ parallel to $e$ is a real projective
automorphism.  We foliate $R^*$ by these segments.
The orbits of $(R,e,v)$ move along the leaves of this foliation in
direction pointing away from $v$.  The attracting fixed point set
of the $(R,e,v)$ is the union of the two half-diagonals of $R$
bounding the upper white triangle. The maps $(R,e',v')$ and $(R,e,v)$ are
inverses of each other.
\newline
\newline
    {\bf Example:\/}
    We take $R=[-1,1]^2$.
and $v=(-1,-1)$ and $e$ the vertical
edge connecting $(-1,-1)$ to $(-1,1)$.
In this case, the map is given by
\begin{equation}
(R,e,v;(x,y))=  f(x,y)=\bigg(x,\frac{x^2+y}{1+y}\bigg).
\end{equation}
Every rhombus map has the form $\psi \circ f \circ \psi^{-1}$ for
some affine transformation $\psi$.
\newline

Figure 1.3 shows two special Rhombus maps which
are relevant to the dodecahedron.

\begin{center}
\resizebox{!}{3.5in}{\includegraphics{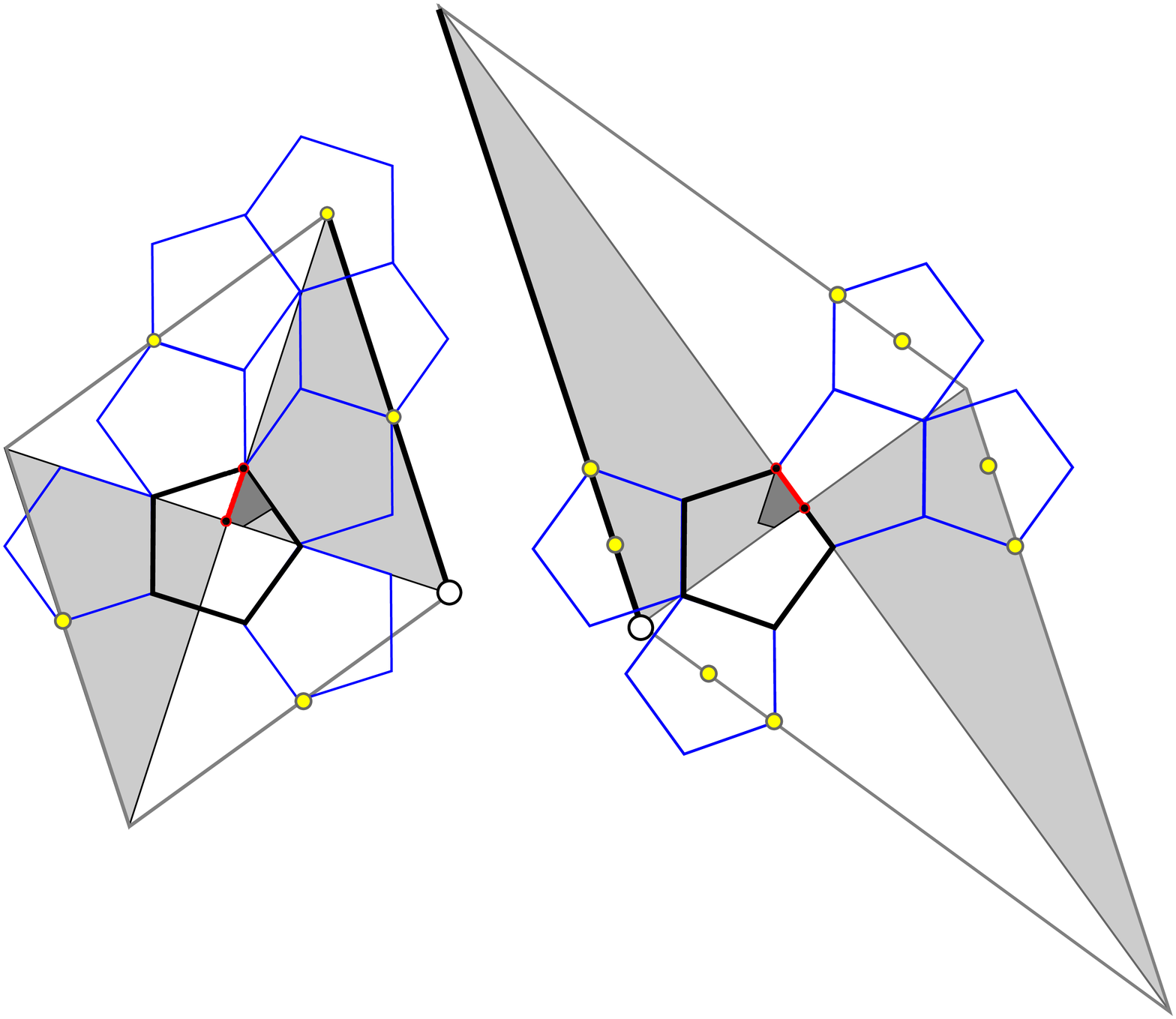}}
\newline
    {\bf Figure 1.3:\/} Two special rhombus maps
    \end{center}

These rhombs have the same shapes as those in a Penrose tiling.
The black central pentagons in Figure 1.3 are $\Pi$.  The blue pentagons are
scaffolding, designed to illustrate the construction of the rhombs.
The darkly shaded regions in the black pentagon are states.
Let $\cal R$ denote the smallest family of rhombus maps
which contains these two and which is closed with respect
to taking inverses and conjugating by the
dihedral symmetry group of $\Pi$.  The family $\cal R$
consists of $40$ maps.  We say that an
${\cal R\/}$-{\it map\/} is a map defined by
one of the members of $\cal R$.

We say that a
state $\Sigma$ is {\it adapted\/} to a
${\cal R\/}$-map $(R,e,v)$ if $\Sigma^o \subset R^*$
and if the diagonals of $R$ contain two consecutive
sides of $\Sigma$.  Each state is adapted to $4$
${\cal R\/}$-maps, and the associated foliations all coincide.
Thus we foliate $\Sigma$ by parallel line segments by
restricting the foliations of the associated maps.
The state shown in Figure 1.3 is
adapted to both of the ${\cal R\/}$-maps shown there.
Each ${\cal R\/}$-map adapted to $\Sigma$
{\it selects\/} the edge of $\Sigma$
that is contained in the attracting fixed
point set of the map.  Each edge of
$\Sigma$ is selected by a unique adapted ${\cal R\/}$-map.
The selected edges are red in Figure 1.3.

\begin{center}
\resizebox{!}{1.7in}{\includegraphics{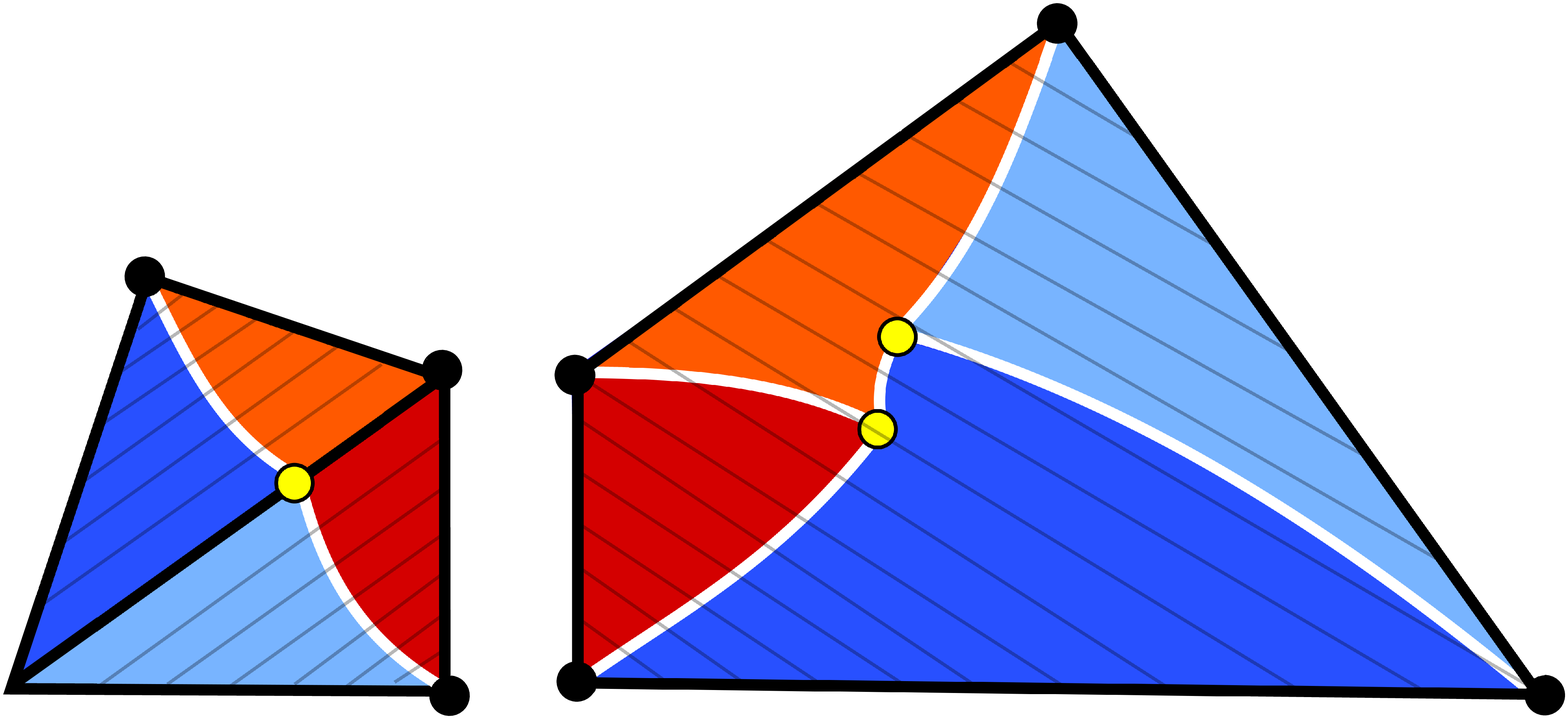}}
\newline
    {\bf Figure 1.4:\/} Combinatoral pattern in two of the cities.
\end{center}

We say that a {\it city\/} is a closed topological disk whose
boundary is a union of $3$ or $4$ algebraic arcs which are either
line segments or cubic curves. We call these segments/curves the
{\it edges\/} of the city.  Figure 1.4 shows a schematic decomposition
of two of the states into $4$ cities each.  The decomposition
on the left is meant to have bilateral symmetry.  The white
edges are nontrivial cubics and the black edges are line segments.
We insist that
each segment in the foliation of $\Sigma$ intersects the
union of edges at the red-blue interface exactly once.
One edge of each city
coincides with an edge of the state.  We call these edges
{\it external\/} and the rest {\it internal\/}.
To each city we associate the unique ${\cal F\/}$-map which
is adapted to the state and which selects the external
edge of the city.

\begin{theorem}
  \label{one}
  \label{main}
  Let $\Sigma$ be a state.
  $\Sigma$ has a decomposition into $4$ cities
  in the combinatorial pattern shown in Figure 1.4.
  If $p \in \Sigma$ then ${\cal G\/}_p$ is the
  union of the images of $p$ under the ${\cal R\/}$-maps
  associated to the cities that contain $p$.
\end{theorem}

Let us unpack and clarify this result.
Each state is isometric to one of the two shown
in Figure 1.4, so by ``a decomposition'', we mean whichever one
in Figure 1.4 corresponds to the isometry type of the state.
If $p \in \partial \Sigma$ then $p$ is fixed
  by the one or two associated maps. Hence $G$ is the identity
  on $\partial \Sigma$. If $p$ lies in the interior of some
  city, then $G(p)$ is the image of $p$ under the
  associated ${\cal R\/}$-map.
Suppose that $p$ lies in the interior of an edge $e$
common to two cities.  If $e$
  is incident to a non-right-angled
  vertex of $\Sigma$, then the two
  associated maps agree on $p$, and
  $G(p)$ is defined by either map.  If $e$
  is incident to a right-angled
  vertex of $\Sigma$, then the two
  associated maps are inverses of each other and
  ${\cal G\/}_p$ is a pair
  of points. If $e \subset \Sigma^o$ then
  ${\cal G\/}_p$ consists of $2$ points, not as clearly related to each other.
  If $p$ lies $3$ or more states, then
  then ${\cal G\/}_p$ is a pair of points.

  Now we show what the cities actually look like.
Figure 1.5 below shows the decomposition of
$\Pi$ into the $60$ cities.  One impressive thing
about the picture is that the internal city edges
are only straight line segments when they are contained
in lines of bilateral symmetry of $\Pi$.  Thus, the
red-blue interfaces look like they are straight line
segments joining non-adjacent edge midpoints of $\Pi$
but this just an illusion.  These are all arcs of
irreducible cubic curves which have the general
form given in Equation \ref{FORM} below. 
My computer program lets you zoom in and see that they are
not straight line segments.

  To give a complete account of the
  farpoint map on $X$ we need to give equations
  for the curves bounding the cities.  The
  line segment edges are all part of the framework
  shown in Figure 1.1.  The cubic edges are all
  solutions of equations of the following form
  \begin{equation}
    \label{FORM}
    \sum_{i+j \leq 3} \bigg(s_{ij} \sqrt{a_{ij}+b_{ij}\sqrt 5}\bigg) x^iy^j=0, \hskip 30 pt
    s_{ij} \in \{-1,0,1\}, \hskip 10 pt
    a_{ij},b_{ij} \in \Z.
  \end{equation}
  We give the precise formulas in \S \ref{formula}.
  The integers involved in the equations are sometimes
  surprisingly large.
  The coordinates for the yellow triple points in Figure 1.4,
  which we also list in \S \ref{formula}, also have crazy equations.

\begin{center}
\resizebox{!}{5.1in}{\includegraphics{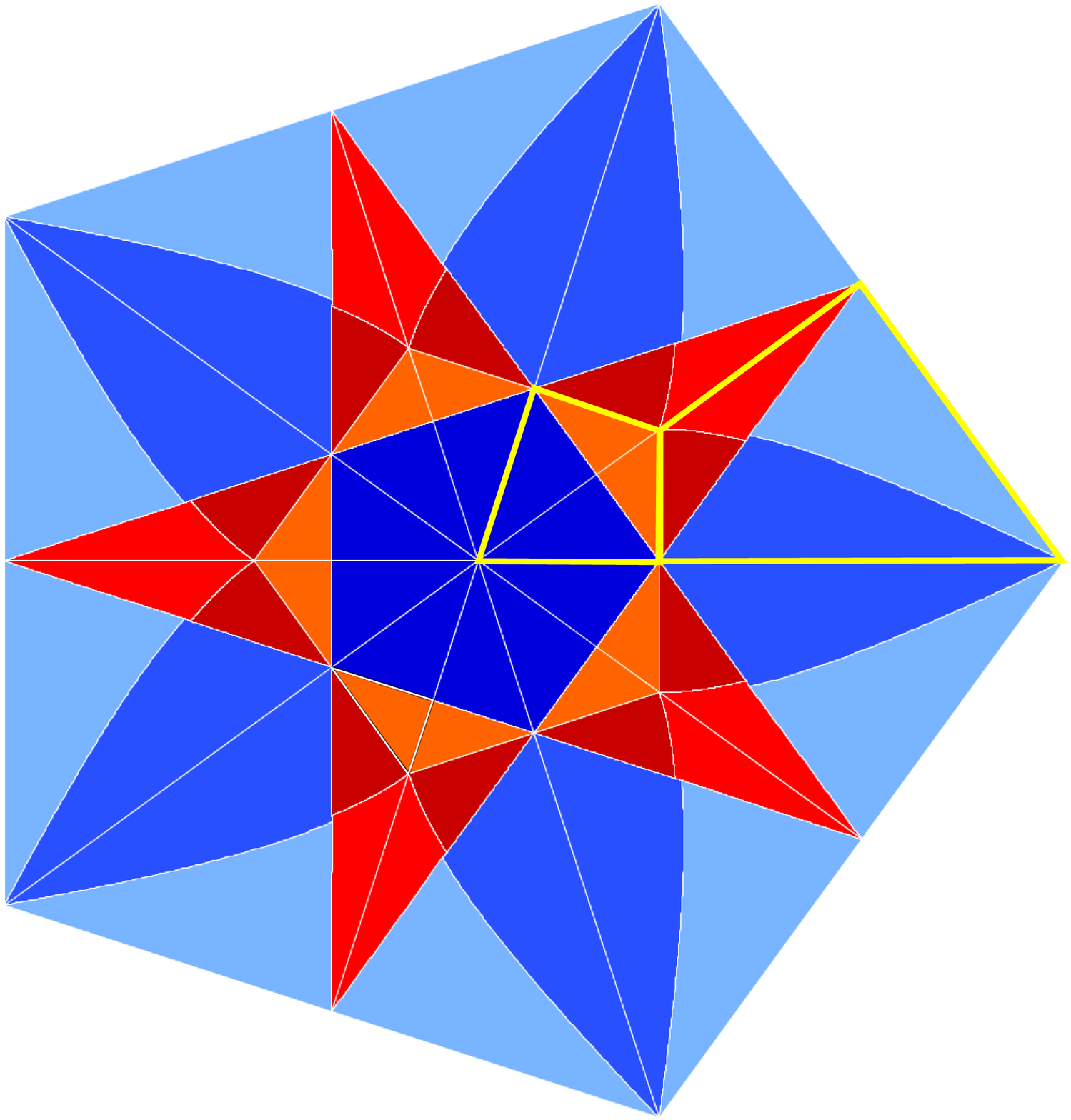}}
\newline
    {\bf Figure 1.5:\/} The decomposition into cities.
\end{center}

  \noindent
      {\bf Proof of Theorem \ref{omega}:\/}
      We see immediately from the description in
      Theorem \ref{main} that $G(p)=p$ if and only if $p$ lies in
      an edge of a state.
      Now let us understand the $\omega$-limit set.
      Inside  the state $\Sigma$, the map $G$ preserves the
  foliation by parallel line segments.  Let
  $\sigma$ be such a line segment.  This
  segment intersects the red-blue interface at
  a single point $p_{\sigma}$.  If $p \in \sigma-p_{\sigma}$
  then $G$ pushes $p$ away from $p_{\sigma}$ and
  towards the endpoint of $\sigma$.  Thus the iterates
  $p,G(p),G^2(p),...$ converge to one endpoint of
  $\sigma$ or the other, depending on the side of
  $\sigma-p_{\sigma}$ which contains $p$.  This shows
  that all orbits in $\sigma$ accumulate on the
  endpoints of $\sigma$.  But then we sweep out
  all of $\partial \Sigma$ as we vary
  $\sigma$ within the foliation.
  \endproof
       
Our proof of Theorem \ref{main}
follows a pattern similar to what we did in [{\bf S\/}],
though the details are much more involved.
I don't completely understand the huge jump
in complexity (of the proof) one sees when
going from the regular octahedron to the
regular dodecahedron but I think that a lot of it
derives from the fact that the regular
pentagon does not tile the plane whereas the
equilateral triangle does.  In the case of the octahedron
we had a global tiling which we used in order to compare
geodesic paths on the octahedron. Instead, we 
have something like a tree of possible
combinatorial types associated to a geodesic
segment on $X$, and we
resort to a computer search to tame the huge
number of combinatorial possibilities.
To give an example, there are $70$ distinct
combinatorial ways that a length minimizing
geodesic segment traveling from the bottom
face of $X$ to the top face of $X$ can
interact with the other faces of $X$.

Even after we narrow down the number of
combinatorial types we need to consider,
the algebra involved in the computations
is formidible.  Typically we consider
cubic polynomials in $2$ variables with
coefficients in a degree $8$ extension of
$\Q$. These are not polynomials that
one can just stare at and understand.
We found it easiest to let Mathematica [{\bf Wo\/}]
deal with these polynomials in an automatic
way.  The main technical gadget that powers
our proof is a positivity certificate for
polynomials in two variables, 
{\it the positive dominance criterion\/},
which we discuss in \S \ref{posdom}.

I have to admit that I am disappointed at
the length and complexity of the paper,
and I may not try to publish it.  However,
I think it is worth having a complete
proof of Theorem \ref{main} on the record.

Here is an outline of the paper.
In \S 2 we describe some preliminary notions,
such as the developing map.
In \S 3 we prove the main results modulo
technical details.  In \S 4-7 we fill in the
details of the outline.  Again, this is a
heavily computer-assisted proof which
freely makes use of the symbolic
manipulation powers of Mathematica.

In addition to getting my Java program,
the reader can also get my Mathematica code
from the same GitHub address.
The directory with the Mathematica code
has an extensive README file explaining
how to run the calculations.

I  thank In-Jee Jeong and Nathan Dunfield for
discussions about this paper.  I thank the Simons
Foundation for their support, in the form of a 2020-21 Simons
Sabbatical Fellowship.  Finally, I think the
Institute for Advanced Study for their support, in the
form of a 2020-21 membership funded by a grant
from the Ambrose Monell Foundation.

\newpage

\section{Preliminaries}

\subsection{A Spatial Argument}

As in the introduction, we let $X$ denote the regular
dodecahedron equipped with its intrinsic path metric.
$X$ is locally Euclidean except for $20$ cone
points.  The cone points each have cone angle
$9\pi/5$. As a polyhedron,
$X$ has $12$ regular pentagonal {\it faces\/}.
We identify one face $\Pi$ of $X$ with the regular
pentagon whose vertices are the $5$th roots of unity.
We think of $\Pi$ as being the bottom face.
The antipodal face $A(\Pi)$ is the top face.
Geometrically, we are normalizing so that the
distance from the center of a face of $X$ to a
vertex of that face is $1$ unit.

Almost all of our paper uses intrinsic $2$-dimensional
arguments, but there is one spatial argument we give,
in order to shorten the overall proof.  Let us do this
first.  Let  $\phi=(1+\sqrt 5)/2$.  The following facts are well known.
\begin{enumerate}
\item The diameter of any face of
  $X$ is $1+(\phi/2)$.
\item The sphere inscribed in $X$ has radius $\phi^2/2$.
\end{enumerate}

\begin{lemma}
  If $p \in \Pi$ then ${\cal F\/}_p$ is disjoint from the faces
adjacent to $\Pi$.
\end{lemma}

\startproof
Suppose this is false, and
$q \in {\cal F\/}_p$ is in a face adjacent to $\Pi$.
By Fact 1, the points
$p$ and $q$ may be connected to a path of length
at most $2+\phi<4$.  On the other hand,
and path in $X$ connecting $p$ to $A(p)$ stays outside
the inscribed sphere and has endpoints which are
antipodally placed with respect to its center.
Hence, $d_X(p,A(p)) \geq \pi \phi^2/2>4$.
This proves that $A(p)$ is farther from $p$ than is $q$.
Hence $q \not \in {\cal F\/}_p$.
\endproof

\noindent
    {\bf Remark:\/}
    We will eventually show that
    ${\cal F\/}_p \subset A(\Pi)$, but a
    crude argument like the one above would not work
    work to rule out the possibility that ${\cal F\/}_p$
    contains points in the interior of a face adjacent to $A(\Pi)$.
The problem is that the vertex antipodal to any vertex
of $\Pi$ lies both in $A(\Pi)$ and an adjacent face.
\newline

Having finished with the spatial argument, we turn to
more $2$-dimensional considerations.

\subsection{The Developing Map}

Figure 2.1 shows a combinatorial diagram for $X$.
We have $\Pi=\Pi_0$.  The antipodal face
$A(\Pi)=\Pi_{11}$ is not shown.  We have colored the
faces of $X$ according to their combinatorial
distance from $\Pi$.  In the pictures below we will
color $\Pi_{11}$ red.

\begin{center}
\resizebox{!}{2.6in}{\includegraphics{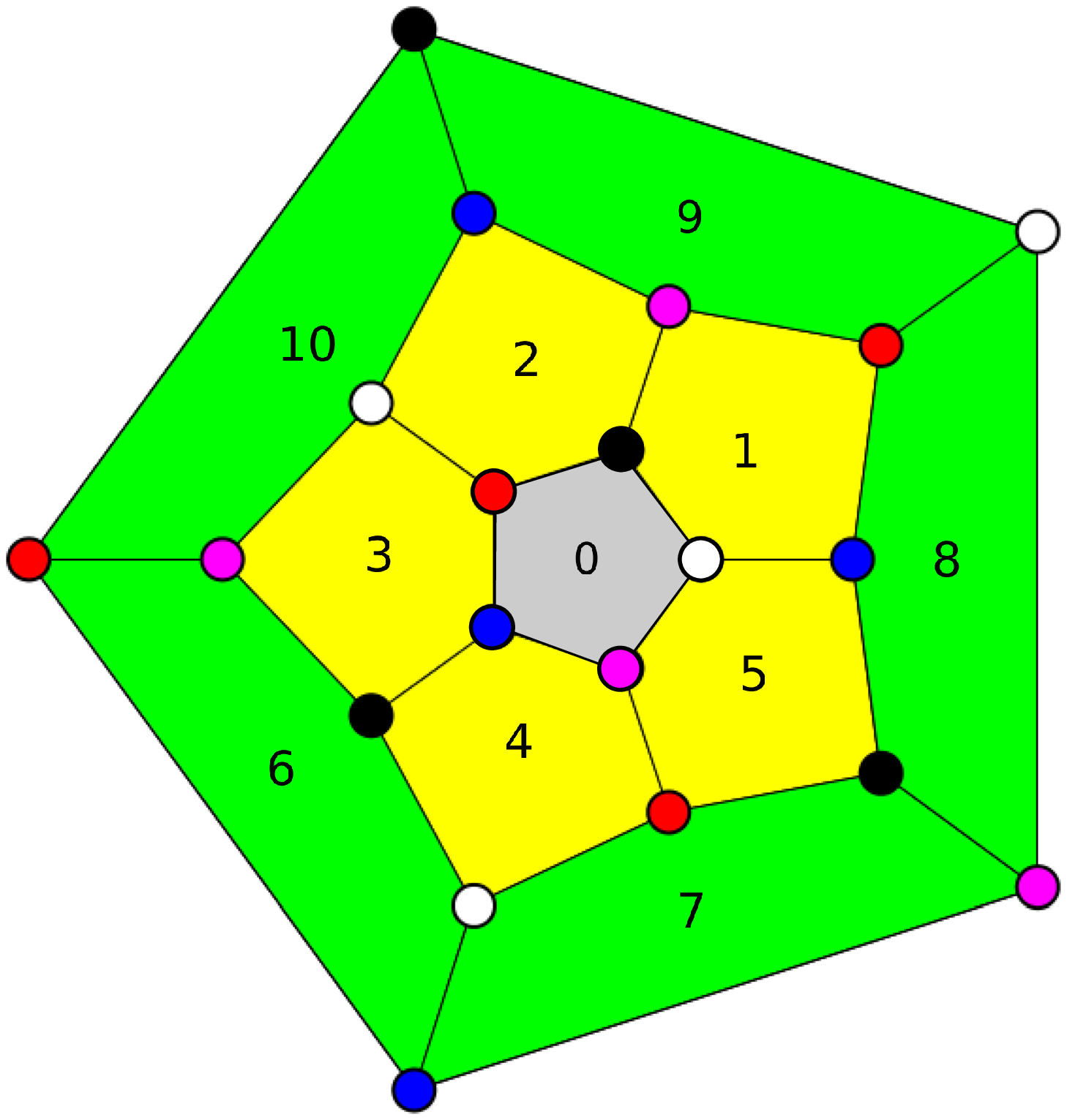}}
\newline
    {\bf Figure 2.1:\/} A diagram for the dodecahedron
\end{center}

Figure 2.1 also shows a particular $5$-coloring of the vertices.
This coloring has the property that the vertices of
the same color are the vertices of a regular tetrahedron.
This coloring will help us keep track of the orientations
of the faces when we develop $X$ out into the plane.

A geodesic segment in $X$ cannot have any cone points
in its interior.  For this reason, any geodesic segment
in $X$ is transverse to the edges of $X$ unless it
lies in a single edge of $X$.  We call such a geodesic
segment {\it transverse\/}.  We ignore the geodesic
segments which lie in a single edge of $X$ because
they never arise in connection with the farpoint map.

Let $\gamma^*$ be a transverse geodesic segment
whose initial endpoint lies in $\Pi=\Pi_0$.
There is a line segment
$\gamma \subset \C$, and an embedded union of
pentagons $$\Pi_{i_0},...,\Pi_{i_k}$$ each sharing a side
with the next, such that $\gamma^*$ and $\gamma$ have
the same length, and $\gamma \cap \Pi_0=\gamma^* \cap \Pi_0$.
To cut down on redundancies, we insist that
$i_0=0$ and that otherwise the pentagon chain is
as short as possible.  If the endpoints of
$\gamma$ lie in the interior of faces of $X$ then
the pentagon chain is unique.  The only potential
non-uniqueness arises when the initial
endpoint of $\gamma$ is a vertex of $\Pi_0$, and here
our ``shortest chain'' condition picks out a
chain uniquely in this case.  We will discuss
an example below.

This rolling process is commonly called the
{\it developing map\/}, and $\gamma$ is commonly
called the {\it developing image\/} of $\gamma^*$.
We call $\Pi_0,...,\Pi_{i_k}$ a {\it pentagon chain\/}
and we sometimes refer to it by its associated
sequence $i_0,...,i_k.$
We call the far endpoint of $q$ of $\gamma$ the
{\it terminal point\/}.  Thus, $\gamma$ is the
segment connecting $p$ to $q$, and the
distance from $p$ to $q^*$ {\it along\/}
$\gamma^*$ equals $|p-q|$.

Figure 2.2 shows the $3$ pentagon chains of length
at most $3$ which are associated to minimal
geodesic segments on $X$ connecting a point in
$\Pi_0$ to a point on $\Pi_4$. The associated
sequences are $04$ and $034$ and $054$.

\begin{center}
\resizebox{!}{1.3in}{\includegraphics{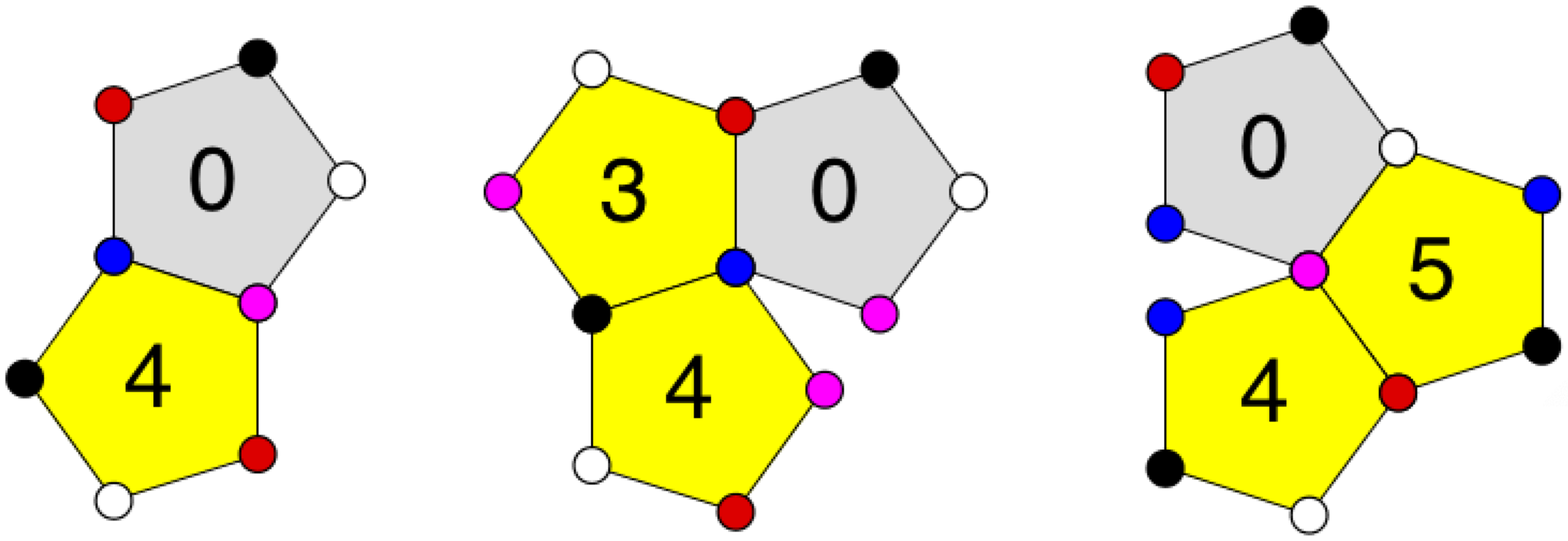}}
\newline
    {\bf Figure 2.2:\/} Chains connecting adjacent faces
\end{center}

The pentagon chains shown in Figure 2.3 do not arise in
connection with a transverse geodesic segment. The only
transverse segment the one on the left could be associated to
starts at the red vertex of $\Pi_0$, but for such geodesic segments the chain does
not have minimal length. The chain $034$ above supports the same segments
and is shorter than $0234$. The chain on the right has a similar story.  The
only potential associated geodesic segments must end at the black
vertex of $X$. The minimal chain in this case would be $036$.

\begin{center}
\resizebox{!}{1.5in}{\includegraphics{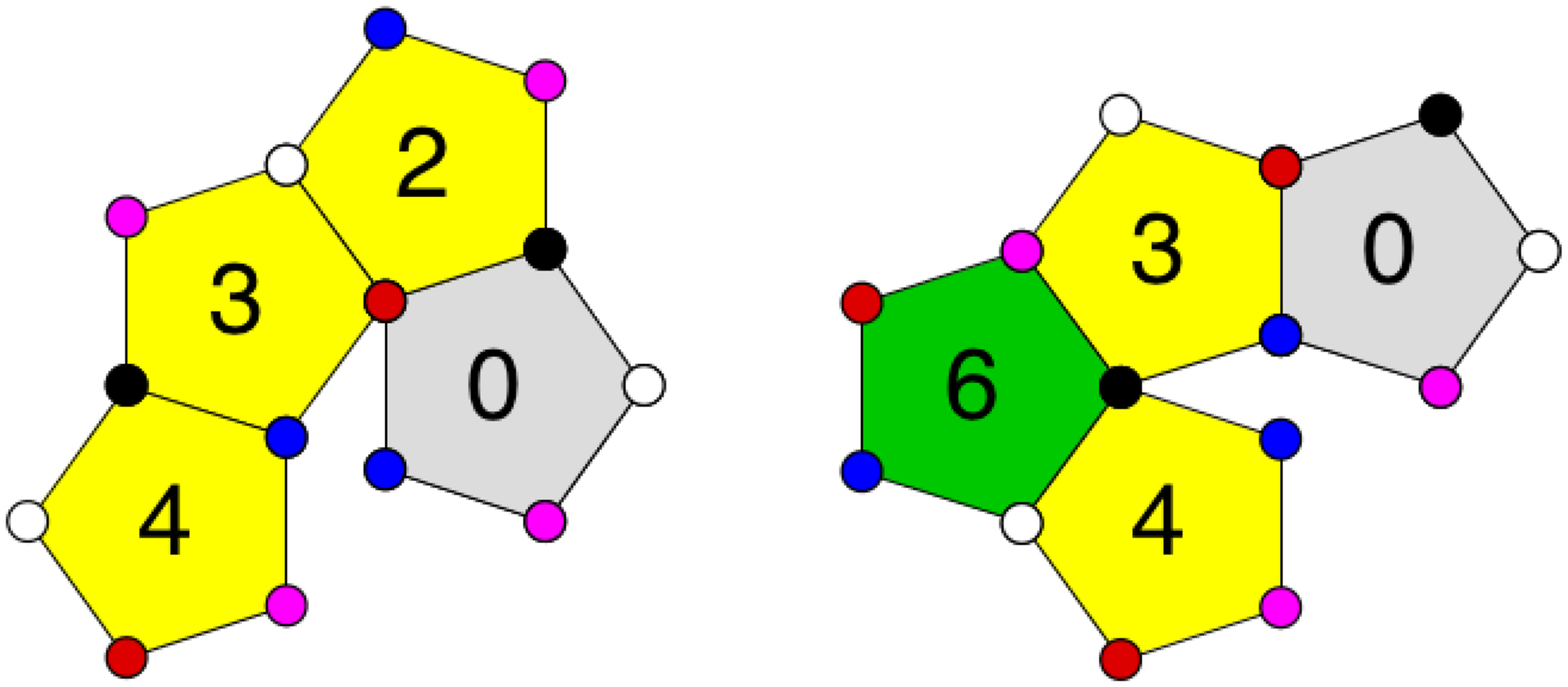}}
\newline
    {\bf Figure 2.3:\/} Two non-minimal chains
    \end{center}

\subsection{Crooked Chains}

More generally, we define a pentagon chain to be any
embedded chain of pentagons with the correct vertex
colorings.  The examples in Figure 2.3 are two
such examples.  We call such a pentagon chain
{\it straight\/} if it contains a line segment
with endpoints in the interiors of the initial
and final faces. We call such a segment
a {\it spanning segment\/}.  We call a chain
{\it crooked\/} if it has no spanning segment.
The chains in Figure 2.3 are
crooked.  A chain arises in connection with
a transverse geodesic segment on $X$ if and only
if it is straight.

One sure-fire way of generating straight chains is to
draw geodesic segments on $X$, develop them out, and then
see what chains we get.  We do not like this method because
it is hard to check that it is exhaustive.  Our approach
is to list out all possible chain sequences, from the
tree of possibilities (up to a certain length), and
then eliminate the crooked ones.  Here we explain
a computational criterion for crookedness.

Each pentagon chain defines
a finite sequence of segments in the plane,
namely the edges common to consecutive
pentagons in the chain.  We call a list of $3$ edges
{\it bad\/} if there is no line which
intersects all three.  If the pentagon chain
contains a bad triple, then it is crooked.
We can test computationally if a triple
$e_1,e_2,e_2$ is
bad in the following way.  Let $e_{k1}$ and $e_{k2}$
be the endpoints of $e_k$.  If
$e_{21}$ and $e_{22}$ both lie on the same side
of all $4$ lines $\overline{e_{1i}e_{3j}}$
then the triple is bad.

If our test does not show that a chain is
crooked it does not necessarily mean that
the chain is straight.  However, in practice,
we can see immediately that all the
remaining chains are indeed straight.

\subsection{Mirror Images}

We will generally be interested in
pentagon chains whose sequences
start with $0$ and end in either $4$, $9$, or $11$.
(We make these choices somewhat arbitrarily.)
To help us cut down on the enumeration,
we note that the symmetry $I$ of $X$
which preserves the faces
$\Pi_0,\Pi_4,\Pi_9,\Pi_{11}$ has the
following action on chains:
\begin{equation}
  \label{interchange}
  (0,1,2,3,4,5,6,7,8,9,10,11) \leftrightarrow (0,2,1,5,4,3,7,6,10,9,8,11).
\end{equation}
What we mean is that the chains associated
to the geodesics $\gamma$ and $I(\gamma)$ are swapped by the
symbolic map in Equation \ref{interchange}.
This $034$ and $035$ are swapped.  We call such
pairs of swapped chains {\it mirror images\/}.

\subsection{An Example Search}

Here we use a computer search to prove a result
which will be the basis for some other results we prove.

\begin{lemma}
  \label{adjacent}
  The only straight chains of the form $0,...,4$ associated to
  minimal geodesic segments are
    $04$ and $034$ and $054$.
\end{lemma}

\startproof
An exhaustive computer search reveals that there are
$8$ straight pentagon chains of the form $0,...,4$ which
have length $\ell \in \{4,5,6\}$.  These are
$$0,1,8,11,6,4, \hskip 15 pt
0,1,9,10,6,4, \hskip 15 pt
0,1,9,11,6,4 \hskip 15 pt
0,1,9,11,7,4$$
and their mirror images.
Figure 2.4 shows these chains.

\begin{center}
\resizebox{!}{2.8in}{\includegraphics{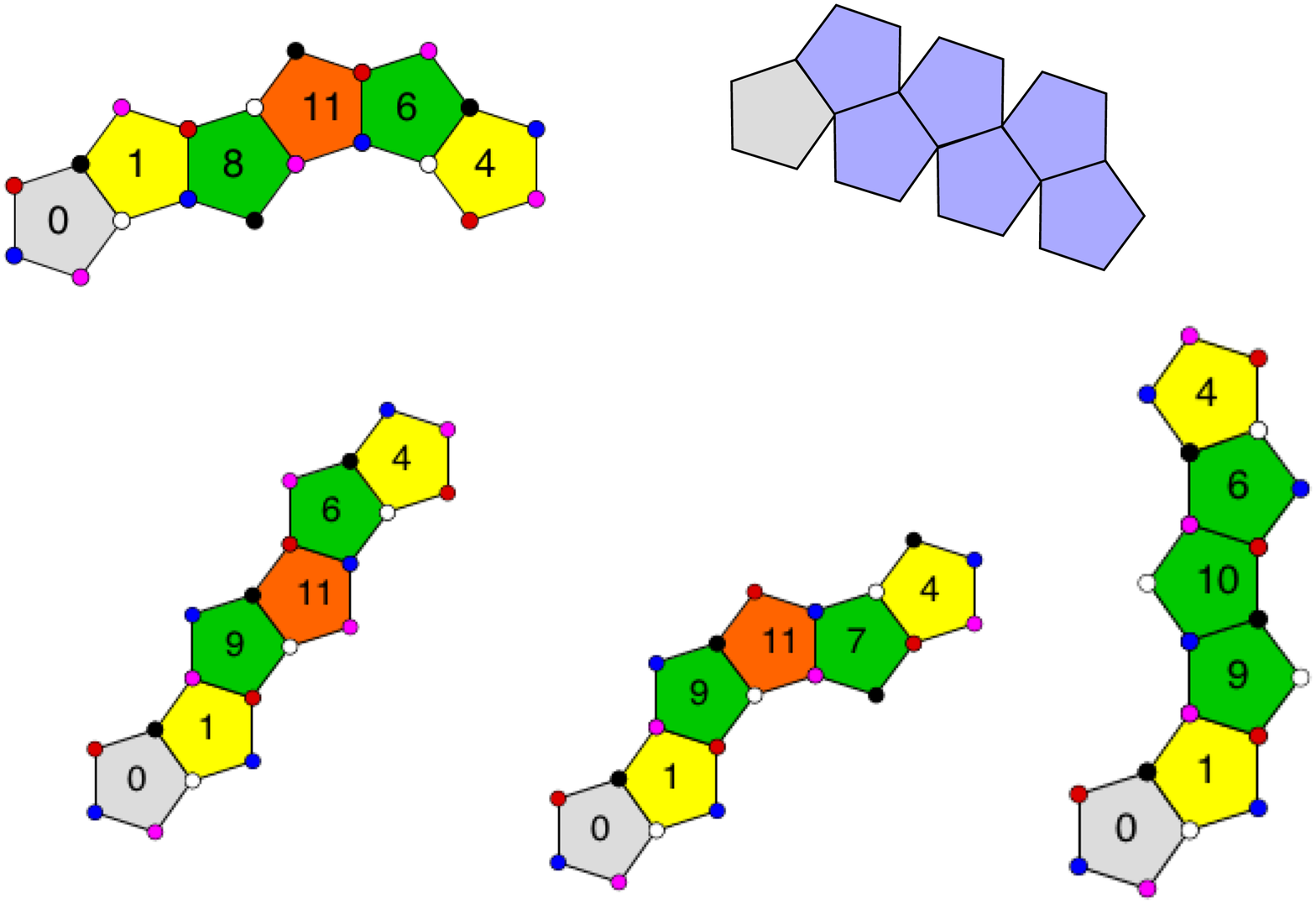}}
\newline
    {\bf Figure 2.4:\/} Two non-minimal chains
    \end{center}

As we have already mentioned, every
two points in $\Pi_0$ and $\Pi_4$ can be
connected by a path of length at most $2+\phi$.
At the same time, no spanning segment for
any of the chains in Figure 2.4 has length less than
$2+\phi$.  Hence, none of the chains in Figure 2.4
corresponds to a minimal geodesic segment.
The same goes for chains of length $7$ or more,
as one can see by considering the ``tightest'' case,
shown in Figure 2.4, in 
which every three consective pentagons
share a vertex.
\endproof

\subsection{Voronoi Decompositions}
\label{voro}

Let $H \subset \C$ denote a convex polygon.
We assume that no three vertices of $H$ are
collinear.  Thus, $H$ is convex in the
strongest possible sense.  We call
$H$ {\it strongly convex\/}.
Let $p_1,...,p_k$ be the vertices of $H$.
Given $q \in H$, let
\begin{equation}
\label{mu}
\mu_H(q)=\min_{j \in \{1,...,k\}} |q-p_j|.
\end{equation}

We say that a {\it minimal index\/} for $q$ is an index $j$ such
that $\mu_p(q)=|q-p_j|$.  The $j$th {\it voronoi cell\/} for $q$ is
the set $C_j$ of points having $j$ as one of their minimal indices.
That is, $\mu_p(q)=|q-p_j|$ if and only if $q \in C_j$.
The list $C_1,...,C_k$ is the
{\it Voronoi decomposition\/} of $H_p$.
The Voronoi cells are convex polygons.  Each
Voronoi cell has $2$ edges in $\partial H$, and
its remaining edges are contained in the union of visectors defined
by pairs of vertices in $H$.

Let $VH$ denote the Voronoi decomposition of $H$.
We say that the {\it graph associated\/} to $VH$ is
the union of the boundarties of the Voronoi cells.
This is a straight-line graph with finite valence.
Here are a few more definitions we make in
connection with this graph.
We say that an
{\it essential vertex\/} of $VH$ is a vertex of some Voronoi
cell that does not lie in $\partial H$.  Each Voronoi
cell has $3$ vertices in $\partial H$, and its remaining vertices
are essential.

We define a {\it triple point\/} to be a point that
is equidistant from at least $3$ vertices of $H$.  We name
triple points by the indices of $3$ equidistant vertices.
Every triple of indices gives rise to a triple point
because $H$ is strongly convex.
All the
essential vertices of $VH$ are triple points but some
triple points need not be essential vertices of $VH$.
We also note that there might be several valid names
for a triple point.  For instance, if $q$ is
equidistant from $p_1,p_2,p_3,p_4$, then
$(123), (124), (134), (234)$ are all valid names
for $q$.

Often we will have a $2$-parameter family
$\{H_p\}$ of strongly convex polygons,
which vary continuously depending on a
parameter $p \in U \subset \C$.  We call
such a family {\it structurally stable\/}
if all the vertices of $VH_p$ have valence
$3$ and if the combinatorics of $VH_p$ is
independent of $p$.  What this means is
these vertices never coalesce as $p$ varies
in $U$.  Put another way, structural stability
means that none of the edges shrinks to a point
as $p$ varies.

We can test for structural
stability computationally.  An edge
$e_p$ of $VH_p$ corresponds to a quadruple
$p_i,p_j,p_k,p_{\ell}$ of vertices, all
depending on $p$.  If these vertices are
never co-circular, then $e_p$ never shrinks
to a point.  We can prove this by showing
that the numerator of the
imaginary part of the cross ratio
of these points is nonzero
on $U$.  Thus we need a way to test that
polynomials in domains are positive.

\subsection{Positivity Certificates}
\label{posdom}

Here I will describe a positivity certificate.
There are many such certificates -- e.g.,
Sturm sequences in one variable, sum-of-squares
methods, the Handelman decomposition.  As far as I know,
I came up with the following one myself. It is quite
easy to implement on a computer.
I call it the {\it Positive Dominance Criterion\/}.  See my monograph
[{\bf S2\/}] for details.

We consider the $n$-variable case. Define
\begin{equation}
x^I=x_1^{i_1}...x_n^{i_n}, \hskip 30 pt
I=(i_1,...,i_n).
\end{equation}
If $I'=(i_1',...,i_n')$ we write
$I' \leq I$ if $i'_j \leq i_j$ for
all $j=1,...,n$.
Consider a polynomial
\begin{equation}
F=\sum A_I X^I, \hskip 30 pt A_I \in \R.
\end{equation}

We call $F$ {\it  positive dominant\/}
if
\begin{equation}
\label{summa}
\sum_{I' \leq I} A_{I'} \geq 0
\hskip 30 pt \forall I,
\end{equation}

\begin{lemma}
\label{PD2}
If $F$ is positive dominant then $F \geq 0$ on $[0,1]^n$.
\end{lemma}

\startproof
We first prove this result in the $1$-variable case.
We suppose that $F(x)=a_0+a_1x+...+a_nx^n$.
The proof goes by induction on the degree.
The case $\deg(F)=0$ follows from the fact that $a_0=A_0$.
Let $t \in [0,1]$.
We have
$$F(t)=a_0+a_1t+t_2t^2+ \cdots + a_nt^n \geq 
a_0t+a_1t+a_2t^2+ \cdots + a_nt^n=$$
$$t(A_1+a_2t+a_3t^2+ \cdots a_nt^{n-1})=tG(t) \geq 0$$
Here $G(t)$ is positive dominant and has degree $n-1$.
In general,
\begin{equation}
\label{slice}
F=f_0+f_1x_n+...+f_mx_n^m,
\hskip 20 pt f_j \in \R[x_1,...,x_{n-1}].
\end{equation} 
Let $F_j=f_0+...+f_j$.
Since $F$ is  positive dominant, we get that
$F_j$ is  positive dominant for all $j$.
By induction on $n$, we get $F_j \geq 0$ on
$[0,1]^{n-1}$. 
But now, if we hold
$x_1,...,x_{n-1}$ fixed and let
$t=x_n$ vary, the polynomial
$g(t)=F(x_1,...,x_{n-1},t)$
is  positive dominant.. Hence
$g \geq 0$ on $[0,1]$.
Hence $F \geq 0$
on $[0,1]^n$.  
\endproof

Now let us restrict our attention to the
$2$-variable case. (Similar remarks apply
in general, however.)
Lemma \ref{PD2} is the vanilla form of the
criterion.  Here we describe some augmentations and
variants:
\newline
\newline
{\bf Subdivision:\/}
It might turn out that $F \geq 0$
on $[0,1]^2$ but that $F$ is not positive dominant.
Given some sub-rectangle $R \subset [0,1]^2$ we
say that $F$ is {\it induced positive dominant\/}
on $R$ if $F \circ \psi$ is positive dominant
on $[0,1]^2$ for some choice of affine
isomorphism $\psi: [0,1]^2 \to R$.
In this case $F \geq 0$ on $R$.
If we want to prove that $F \geq 0$ on $[0,1]^2$
and $F$ is not positive dominant, we can check
that $F$ is induced positive dominant on
$[0,1/2] \times [0,1]$ and $[1/2,1] \times [0,1]$.
In practice this will mean checking that the
functions
$F_1(x,y)=F(x/2,y)$ and
$F_2(x,y)=(1-x/2,y)$
are both positive dominant.
We could also subdivide in the $Y$-direction.
Also, this trick can be iterated.
\newline
\newline
{\bf Triangular Domains:\/}
Sometimes we will want to know that $F \geq 0$ on
a triangle $\Upsilon$.  To do this, we produce
a polynomial map $\phi: [0,1]^2 \to \Upsilon$
and then consider the polynomial
$F \circ \phi$ on $[0,1]^2$.  Let $T_0$ be the
triangle with vertices $(0,0)$, $(1,0)$ and $(1,1)$.
The map
$\phi$ is the composition $\phi_1 \circ \phi_2$
where $\phi_1$ is an affine map from
$T_0$ to $\Upsilon$ and $\phi_1(x,y)=(x,xy)$ is
a map from $[0,1]^2$ to $T_0$.  The map
$\phi$ is a surjective polynomial map which
induces a homeomorphism from $(0,1)^2$ to
the interior $\Upsilon^o$.
\newline
\newline
{\bf Strict Positivity:\/}
Sometimes we will want to check that
$F>0$ on $[0,1]^2$.
If all the coefficient sums in
Equation \ref{summa} are
positive then we call $F$
{\it strongly positive dominant\/}.
The same argument as in Lemma \ref{PD2}
shows that $F>0$ on $[0,1]^2$ when
$F$ is strongly positive dominant.

Even if $F$ vanishes on some points on
the boundary of $[0,1]^2$ we might want to
know that $F>0$ on $(0,1)^2$.
Let $F_{\Sigma}$ denote the sum of all
the coefficients of $F$.  We call $F$
{\it solidly positive dominant\/} if $F$ is
positive dominant and $F_{\Sigma}>0$.
Essentially the same argument
as in Lemma \ref{PD2} shows that
$F>0$ on $(0,1)^2$ provided that
$F$ is solidly positive dominant.

We can combine these definitions with
our subdivision approach.  Suppose that we
suspect $F>0$ on $(0,1)^2$.  If we can
show that $F_1$ and $F_2$ above are both
solidly positive dominant it means that
$F>0$ on $(0,1)^2$ except perhaps on the
vertical segment $V=\{1/2\} \times (0,1)$.
We then test the function $F_3(y)=F(1/2,y)$
and show that it is solidly positive dominant.
This shows that $F>0$ on $V$ as well.

\newpage

\section{The Proof in Broad Strokes}

\subsection{The Antipodal Face}

In \S 4 we will prove the following result.
\begin{lemma}[Antipodal]
  Given $p \in \Pi$, we have ${\cal F\/}_p \subset A(\Pi)$.
\end{lemma}
This result involves a search through the tree of
possible combinatorial
types of length-minimizing segments on $X$.
The analogous result for the octahedron is
Statement 1 of [{\bf S\/}, Octahedral Plan Lemma].  There the
proof is easy because the regular octahedron
develops out onto a global equilateral tiling
of the plane.

Our proof of the Antipodal Lemma will reveal some
additional structure of $X$:  It will turn out
that there are $70$ pentagon chains associated
to minimal geodesics starting in $\Pi=\Pi_0$ and ending
in $A(\Pi)=\Pi_{11}$.   Among these $70$ chains,
there are $10$ of length $4$, and the remaining $60$
have length at least $5$.
Some computer experimentation reveals that each of
these $70$ chains does in fact arise in
connection with some minimal geodesic segment.

\subsection{Eliminating Combinatorial Types}

We call a geodesic segment $\gamma$
{\it straightforward\/} if it connects
a point in $\Pi$ to a point in $A(\Pi)$ and
has an associated pentagon chain of length $4$.
A glance at Figure 3,1 below shows that any
pair $(p,q) \in \Pi \times A(\Pi)$ has
some straightforward geodesic segment
connecting it.

It might be nice if we could simply say that
for every $(p,q) \in \Pi \times A(\Pi)$ the
distance $d_X(p,q)$ is realized by the
length of a straightforward geodesic segment
connecting them.  Call this property S.
Unfortunately property S can fail.
What makes our proof of Theorem \ref{main}
work is that the failures
of property S occur when $q$ is far from
${\cal F\/}_p$.   To formalize this idea,
we define
\begin{equation}
  \widehat f_X: \Pi \times A(\Pi) \to \R
\end{equation}
as follows: $\widehat d_X(p,q)$ is the minimal length
of all straightforward geodesic segments in $X$ which join $p$ to $q$.
We are not claiming that $\widehat d_X$ is a metric, and
the failure of property S tells us that
sometimes $d_X<\widehat d_X$.

Still, we can use $\widehat d_X$
to define a ``new'' dynamical system by simply substituting
$\widehat d_X$ for $d_X$ in all the definitions of the
farpoint map:
\begin{itemize}
\item $\widehat {\cal F\/}_p$ is the set of points in
$A(\Pi)$ which maximize the function $\widehat d_X(p,*)$.
\item We set $\widehat F(p)=q$ when  $\widehat {\cal F\/}_p=\{q\}$.
\item $\widehat {\cal G\/}_p=A(\widehat {\cal F\/}_p)$ and $\widehat G=A \circ \widehat F$.
\end{itemize}

In \S 4 we prove the following result:
\begin{lemma}[Comparison]
Assume that $\widehat {\cal G\/}$ has the description given
by Theorem \ref{main}.  Then $\widehat d_X(p,\widehat q)=d_X(p,\widehat q)$
when $\widehat q \in \widehat {\cal F\/}_p$.
\end{lemma}
The Comparison Lemma is the
analogue of Statement 2 of [{\bf S\/}, Octahedral Lemma],
though the proof is different.  Our proof here
involves applying our positivity certificate
to $360$ different polynomials that arise when
we compare $\widehat d_X$ and $d_X$.

Now we derive a corollary.
\begin{corollary}
  \label{compare2}
  Assume that $\widehat {\cal G\/}$ has the description given
  by Theorem \ref{main}.  Then
  ${\cal G\/}_p=\widehat {\cal G\/}_p$ for all $p \in \Pi$.
\end{corollary}

\startproof
Suppose that $p \in \Pi$ is such that
${\cal G\/}_p\not =\widehat {\cal G\/}_p$.
Then ${\cal F\/}_p \not = \widehat {\cal F\/}_p$.
Let $q \in {\cal F\/}_p$ and $\widehat q \in \widehat {\cal F\/}_p$.
Since our sets mismatch, we can assume without loss of generality that either
$q \not \in \widehat {\cal F\/}_p$ or that
$\widehat q \not \in {\cal F\/}_p$.

In the first case we have
$$d_X(p,\widehat q) \leq d_X(p,q) \leq
\widehat d_X(p,q)<\widehat d_X(p,\widehat q).$$
In the second case, we have
$$d_X(p,\widehat q) < d_X(p,q) \leq
\widehat d_X(p,q) \leq \widehat d_X(p,\widehat q).$$
In both cases, this middle inequality comes
from the fact that $d_X \leq \widehat d_X$.
Both these equations contradict the Comparison Lemma.
\endproof

In the rest of the chapter, we explain how we prove that
$\widehat {\cal G\/}$ has the description given by
Theorem \ref{main}.  Then, at the end, we
invoke Corollary
\ref{compare2} to conclude that
Theorem \ref{main} equally well describes
$\cal G$.

\subsection{The Decagon and the Hexagon}

We first describe a coloring of $X$.
We divide $\Pi$ and $A(\Pi)$ into $5$ triangles
and color them so as to be invariant under $A$.
We color the other $10$ pentagons grey.
This coloring is not so directly related to the
vertex coloring discussed in connection with
Figure 2.1, but nonetheless it is useful to us.

\begin{center}
\resizebox{!}{5in}{\includegraphics{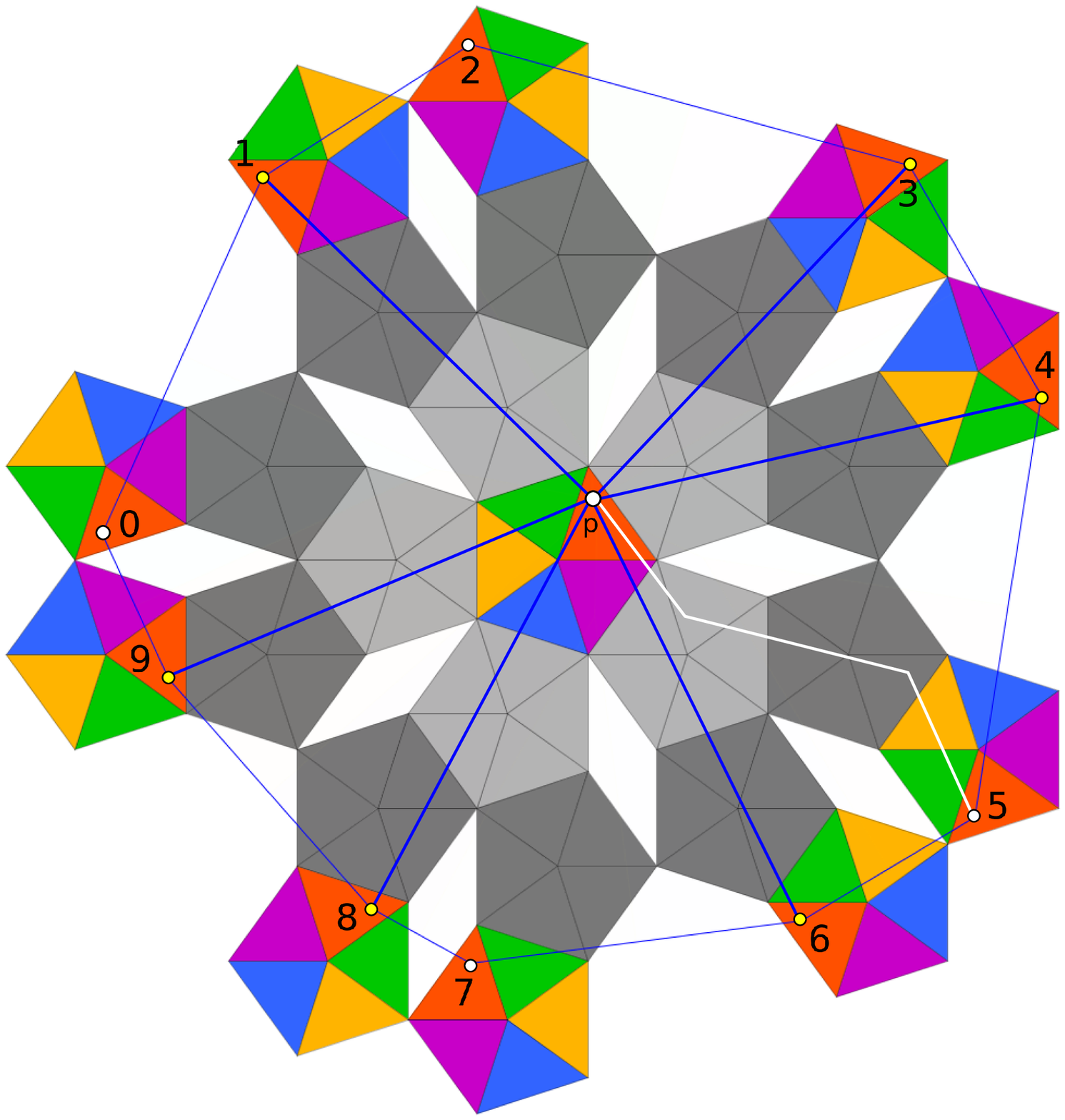}}
\newline
    {\bf Figure 3.1:\/} The dodecahedral plan.
\end{center}

Figure 3.1 shows the pentagon chains associated
to the straightforward geodesics mentioned above.
These chains have been superimposed over each other.
Let $P$ be the union of all these planar pentagons.
The inner and outer pentagons are colored so as to
be compatible with the developing map and with the
coloring of $X$ just described.

Let $P_{10} \subset P$ be the
union of the outer $10$ pentagons. 
There is a color-preserving $10$-to-$1$ map
$\Psi: P_{10} \to A(\Pi)$.
The decagon $D_p$  in Figure 3.1 has vertices 
$$\Psi^{-1}(A(p))=\{p_0,...,p_9\}.$$
The number $k$ in Figure 3.1 denotes $p_k$.
By construction and by symmetry
\begin{equation}
  \label{dec}
  \widehat d_X(p,A(q)) \geq \mu_{D_p}(q).
\end{equation}
Here $\mu_{D_p}$ is as in Equation \ref{mu}.
\newline
\newline
    {\bf Remarks:\/} \newline
    (i)  The reason we could have strict inequality is
that perhaps the line segment joining $q$
to the closest vertex of $D_p$
does not lie in $P$.  In that  case it
would not correspond to a geodesic segment in $X$.
The white zigzag in
Figure 3.1 highlights an example where
$\overline{pp_5}$ does not correspond to a
geodesic segment in $X$.
\newline
(ii)
We will not bother to prove that $D_p$ is
strongly convex, even though it is. Equation
\ref{mu} makes sense even for non-convex polygons.
\newline
(iii)  We mention
one beautiful piece of structure.  For each
index $i$, we can consider the bisector
$\beta_i$ for the points $(p_i,p_{i+5})$,
with indices taken mod $5$.  Thus
$\beta_i$ is the set of points equidistant
from these two points.  The $5$ bisectors
$\beta_0,\beta_1,\beta_2,\beta_4,\beta_5$ all
cross at $p$ and are parellel to the
$10$th roots of unity. They make a perfect
asterisk at $p$.  This does not just
follow from symmetry:  Rotation by
$\pi/5$ about $p$ is not generally a symmetry of $D_p$.
\newline

Let $\Delta \subset \Pi$ denote the central red triangle in
Figure 3.1. This triangle has vertices
$0,1,\omega$, where $\omega=\exp(2 \pi i/5)$.
Let $\Delta_k$ be the outer red triangle
labeled $k$ in Figure 3.1.   
Looking at Figure 3.1 we can see the every point in
$\Delta$ can be joined to every point of
$\Delta_k$ for $k=1,3,4,6,8,9$ by a line
segment that remains in $P$. Put another way when 
$p,q \in \Delta$ the segment
$\overline{qp_k}$ lies in $P$ for all
$k=1,3,4,6,8,9$.  This motivates us to define
$H_p$ denote the hexagon whose vertices
are $p_1,p_3,p_4,p_6,p_8,p_9$.  These are the
vertices joined to $p$ by line segments in Figure 3.1.
Given the properties of $H_p$ just mentioned, we have

\begin{equation}
  \label{hex}
  \widehat d_X(p,A(q)=\widehat d_X(q,A(p)) \leq \mu_{H_p}(q)
\end{equation}
Once again, we are referring to
Equation \ref{mu}.
\newline
\newline
\noindent
    {\bf Remarks:\/} \newline
    (i)
We might have inequality because the
minimal geodesic joining $p$ to $A(q)$ might
develop out to a line segment connecting
$q$ to a vertex of $D_p$ which is not a vertex
of $H_p$. \newline
(ii)  Here is the proof that
  $H_p$ is strongly convex for all $p \in \Delta$.
Let $A,B,C$ be three red triangles containing
consecutive vertices of $H_p$.  We can see directly
that any line $L$ that intersects both $A$ and $B$ separates
all points of $B$ from the origin.   One just has to
check the extreme cases where $L$ goes through a
vertex of $A$ and a vertex of $C$. 
\newline

We bothered to prove that $H_p$ is strongly
convex because we want to consider the
Voronoi decomposition $VH_p$.   We prove
the following in \S \ref{struct}.

\begin{lemma}[Voronoi Structure]
  Let $p \in \Delta$.
  \begin{enumerate}
  \item The essential vertices of $H_p$ lie $\Delta$.
  \item If $r$ is an essential vertex of $H_p$, then
    $\mu_{H_p}(r)=\mu_{D_p}(r).$
    \end{enumerate}
\end{lemma}
Figure 3.2 shows a typical picture of $VH_p$ for
$p \in \Delta$.  The right side shows a closeup of the
left side.  The edge between the pink and purple cells
is contained in the perfect asterisk remarked on above.

\begin{center}
\resizebox{!}{2.6in}{\includegraphics{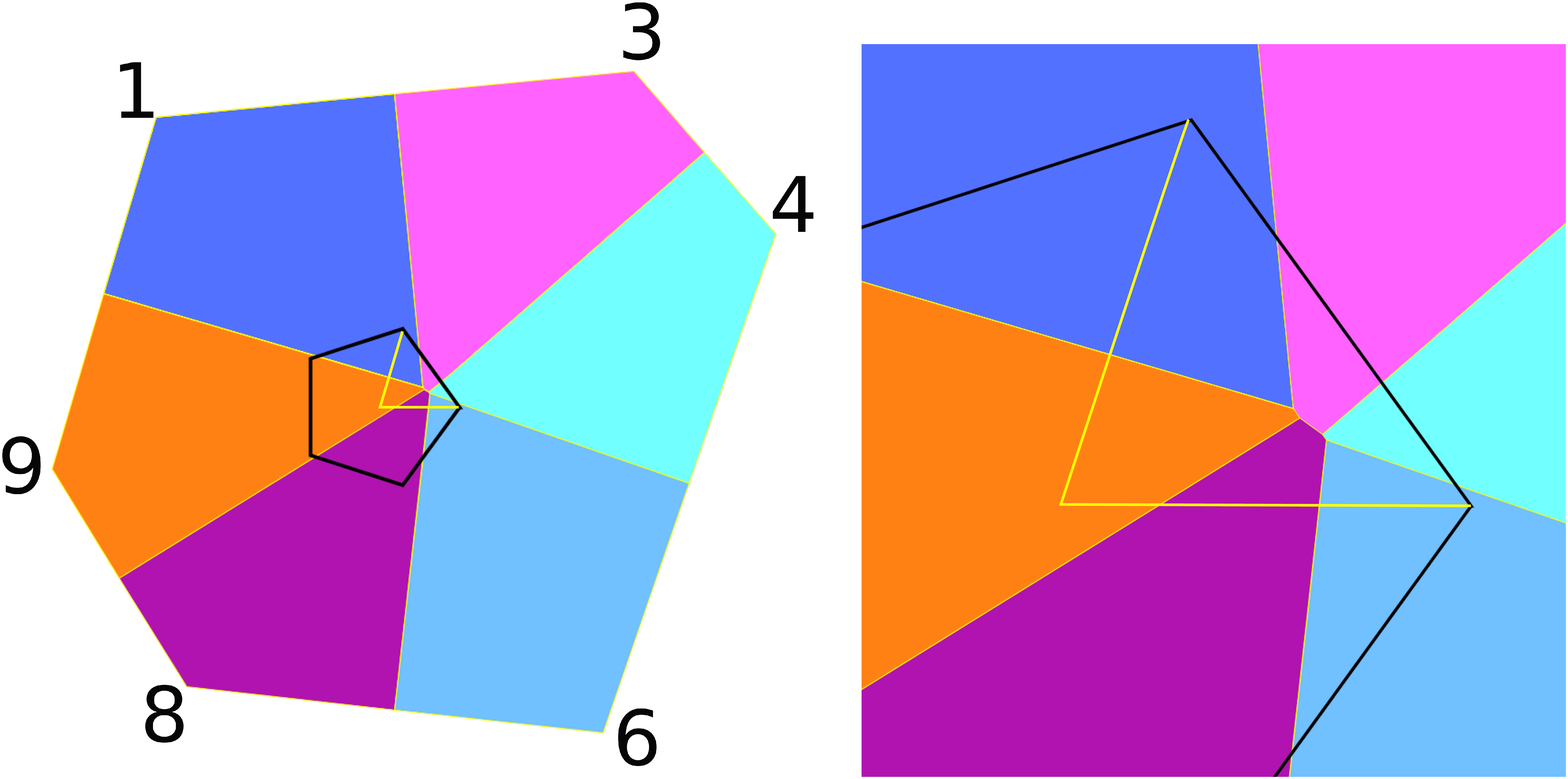}}
\newline
    {\bf Figure 3.2:\/} The Voronoi cell decomposition $VD_p$ for typical $p$.
\end{center}

\subsection{Setting up the Vertex Competition}

By definition, $\widehat {\cal G\/}_p \in \Pi$
when $p \in \Pi$.  In \S \ref{struct} we will deduce
the following result from the Voronoi Structure Lemm:

\begin{lemma}[Selection]
  If $p \in \Delta$ then $\widehat {\cal G\/}_p \in \Delta$.
\end{lemma}

Let $D_p$ and $H_p$ be the decagon and hexagon associated to
$p \in \Delta$ as in the previous section.  Recall that
$\Pi$ is the central pentagon.  We have $\Delta \subset \Pi$.

  \begin{lemma}
  \label{voroX}
  The function $\mu_{H_p}$ takes its maximum exactly on some
  sub-collection of the essential vertices of $VH_p$.
  \end{lemma}

  \startproof
  Let $q \in H_p$ be some point, not necessarily in $\Delta$.
  There is some Voronoi cell $C_i$
  such that $q \in C_i$.  The function
  $f(q)=|q-p_i|^2$ is a convex function defined
  on $C_i$, and hence it is maximized exactly
  on some collection of the vertices of $C_i$.
  If $v=p_i$ then obviously $f$ is not maximized at $v$.
  If $v \not = p_i$ is some inessential vertex of $C_i$
  then $v$ is the endpoint of
  a bisector between $C_i$ and $C_{i \pm 1}$.
  In this case, we increase $f(q)$ by pushing
  $q$ along the bisector into $H_p$.  This is
  to say that the vertices where $f$ is maximized
  are essential vertices.  But then
  $f(q) \leq f(r)$ for some essential vertex $r$,
  and the inequality is strict
  unless $q$ is also an essential vertex.
  \endproof
  
Our next result is closely related to
[{\bf R2\/}, Lemma 3] though it is stated in very different language.
  
\begin{lemma}[Vertex]
  If $q \in \widehat {\cal G\/}_p$, then $q$ is an essential vertex of $VH_p$.
\end{lemma}

\startproof
Let $\mu_{H_p}$
By the Selection Lemma, we have $q \in \Delta$.
We also have $\Delta \subset H_p$, by a wide margin.
Hence $q \in H_p$.  If $q$ is not an essential
vertex of $VH_p$ then $\mu_{H_p}(p)<\mu_{H_p}(r)$ for some
essential vertex $r$ of $VH_p$.
We have
\begin{equation}
  \label{major}
  \widehat d_X(p,A(q)) \leq \mu_{H_p}(p,q) < \mu_{H_p}(p,r)=\mu_{D_p}(p,r) \leq
  \widehat d_X(p,A(r).
\end{equation}
The first inequality is Equation \ref{hex}. The
equality is the Statement 2 of the
Voronoi Structure Lemma.  The last inequality
is Equation \ref{dec}.
The fact that
$\widehat d_X(p,A(q))<\widehat d_X(p,A(r))$
is a contradiction.
\endproof

\subsection{The Vertex Competition}

Figure 3.3 shows a close-up of the cities contained in
the triangle $\Delta$.
There are $3$ states contained in $\Delta$, which we
label $\Sigma_0,\Sigma_1,\Sigma_2$ as indicated.
Let $\Upsilon_0$ denote the half of $\Sigma_0$ lying
beneath the dotted line. Our convention is
that $\Upsilon_0$ is closed, as is $\Sigma_j$ for
$j=0,1,2$.
In proving Theorem \ref{main} it suffices by symmetry
to take $p \in \Sigma_0 \cup \Sigma_1$.  Henceforth
we do this.

\begin{center}
\resizebox{!}{4.2in}{\includegraphics{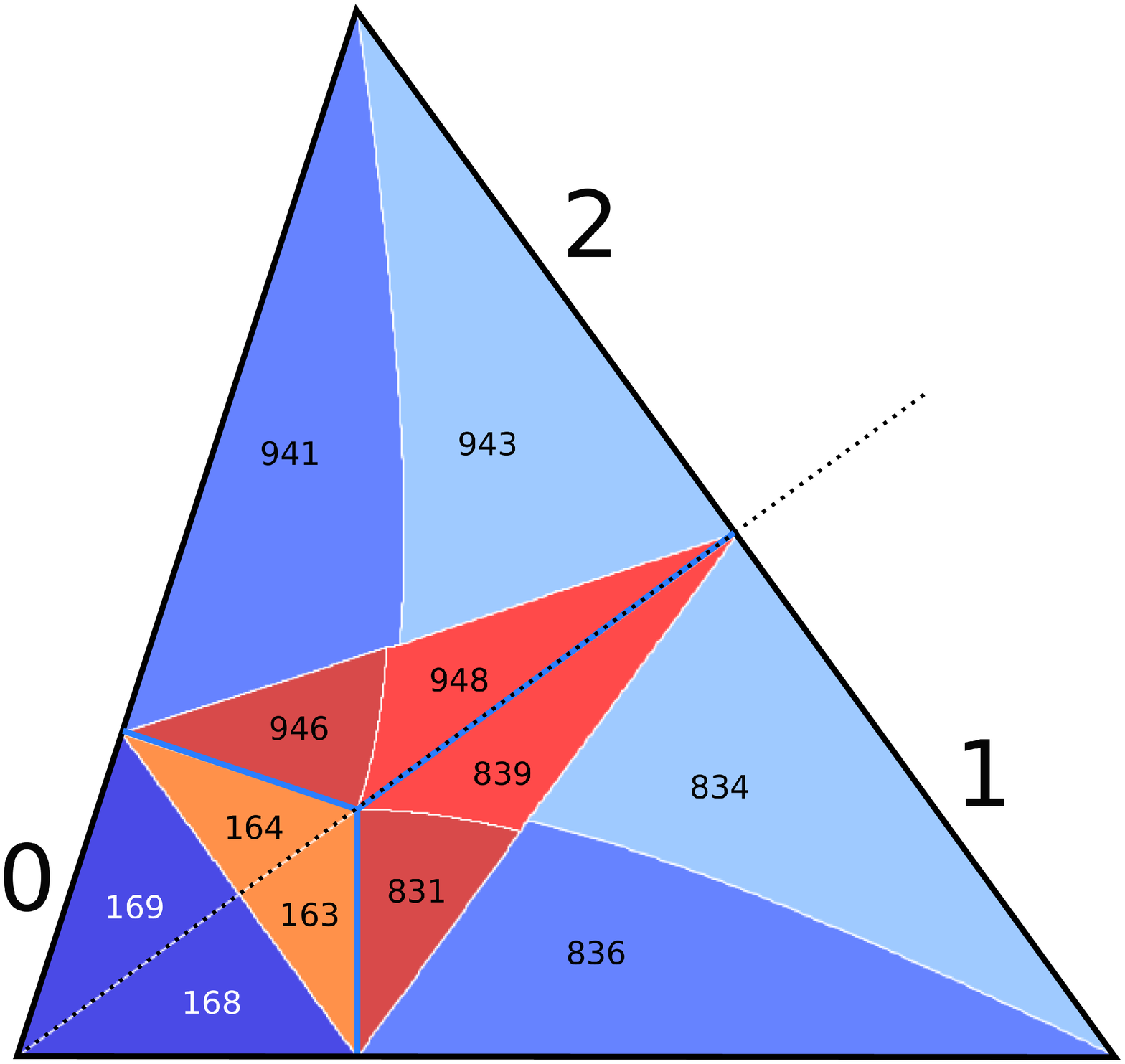}}
\newline
    {\bf Figure 3.3:\/} The cities in $\Delta$.
\end{center}

Let $(ijk,p)$ denote the unique point in $H_p$ which is
equidistant from the vertices $p_i,p_j,p_k$.  This point
may or may not be an essential vertex of $VH_p$, and
 perhaps there are other vertices of $H_p$ that have
the same distance to this point.  From the
Vertex Lemma, we know that every point
of $\widehat {\cal G\/}_p$ has the form
$(ijk,p)$ for some triple of indices and
moreover the point in question must be an
essential vertex of $VH_p$.
The rest of our proof just amounts to calculating
which triple is assigned to which point.
This boils down to algebra.
In \S \ref{struct2} we prove the following result.

\begin{lemma}[Competition]
  The following is true.
  \begin{enumerate}
    \item When $p \in \partial \Sigma_j$ for $j=0,1$ we have $\widehat {\cal G\/}_p=\{p\}$.
  \item When $p \in \Upsilon_0^o$ we have $\widehat {\cal G\/}_p \subset \{(163,p),(168,p)\}$.
  \item When $p \in \Sigma_{0}^o$ we have $\widehat {\cal G\/}_p \subset \{(163,p),(168,p),(164,p),(169,p)\}$.
  \item When $p \in \Sigma_{1}^o$ we have $\widehat {\cal G\/}_p \subset \{(831,p),(834,p),(836,p),(839,p)\}$.
  \end{enumerate}
\end{lemma}
Curiously, after all the algebra we do, this one result
has an easy geometric proof.

Directly computing all these maps we see that whenever
$\widehat {\cal G\/}_p$ is a singleton, the map
$\widehat G$ is an ${\cal R\/}$-map adapted to $\Sigma$.
Next, we
identify the domains in $\Sigma_0$ and $\Sigma_1$ which
correspond to each possible map and to verify that
we have the combinatorial structure shown in Figure 1.4.
This amounts explicit
calculations involving polynomials.  We carry this out
in \S \ref{struct2}.  

Given the Competition Lemma, the map $\widehat G$ from
Theorem \ref{main}, when defined in terms of the function
$\widehat d_X$, always has the form
\begin{equation}
\widehat G(p)=(ijk;p).
\end{equation}
According to the Competition Lemma
and symmetry there are $4$ possibilities each
within $\Sigma_0$ and $\Sigma_1$.
When we explicitly
compute all these maps, we find that the
coincide with the ${\cal R\/}$-maps described
in connection with Theorem \ref{main}.
Now we apply symmetry to get similar
results in all the states.
All this shows that Theorem \ref{main}
really does describe
$\widehat {\cal G\/}$.

  Everything we have said so far concerns
  the map $\widehat {\cal G\/}$, which is
  defined in terms of our function
  $\widehat d_X$.  But now we conclude
  from Corollary \ref{compare2} that
  ${\cal G\/}=\widehat {\cal G\/}$.  This completes
  the proof of Theorem \ref{main}.

\newpage

\section{The Antipodal Lemma}
\label{bigangle}

\subsection{The Basic Chains}

Let $p \in \Pi$ and $q^* \in {\cal F\/}_p$.
In this chapter we prove that
$q^* \in A(\Pi)$. We have already ruled out the
case that $q^*$ lies in a face of $X$ adjacent to $\Pi$.
By symmetry we just have to rule out the
possibility that $q^* \in \Pi_9-\Pi_{11}$.
We argue by contradiction.
We first consider the following $7$ pentagon chains and their mirror images.
$$
0,2,9 \hskip 15 pt
0,2,1,9 \hskip 15 pt
0,2,10,9 \hskip 15 pt
0,3,2,9 \hskip 15 pt
0,3,2,10,9$$
\begin{equation}
  \label{basic}
0,3,10,9 \hskip 15 pt
0,3,4,10,9. \hskip 15 pt
\end{equation}
We draw these pentagon chains in Figure 4.1.
In each of the first $5$ cases we add in the
magenta line which goes through the magenta vertex
of the the final pentagon and which is parallel to
the opposite side of this pentagon.  In the last
$2$ cases we draw not just this magenta line
but also the parallel blue line which goes
through the blue vertex of the final pentagon.

\begin{center}
\resizebox{!}{3.3in}{\includegraphics{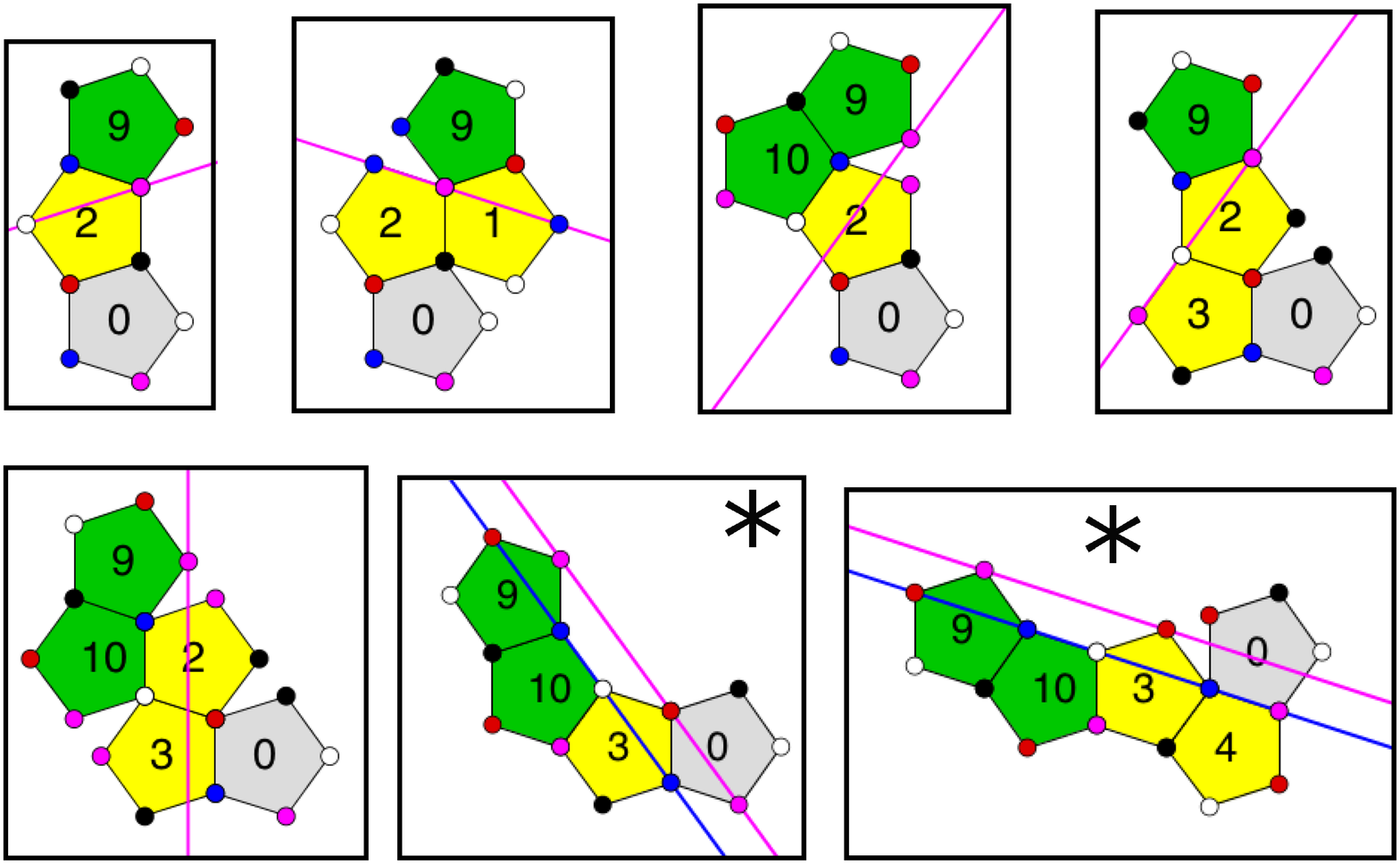}}
\newline
    {\bf Figure 4.1:\/} Seven basic pentagon chains
\end{center}

Below we will prove the following result.

\begin{lemma}
  \label{minimal}
  The only chains of the form $0,...,9$ associated to
  minimal geodesic segments are the
  $7$ basic ones and their mirror images.
\end{lemma}

\begin{lemma}
  Lemma \ref{minimal} implies the Antipodal Lemma.
\end{lemma}

\startproof
Consider the $5$ cases which are not marked by
$(*)$.  Again, $q \in \C$ is the terminal point in the chain
corresponding to $q^* \in X$.  Let
$r \in \C$ denote a point very near $q$ which we
reach by pushing $q$ away from the
magenta line and perpendicular to it.
We call this the {\it magenta variation\/}.
We have $|q-r|>|p-q|$ because the magenta
line separates all points in the the final
pentagon from all points in the
initial pentagon, $\Pi_0$.
In short, the magenta variation
increases distances.

Consider the two remaining cases, the
ones marked $(*)$.  In these cases, the
magenta variation may not increase
distances, because the magenta line
does not separate the final pentagon
from $\Pi_0$.  However, notice that in
each of these cases, $q$ cannot lie
between the blue and magenta lines,
because no line segment incident to
such a point can connect to a point in $\Pi_0$ and
yet remain in the pentagon chain.
Hence $q$ lies below the blue line.
But then the blue line separates $q$
from all of $\Pi_0$.  So, once
again, the magenta variation increases
distances for all {\it relevant\/}
choices of $q$.

By symmetry, the magenta variation also
increases distances in the $7$ mirror image chains.

Now consider all possible minimal
geodesic segments connecing $p$ to $q^*$.
Even if there is more than one such, we
can perform the magenta variation
simultaneously with respect to all the
pentagon chains.  This gives rise
to the same point $r^* \in X$ in all
case.  By compactness and Lemma \ref{minimal},
we can choose $r^*$ close enough to
$q^*$ so that each minimal geodesic
connecting $p$ to $r^*$ gives rise
to one of the pentagon chains associated
to minimal geodesics connecting $p$ to
$q$. But then each of these $p$-to-$r^*$
geodesics is longer than the corresponding
$p$-to-$q^*$ minimal geodesics.
In short, $d_X(p,r^*)>d_X(p,q)$.
This proves that $q \not \in {\cal F\/}_p$
under the assumption that Lemma \ref{minimal} is true.
\endproof

The rest of the chapter is devoted to proving
Lemma \ref{minimal}.
The general idea of our proof of Lemma \ref{minimal}
is to take the other candidate pentagon chains,
which we will discuss in the next section, and show
in each case that they contain a certain quadrilatral
$Q$ which can be replaced by a ``smaller'' quadrilateral
$\overline Q$ made from the edges of one of the
$14=7+7$ pentagon chains discussed above.
We first explain what we mean by this, and
then we carry out the analysis.

\subsection{Compressing Quadrilaterals}
\label{compress}

The edges of the pentagons in a pentagon chain have
length $\ell = 2 \sin(2 \pi/5)$.
We consider quadrilaterals $Q=(A_1,A_2,B_1,B_2)$
whose sides $A_1A_2$ and $B_1B_2$ have length $\ell$.
We call these sides {\it distinguished\/}.
These quadrilaterals need not be embedded.
Given a point $(u,v) \in [0,1]^2$ the segment
$Q(u,v)$ is the one which connects the points
$$(1-u) A_1 + u A_2, \hskip 30 pt
(1-v) B_1 + v B_2.$$
The special segments $Q(0,0)$ and $Q(1,1)$ are
the other edges of $Q$.  The special
segments $Q(0,1)$ and $Q(1,0)$ are the diagonals of $Q$.
We denote the length of $Q(u,v)$ by $|Q(u,v)|$.

Given a second quadrilateral $\overline Q$ of the
same form, we write
$Q \geq \overline Q$ if $|Q(u,v)| \geq |\overline Q(u,v)|$
for all $u,v \in [0,1]^2$.  We call $\overline Q$ a
{\it compression\/} of $Q$ in this case.

\begin{lemma}
  \label{dominance}
  $Q \geq \overline Q$ provided that
  $|Q(i,j)| \geq |\overline Q(i,j)|$ for
  all $i,j \in \{0,1\}$.
\end{lemma}

\startproof
The function $|Q(u,v)|$ is quadratic in $u$ and $v$.
Setting $v=0$ and letting $u \to \infty$ we see that
the coefficient of $u^2$ in this expression is $\ell^2$.
Likewise, the coefficient of $v^2$ in this expression is
$\ell^2$.
Therefore
\begin{equation}
 g(u,v)= |Q(u,v)|-|\overline Q(u,v)| = A u + B v + Cuv,
\end{equation}
for some constants $A,B,C$.
The restriction of $g$ to any horizontal line
in $[0,1]^2$ is a linear function.  Likewise the
restriction of $g$ to any vertical line in
$[0,1]^2$ is a linear function.
But a linear function on a segment which is non-negative at its
endpoints is non-negative on the whole segment.
Since $g(0,0), g(0,1) \geq 0$ we see that $g(0,v) \geq 0$ for all $v \in [0,1]$.
Likewise $g(1,v) \geq 0$ for all $v \in [0,1]$.  But now we restrict
$g$ to the line segment $v=v_0$.  Since $g(0,v_0), g(1,v_0) \geq 0$ we
see that $g(u,v_0) \geq 0$ as well.
\endproof

Each pentagon chain $C=(\Pi_0,...,\Pi_{i_k})$ has a {\it bottleneck quadrilateral\/}
$Q$ whose distinguished
sides are $\Pi_0 \cap \Pi_{i_1}$ and $\Pi_{i_{k-1}} \cap \Pi_9$.
Our line segment which starts in $\Pi_0$ and ends in $\Pi_9$ must
cross both edges $Q$.  If we have a second pentagon
chain $\overline C$ we say that the
{\it comparison quadrilateral\/} is the quadrilateral
$\overline Q$ having distinguished sides in the first and last pentagons
with the same vertex colors as $Q$.  The examples below will make
thie definition more clear. See e.g. Figure 4.4.
We call the second chain the
{\it comparison chain\/}.

\subsection{A Computer Search}

We call a pentagon chain {\it inefficient\/} if
its sequence $i_0,...,i_k$ has the property that
there is some index $j$ such that
$i_j \geq 6$ and $i_{j+1} \leq 5$.
Otherwise we call the chain {\it efficient\/}.

\begin{lemma}
  A pentagon corresponding to a
  length-minimizing geodesic segment is efficient.
\end{lemma}

\startproof
pentagon chain corresponds to a
length minimizing geodesic segment,
then every initial portion of the chain
does as well.  Thus, using symmetry,
we could find an inefficient chain
of the form $0,...,4$ corresponding
to a length minimizing geodesic segment.
This contradicts Lemma \ref{adjacent}.
\endproof

Now we describe the results of a $3$-step computer search.
\newline
\newline
{\bf Step 1:\/} We do a search over all efficient
pentagon chains of length at most $8$ which do not
contain $\Pi_{11}$.  We discover that all such
pentagon chains of length $8$ are bad and hence crooked.  We retain
the list of all non-bad pentagon chains and we
check that each of these is straight.  We call
such pentagon chains {\it short\/}.  We discard
the $14$ basic chains from above.
\newline
\newline
{\bf Step 2:\/} From Step 1 we see that any straight
pentagon which has no $11$ in its sequence must have
length at most $7$.
Now we do a search over all chains
of length at most $9$ whose sequence ends in $11,9$.
We retain the 
list of all non-bad pentagon chains, and we
check that all these are straight.  We call
such pentagon chains {\it long\/}.  From the
remarks about Step 1, we know that we have found all efficient
straight pentagon chains whose sequence
ends in $11,9$. 
\newline
\newline
    {\bf Step 3:\/} We merge the list of short chains with the
    list of long chains.  There are
    $38=2 \times 19$ chains on the list.  We choose one
    representative from each pair of mirror chains.   This leaves us
    with $19$ {\it candidates\/}.
\newline

We will show that no candidate can be a associated
to a distance minimizing segment in $X$ connecting
a point in $\Pi_0$ to a point in $\Pi_9$, and
this result finishes
the proof of Lemma \ref{minimal}.
The reason is that any other straight chain
is obtained from a candidate or its mirror by
appending some pentagons.  Any
minimal geodesic giving rise to this even longer
chain would have a sub-arc giving rise to a
candidate or its mirror image.

\subsection{The Isometric Cases}

Figure 4.2 shows the first $3$ candidates and the $3$ basic chains
which serve as comparison chains.  Each candidate is on the left
and the comparison chain is on the right.  We have also drawn
the bottleneck quadrilaterals on the left and the comparison
quadrilaterals on the right.  The comparison quads are not
in the same orientation as the bottleneck quads, but each
comparison quad is isometric to the corresponding bottleneck
quad in a color-preserving way.  This allows us to perform a length-decreasing
surgery on any geodesic segment $\gamma^*$ that give rise to the candidates.
We will perform the surgery using the first candidate, and the
operation works exactly the same way for the other two candidates.

\begin{center}
\resizebox{!}{5in}{\includegraphics{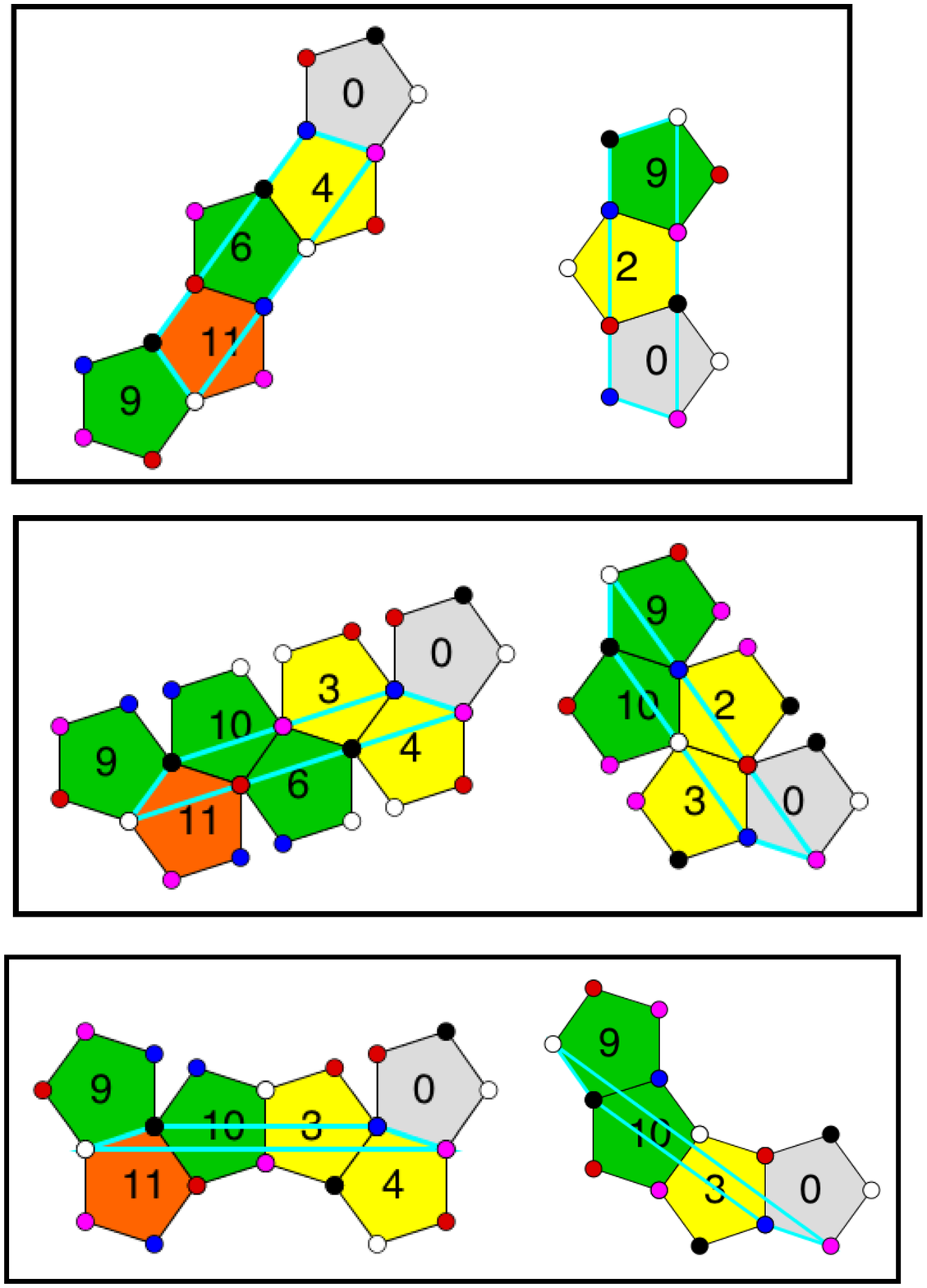}}
\newline
    {\bf Figure 4.2:\/} Three candidates and their comparison chains
\end{center}

Figure 4.3 illustrates our surgery operation.  On the left we have
the line segment $\gamma=\overline{pq}$ corresponding to $\gamma^*$.
This segment goes through points $p,d,c,b,a,q$ in order.  
The point $a$ lies just a tiny bit inside $\Pi_9$, very near the
the distinguished edge of $Q$.
The points $a,c$ are swapped by reflection in the cyan line through
the distinguished edge of $Q$.  Notice that
the blue segment $\sigma_L'=\overline{cd}$ is shorter than the black segment
$\sigma_L=\overline{ad}$.

\begin{center}
\resizebox{!}{3in}{\includegraphics{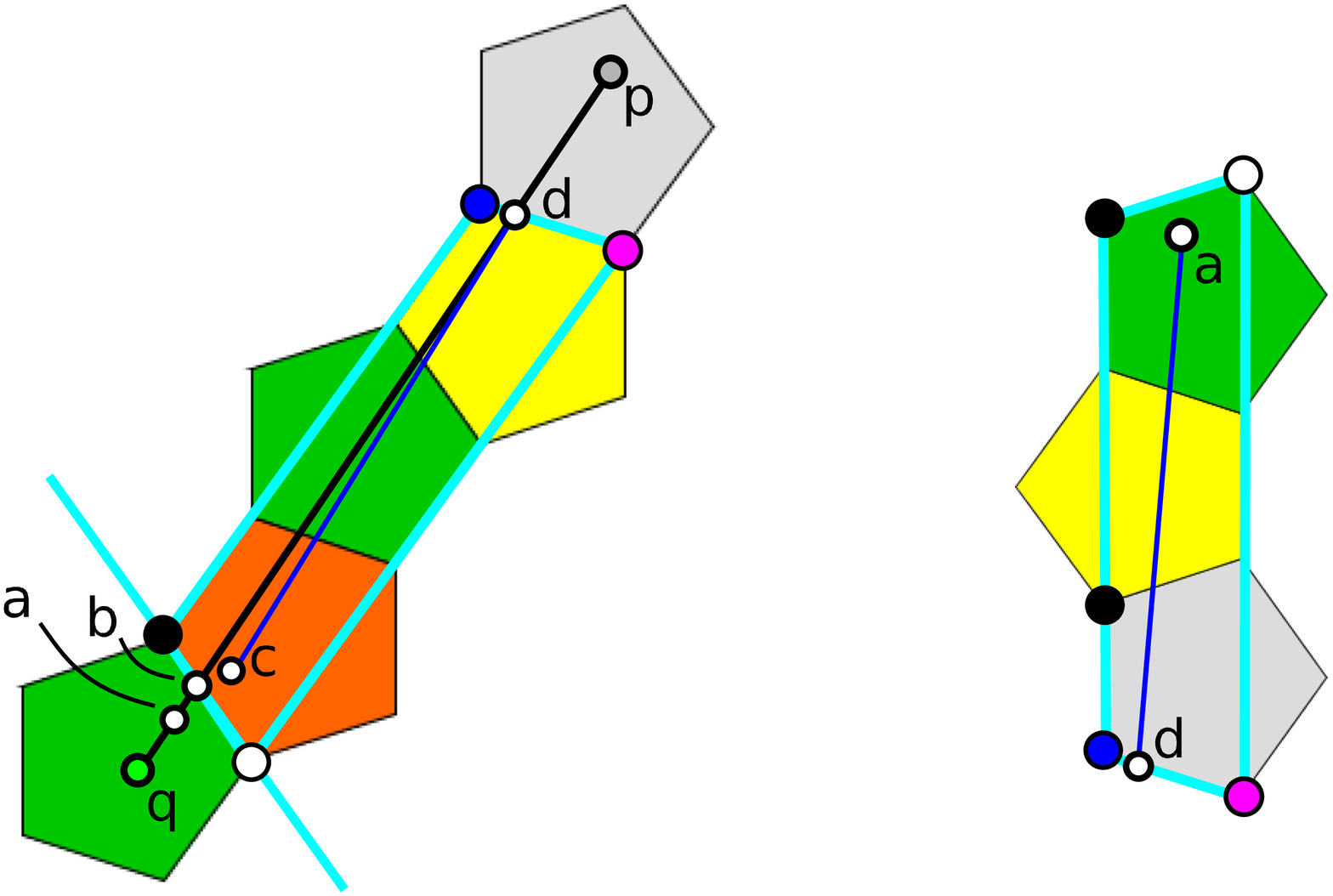}}
\newline
    {\bf Figure 4.3:\/} Three candidates and their comparison chains
\end{center}

On the right, the point $a$ is in the same position in $\Pi_9$
as is the point $a$ on the left.  In other words, the
points $a^*_{\rm left\/}$ and $a^*_{\rm right\/}$ in $X$ corresponding to these points are the same point.
Let $a^*$ be this common point.
Likewise, the $d$-point on the right is the same point
as the $d$-point on the left. Let $d^*$ be the
corresponding point in $X$.  We have $a^*,d^* \in \gamma^*$.
By construction,
the blue segment $\sigma_R=\overline{ad}$ on the right is
isometric to the blue segment $\sigma'_L=\overline{cd}$ on the left.
Hence $\sigma_R$ is shorter than $\sigma_L$.  
But the corresponding segments
$\sigma_L^*$ and $\sigma_R^*$ have the same endpoints in $X$, namely
$a^*$ and $d^*$.
Hence, if we cut out
$\sigma_L^*$ from $\gamma^*$ and replace it with
$\sigma_R^*$ we have a shorter polygonal path
with the same endpoints as $\gamma^*$.  This shows that
$\gamma^*$ is not a minimal geodesic segment.

One final word:
The only way this argument could fail is if we have
no choice of $a$ which places $c$ inside $Q$.
This happens only if $\gamma$ lies
the line extending the edge of $Q$ having white and pink
vertices.  But then $\gamma^*$ contains a cone point in
its interior and is not distance minimizing.  

Now we move on to more candidates.
Figure 4.4 shows $3$ more candidates on the left and the
basic comparison chains on the right.

\begin{center}
\resizebox{!}{6.4in}{\includegraphics{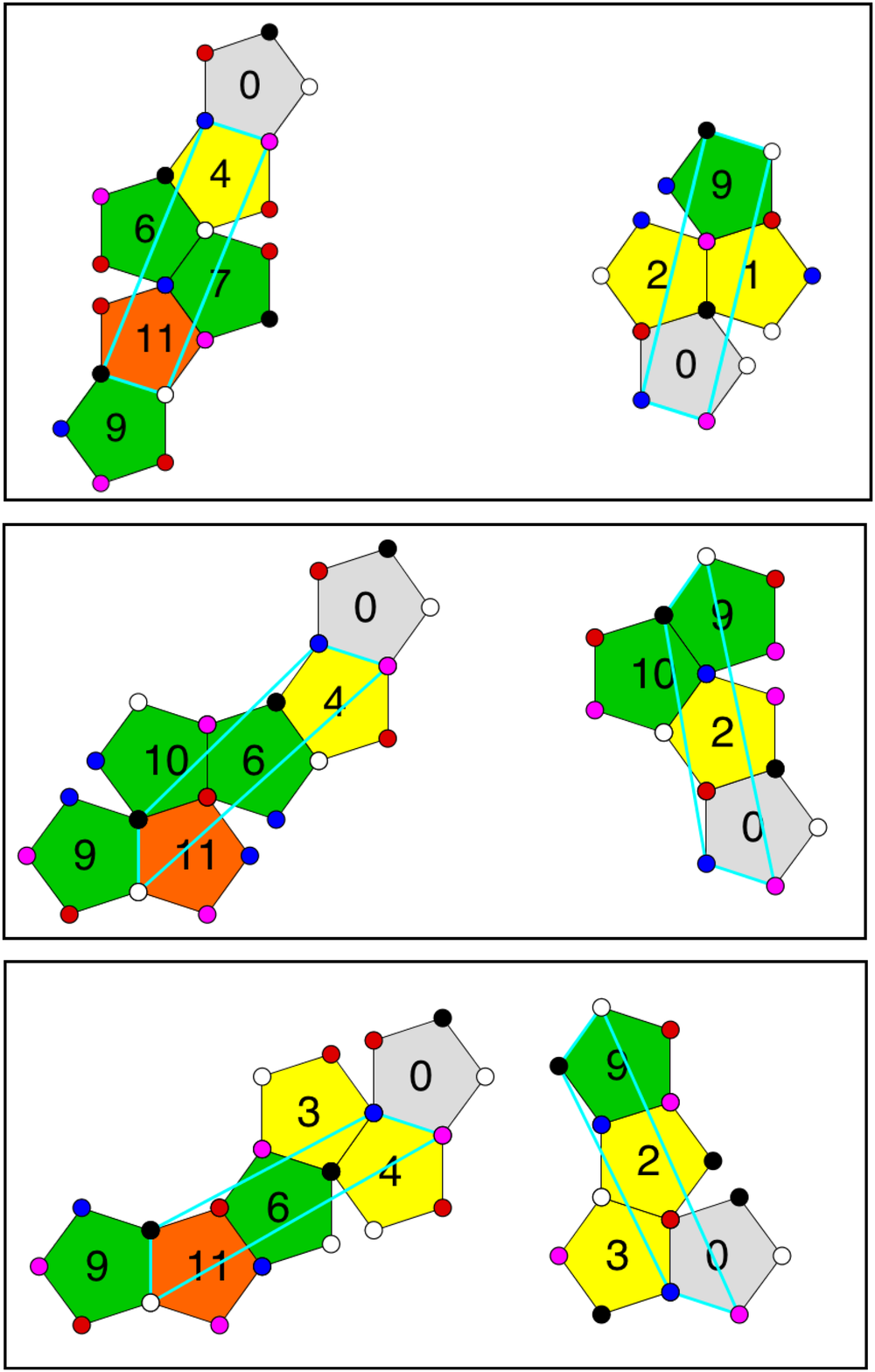}}
\newline
    {\bf Figure 4.4:\/} Three candidates and their comparison chains
\end{center}

As in the previous case, the bottleneck quads and the corresponding
comparison quads are isometric in a vertex-color-preserving way.
The difference here is that not every segment connecting
the distinguished edges in the comparison quad lies inside the
pentagon chain. Thus, we might have trouble drawing the blue
segment $\sigma_R=\overline{ad}$ on the right in Figure 4.3.  Let us look at this closely.
Say that a {\it spanning segment\/} in the bottleneck quad is
one which has its endpoints in the distinguished segments of the quad.
Call a spenning segment in the bottleneck quad {\it realizable\/} if it
lies in the pentagon chain.  Make the same definitions for the
comparison quad.

A close look at the pictures (or a computer
plot, as we did) reveals that the vertex-color-preserving
isometry from the bottleneck quad to the comparison quad
maps realizable spanning segments to realizable spanning segments.  Referring
to Figure 4.3, the spanning segment $\overline{bd}$ on the left is
realizable.  Hence, by compactness, the segment
$\overline{cd}$ lies in the comparison chain
provided that we choose $a$ sufficiently close to $b$.
The only way this could fail is if $\gamma=\overline{pq}$
contains a vertex of the chain.  But in this case, $\gamma^*$ contains
a cone point in its interior.  So, once again, we can shorten
$\gamma^*$ by surgery.

\subsection{Compressing Cases}

For the next group of candidates, the bottleneck quadrilateral
$Q$ is not isometric to the comparison quadrilateral $\overline Q$.
However, two nice things are true.
\begin{enumerate}
\item We have $Q \geq \overline Q$.  We test this using the
criterion in Lemma \ref{dominance}.
\item Every spanning
  segment of the comparison quad $\overline Q$ is realizable.
\end{enumerate}
Even though $Q$ and $\overline Q$ are not isometric, there is a
canonical correspondence between spanning segments with respect to $Q$ and
spanning segments with respect to $\overline Q$.  Corresponding segments cut
the distinguished edges at the same places -- i.e., they correspond to the
same point in the unit square from  \S \ref{compress}.
This correspondence
lets compare the segments $\sigma_L'$ and $\sigma_R$ which arise in
the surgery.

Referring to the surgery in Figure 4.3,
$\sigma_R=\overline{ab}$ on the right
might not be isometric to
$\sigma_L'=\overline{dc}$ on the left.
However, the inequality $Q \geq \overline Q$ in each case
quarantees
that $\sigma_R$ is not longer than $\sigma_L'$.
So, again, $\sigma_R$ is shorter than $\sigma_L$
and we may do our surgery operation successfully.

Figure 4.5 shows each of the candidates in this group,
and their comparison basic chains.  The reader can
see bigger pictures using our program.

\begin{center}
\resizebox{!}{7.3in}{\includegraphics{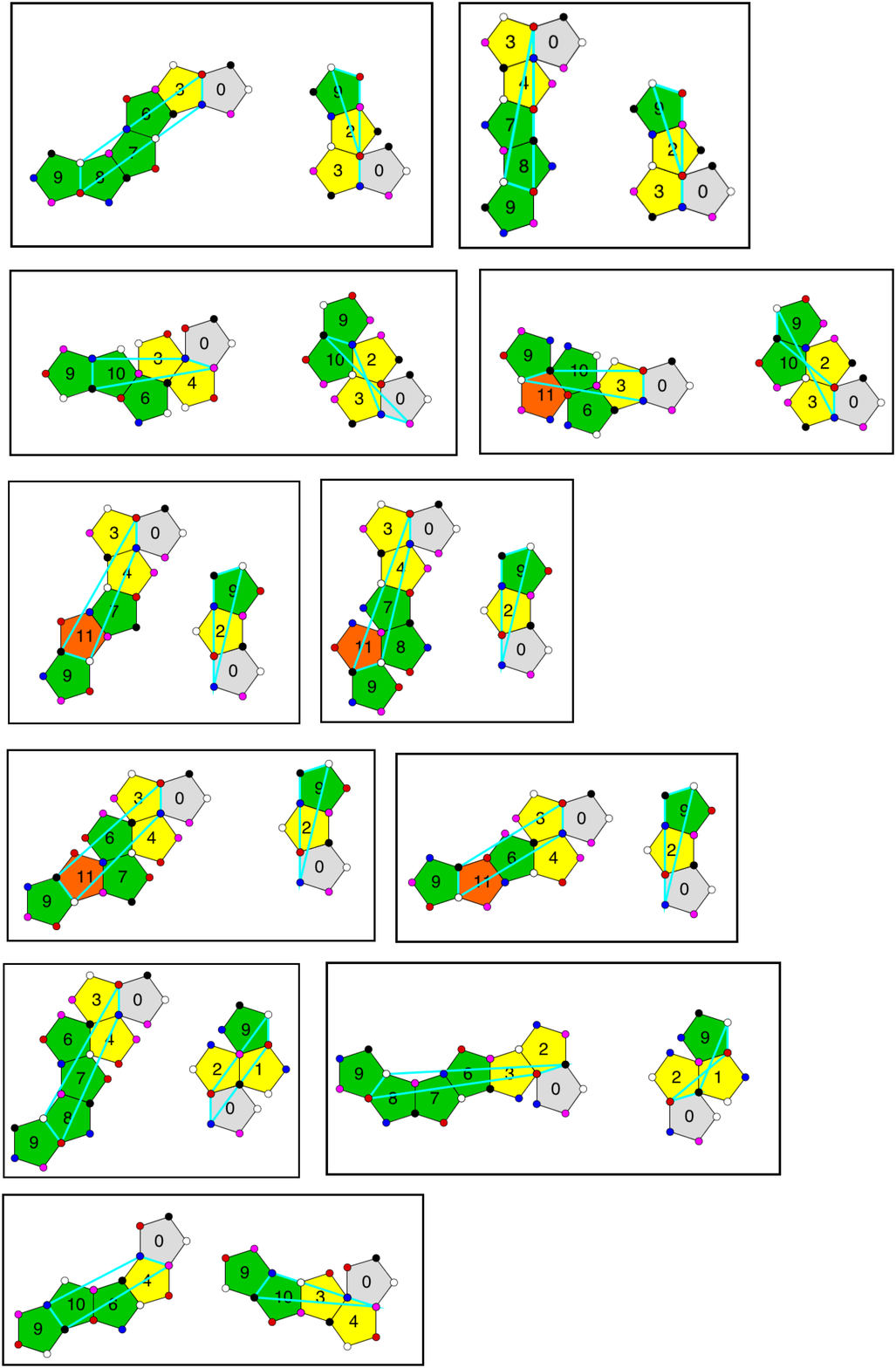}}
\newline
    {\bf Figure 4.5:\/} The candidates with compressing bottlenecks
 \end{center}

Figure 4.6 shows the next case, which is more subtle than
the ones above.  We still have $Q \geq \overline Q$.
This time not all spanning segment in $\overline Q$
are realized. 

\begin{center}
\resizebox{!}{2in}{\includegraphics{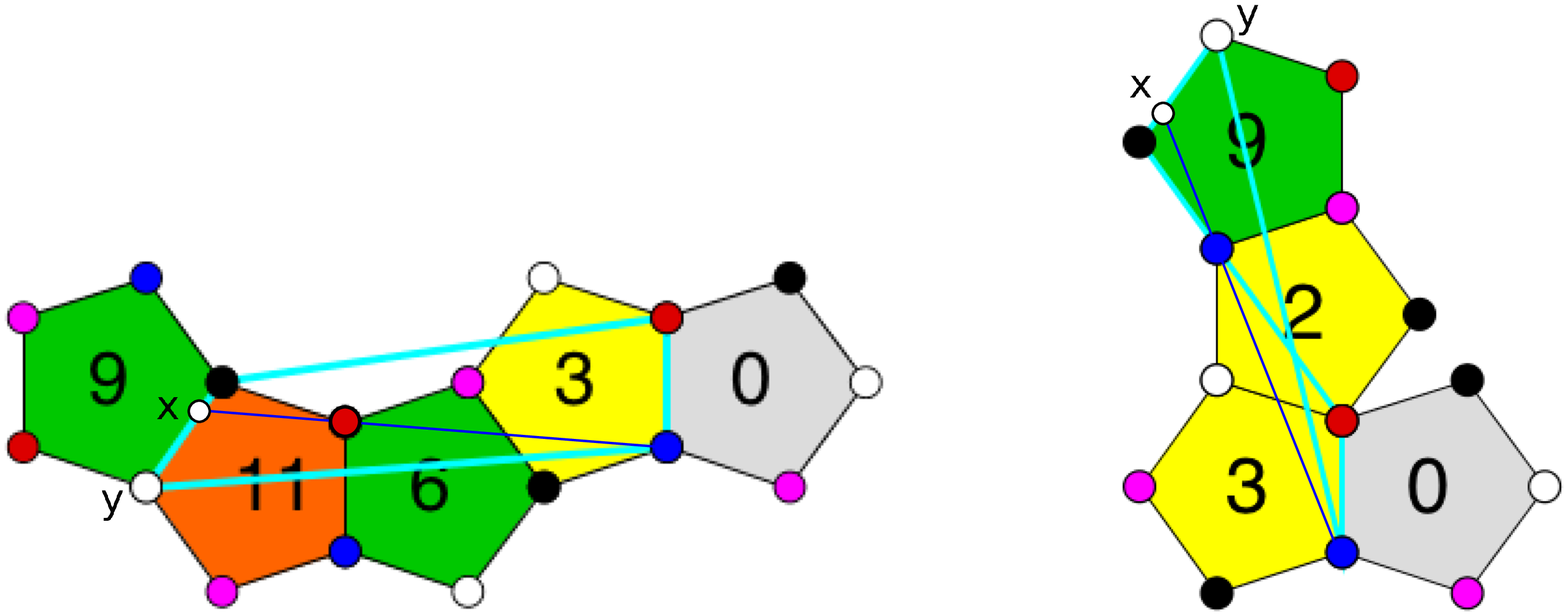}}
\newline
    {\bf Figure 4.6:\/} A subtle compressing case
     \end{center}

Here is the fact that makes this case work:
If $s$ is a spanning segment with respect to
$Q$ that is realizable then the corresponding segment
with $\overline s$ is realizable with respect to
$\overline Q$.  To see this, note that
the point $x$ in both figures lies in the same
position with respect to the distinguished segment that
contains it.  Any realizable spanning segment with respect
to $Q$ must have an endpoint in $\overline{xy}$.  At the
same time, every spanning segment with respect to
$\overline Q$ having an endpoint in $\overline{xy}$ is
realizable.  That is, having an endpoint in
$\overline{xy}$ is necessary for spanning on the
left, and sufficient for spanning on the right.
Given this property of spanning segments, the surgery
operation works just as it does for the other cases
with a compressing bottleneck.

\subsection{One Final Case}

The final case is the trickiest one.
Figure 4.7 shows the final candidate and
a comparison chain for which $Q \geq \overline Q$.
Unfortunately, our luck runs out.
All spanning segments with respect to $Q$ are
realized, but some spanning segments with respect
to $\overline Q$ are not realized.
This means that we cannot always do our
surgery.  We need to scramble to get this
case to work.

Let us choose our colors so that the coordinates
on $[0,1]^2$, corresponding to spanning segments,
have the following meaning.
\begin{itemize}
  \item $(0,0)$ is the segment joining the red vertex to the white vertex.
  \item $(1,0)$ is the segment joining the blue vertex to the white vertex.
  \item $(0,1)$ is the segment joining the red vertex to the black vertex.
  \item $(1,1)$ is the segment joining the blue vertex to the black vertex.
\end{itemize}

\begin{center}
\resizebox{!}{1.8in}{\includegraphics{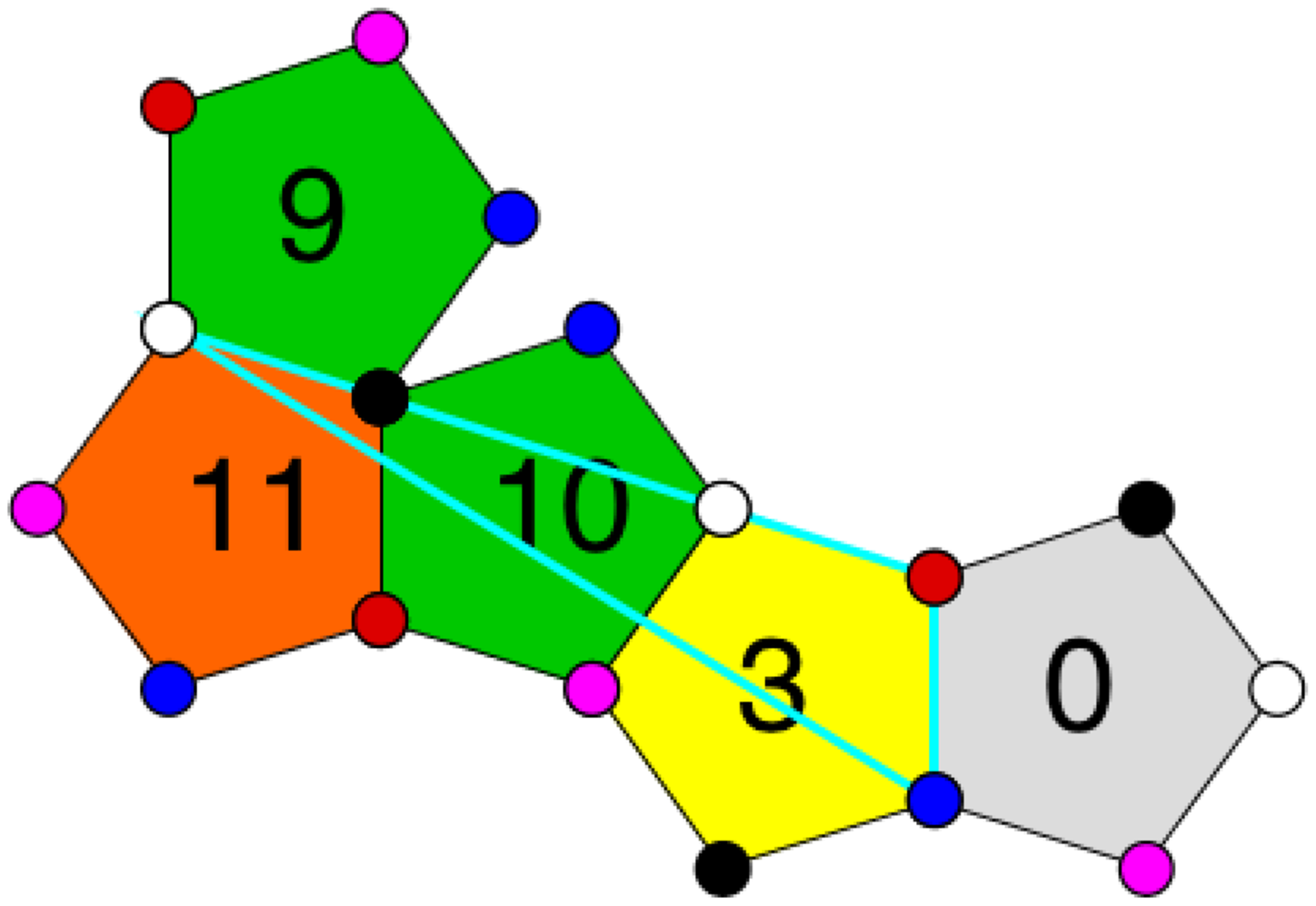}}
\newline
    {\bf Figure 4.7:\/} A compressing case that does not always work
    \end{center}

Notice that the spanning segment $\overline Q(0,0)$ is not realizable
but $\overline Q(1,0)$ and $\overline Q(0,1)$ and $\overline Q(1,1)$ are
all realizable.  A calculation shows that $\overline Q(u,v)$ is realizable
provided that $u+v \geq 1$.  This means that we can apply
our surgery whenever the intersection with the segment $\gamma=\overline{pq}$
with $Q$ is a segment of the form $Q(u,v)$ with $u+v \geq 1$.

Figure 4.8 shows another comparison chain.  In this case, the two
comparison quads are isometric, but the isometry does not preserve
the colors.  A calculation shows that
$\overline Q(u,v)$ is not longer than
$Q(u,v)$ provided that $u+v \leq 1$. Furthermore, every spanning
segment with respect to $\overline Q$ is realizable.  It is not
true that $Q \geq \overline Q$.  So, if our segment
$\gamma$ intersects $Q$ in a segment of the form $Q(u,v)$ with
$u+v \leq 1$ we can apply our surgery.  So, in all cases, we
can perform the surgery and shorten $\gamma$.

\begin{center}
\resizebox{!}{1.8in}{\includegraphics{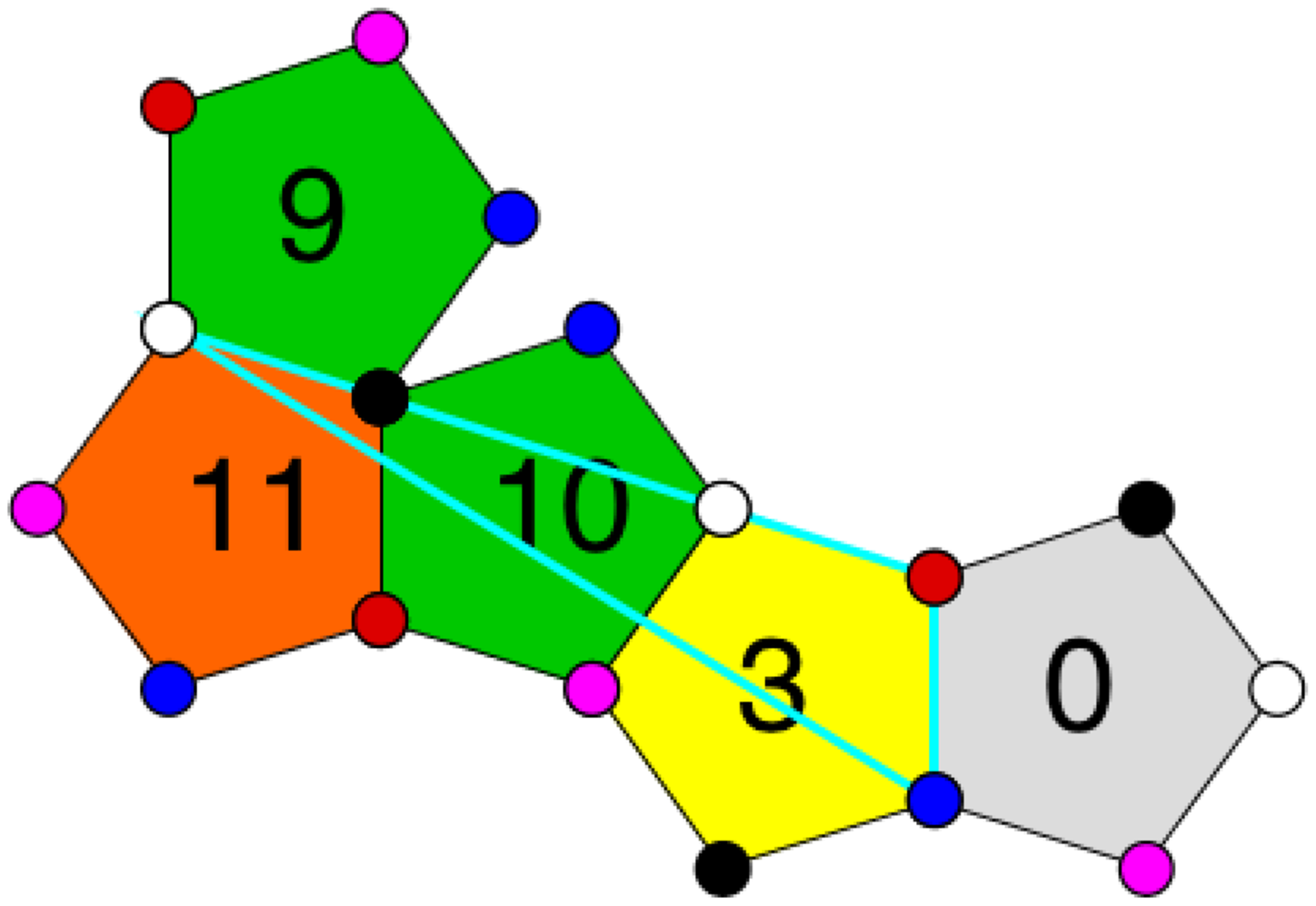}}
\newline
    {\bf Figure 4.8:\/} A partially compressing case
    \end{center}

\newpage

\section{The Comparison Lemma}

\subsection{Confining the Cities}

Our goal in this chapter is to prove the Comparison Lemma.
It is important to mention the logical structure of
our argument.  Even though, at this point in the paper,
we have not yet proved Theorem \ref{main}, we nonetheless
can study properties of the maps defined in
Theorem \ref{main}.  As we have already discussed, one
hypothesis of the Comparison Lemma is that
the description in Theorem \ref{main} really does
describe the map $p \to \widehat {\cal G\/}_p$, the
map defined in terms of our function
$\widehat d_X$.

We denote cities by their corresponding
triple points.  Thus in $C_{834}$,
the map $G$ is given by 
$G(p)=(834;p)$.  This is the point
equidistant from the vertices
$p_8,p_3,p_4$ of the decagon $D_p$.
Our ordering of the points is such
that the first two digits are common
to all $4$ cities within a state.

Now we discuss the relevant states.
Let $\Sigma_0$ be the state containing the
cities $$C_{163}, C_{168},  C_{164}, C_{169}.$$
Let $\Upsilon_0$ denote the bottom half
of $\Sigma_0$, namely the union of
$C_{163}$ and $C_{168}$.
Let $\Sigma_1$ be the state containing the
cities $$C_{831}, C_{834},  C_{836}, C_{839}.$$
Figure 5.1 shows the these cities, and also
some auxiliary line segments.

\begin{center}
\resizebox{!}{2.5in}{\includegraphics{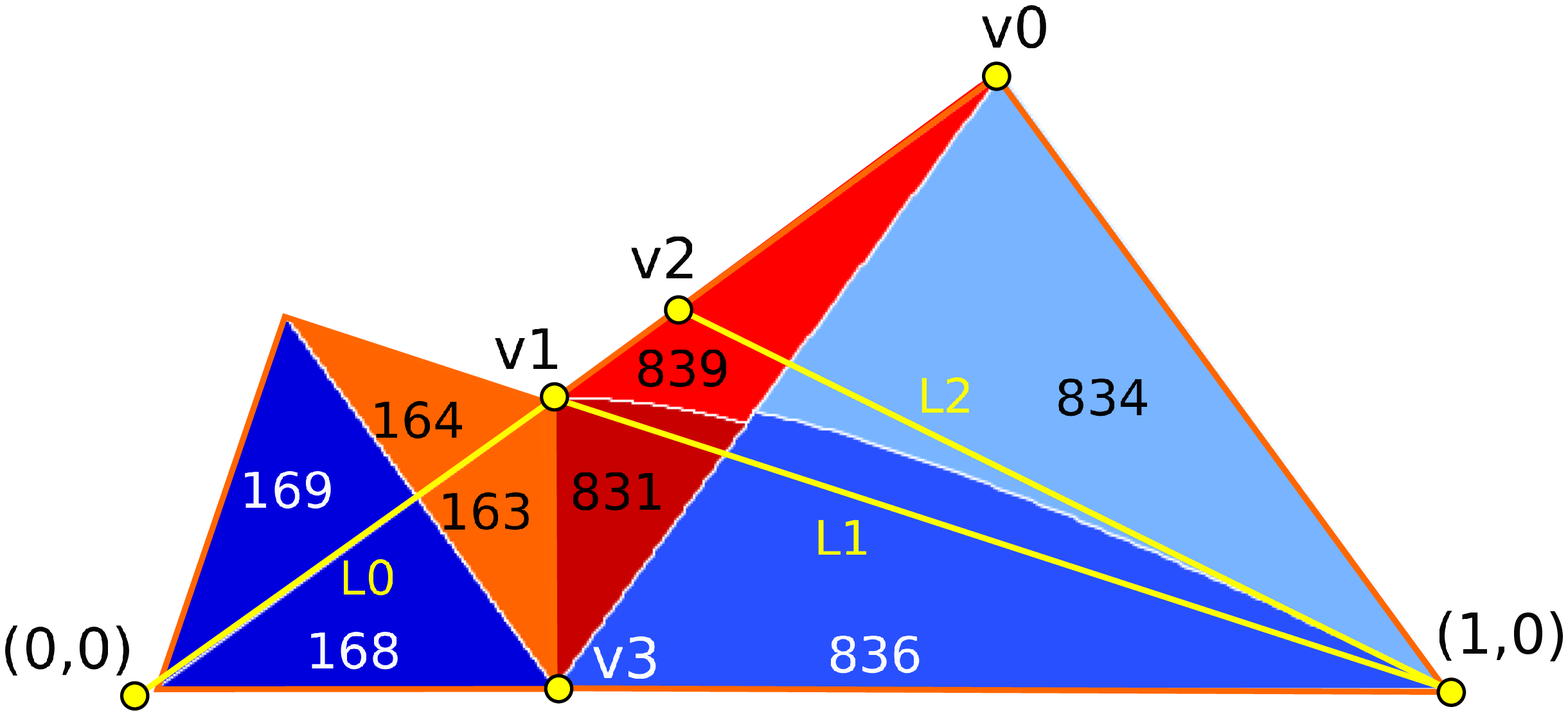}}
\newline
    {\bf Figure 5.1:\/} The cities contained in states $\Sigma_0$ and $\Sigma_1$.
\end{center}

It suffices to prove the Comparison
Lemma when $p \in \Upsilon_0 \cup \Sigma_1$.
The cities are curvilinear regions within their
states, and for the purposes of computation we
would prefer to deal with entirely polygonal
(and in fact triangular) domains.  In this section
we set this up.
Figure 5.1 also shows $3$ vertices with the following 
coordinates:
\begin{equation}
  v_0=\frac{\phi \exp(\pi i/5)}{2}, \hskip 30 pt
  v_1=\frac{\exp(\pi i/5)}{\phi^2}, \hskip 30 pt
  v_2=\frac{\exp(\pi i/5)}{2}.
\end{equation}
We define the following $3$ triangles:
\begin{itemize}
  \item $\Upsilon_{1638}$ is the right triangle with vertices $(0,0)$, $v_3$, and $v_1$.
  \item $\Upsilon_{8349}$ is the triangle whose vertices are $(1,0)$, $v_0$, and $v_1$.
  \item $\Upsilon_{8316}$ is the triangle whose vertices are $(0,0)$, $(1,0)$ and $v_2$.
\end{itemize}

By symmetry we have

\begin{lemma}[Confinement]
We have 
\begin{itemize}
\item $C_{163} \cup C_{168} \subset \Upsilon_{1638}$.
\item $C_{834} \cup C_{839} \subset \Upsilon_{8349}$.
\item $C_{831} \cup C_{836} \subset \Upsilon_{8316}$.
\end{itemize}
\end{lemma}

\startproof
The first of these results follows from symmetry.
For the other statements, define
\begin{equation}
  \label{disc0}
  f_{ijk\ell}=|(ijk,p)-p_i|^2 - |(ij\ell,p) - p_i|^2,
\end{equation}
We check that the functions $f_{8346}$ and $f_{8319}$ respectively
do not vanish on $L_1$ and $L_2$, except at the endpoints.
We then check that the signs are correct.  For instance,
$f_{8319}>0$ on $C_{831}$, and we check that
$f_{8319}>0$ on the interior of $L_1$ and
$f_{8319}<0$ on the interior of $L_2$.
Since the description in Theorem \ref{main} gives these
cities as connected sets, we see that the facts
just established prove our claims.
\endproof

We note that $C_{163},C_{168} \subset \Upsilon_{8316}$ as well,
but $\Upsilon_{1638}$ gives a tighter fit.

\subsection{The Proof}

Given the analysis in the previous chapter,
the only pentagon chains corresponding
to minimal geodesic segments joining
$p$ to some point of $A(\Pi)$ are obtained
from appending an $11$ to each of the
basic chain sequences shown in Figure 4.1 or 
else taking the dihedral image of such a
chain. 
This gives us a collection
of $70$ {\it admissible chains\/}.
Of the $70$ admissible chains, $10$ have
length $4$.  These correspond to
the straightforward geodesic segments.
We call these $10$ chains {\it straightforward\/}.
Let ${\cal S\/}$ denote the collection of $10$
straightforward chains.
Let ${\cal B\/}$ denote the collection of the $60$ other chains.

As usual we take
$p \in \Sigma_0 \cup \Sigma_1$.
For concreteness we will describe our proof for $p \in C_{834}$.  We
treat the other $5$ cities exactly the same way.
For $p \in C_{834}$ we have
\begin{equation}
\widehat d_X(p,\widehat G(p))=|(834;p)-p_8|=
|(834;p)-p_3|=
|(834;p)-p_4|.
\end{equation}
The point $(834;p)$ is equidistant from the
$3$ points $p_8,p_3,p_4$.  
The $3$ points $p_3,p_4,p_8$ are three
of the vertices of the decagon in Figure 3.1.
We give formulas in \S \ref{code}.
At the same time, there are $60$ points
$q_k=\langle C_k,p \rangle$ corresponding to
each of the $60$ chains $C_1,...,C_{60} \in {\cal B\/}$.
Again, we give formulas in \S \ref{code}.
To show that 
$$\widehat d_X((834;p),p)=d_X((834;p),p)$$
it suffices to show that
\begin{equation}
\label{makor}
|(834;p)-p_j| \leq |(834;p)-r_k| \hskip 30 pt
\forall j \in \{8,3,4\}, \hskip 20 pt
\forall k \in \{1,...,60\}.
\end{equation}

We use a trick to avoid computing the triple points.
We consider the polynomials
\begin{equation}
  \label{cr0}
P_k=\Im \frac{(p_8-p_4)(p_3-r_k)}{(p_8-p_3)(p_4-r_k)}, \hskip 30pt
k=1,...,60.
\end{equation}
Here we are taking the imaginary part of the cross ratio.
The function $P_k$ is a rational function of $(x,y)$, and
positive when $p_k$ lies outside the
disk bounded by the circle containing $p_8,p_3,p_4$.
Also $P_k=0$ if and only if the points are co-circular.
So, it suffices to show $P_k \geq 0$ on
$\Upsilon_{8346}$ for all $k=1,...,60$.
We complete the proof by doing this.
We do all the same steps for
the remaining five city-triangle pairs
$(C,\Upsilon)$ with $C \subset \Upsilon$.

To analyze the $360$ functions
of interest to us, we use the
techniques described in \S \ref{posdom}.
In \S \ref{trianglemap} we give formulas for
surjective polynomial maps from
$[0,1]^2$ to each triangle $\Upsilon$ mentioned above.
For each relevant rational function $P_k$ and each
relevant triangle map $F$ we consider
\begin{equation}
Q_k={\rm numerator\/}(P_k \circ F).
\end{equation}
We choose the numerator so that
the sum of all the coefficients of $Q_k$ is positive.
That is, $Q_k(1,1)>0$.

By construction $P_k \geq 0$ on
$\Upsilon_{8349}$ if and only if
$Q_k \geq 0$ on $[0,1]^2$.
Each polynomial $Q_k$ has the form
\begin{equation}
Q_k(x,y)=\sum_{i,j \leq 3} c_{ijk} x^i y^j, \hskip 30 pt
c_{ijk} \in \Q(\sqrt 2,\sqrt 5,\sqrt{5-\sqrt 5},\sqrt{5 + \sqrt 5}).
\end{equation}
In all $360$ cases, we check that
both polynomials $$Q_1(x,y)=Q(x/2,y), \hskip 30 pt
Q_2(x,y)=(1-x/2,y)$$ are positive dominant.
This completes the proof.

\newpage

\section{The Voronoi Decomposition}
\label{struct}

In this chapter we prove the
Voronoi Structure Lemma.  At the end
of the chapter we deduce the
Selection Lemma from the
Voronoi Structure Lemma.

\subsection{Four Triangles}

The Voronoi Structure Lemma makes two statements about the structure of
$VH_p$, the Voronoi decomposition of the hexagon $H_p$.
By symmetry it suffices to
take our point $p$ in the union
$\Upsilon_0 \cup \Sigma_1$ discussed in the last chapter.
In the previous chapter we covered
$\Sigma_1$ with two other triangles.
In this chapter we make a sharper partition.
We write $\Sigma_1 = \Upsilon_1 \cup \Upsilon_2 \cup \Upsilon_3$ where
\begin{enumerate}
\item $\Upsilon_1$ is the triangle with
  vertices $v_1$, $v_3$ and $(1,0)$.
\item $\Upsilon_2$ is the triangle with
  vertices $v_1$, $v_2$ and $(1,0)$.
\item $\Upsilon_3$ is the triangle with
  vertices $v_0$, $v_2$ and $(1,0)$.
\end{enumerate}
Figure 6.1 shows these $4$ triangles.

\begin{center}
\resizebox{!}{1.6in}{\includegraphics{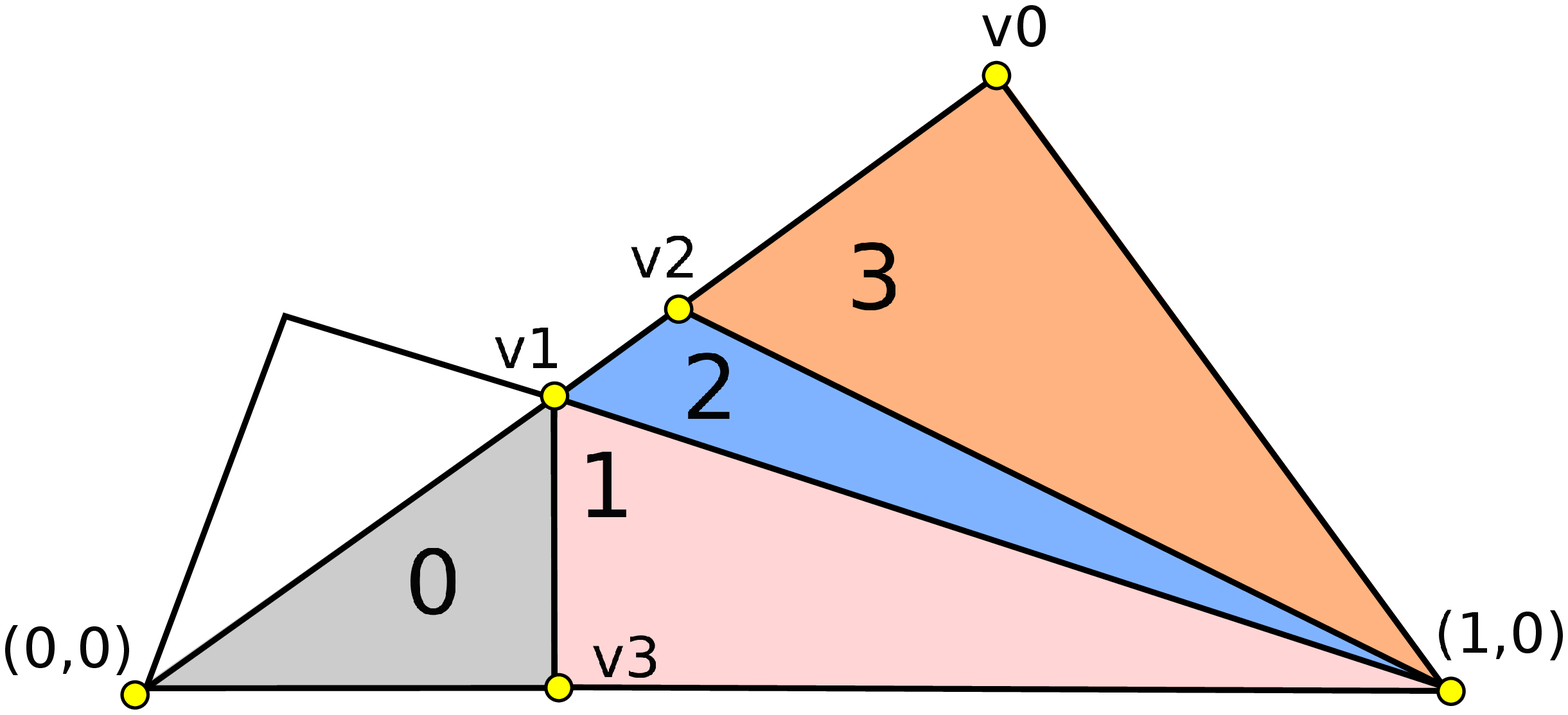}}
\newline
    {\bf Figure 6.1:\/} The $4$ triangles $\Upsilon_j$ for $j=0,1,2,3$.
\end{center}

In \S \ref{trianglemap} we describe
polynomial maps $F_j: [0,1]^2 \to \Upsilon_j$
for $j=0,1,2,3$.  These maps play an
important role in our proof, just as similar
maps played an important role in the previous chapter.
We note that $F_j$ gives a homeomorphism
between $(0,1)^2$ and the interior of $\Upsilon_j$.

\subsection{A Special Edge}
\label{edge}

The vertices of the hexagon $H_p$
are $p_1,p_3,p_4,p_6,p_8,p_9$.
We give formulas in \S \ref{code}.
A direct calculation shows that
\begin{equation}
  \label{ast}
  |p-p_j| = |p-p_{j+5} |, \hskip 30 pt j=1,3.
\end{equation}
Thus $p$ lies in the bisectors $b_{16}$ and $b_{38}$.
\newline
\newline
{\bf Remark:\/}
We called these two bisectors $\beta_1$ and $\beta_3$ in
Remark (iii) after Equation \ref{dec}.  Equation
\ref{ast}, which holds for all $j=0,1,2,3,4$,
is partly responsible for
``asterisk'' mentioned in Remark (iii).
\newline

\begin{lemma}
  If $p \in \Upsilon_0^o$ then
    $VH_p$ has an edge $e_{16}$ which
    contains $p$ in its interior and
    is contained in $b_{16}$.
\end{lemma}

\startproof
This result follows from the claim that
\begin{equation}
\label{bisector}
\Phi_j=|p-p_j|^2 - |p-p_1|^2, \hskip 30 pt j=3,4,8,9.
\end{equation}
is positive on $\Upsilon_0^o$, because then
$p$ will lie only in
the Voronoi cells $C_1$ and $C_6$ and
$e_{16}$ is the intersection of these cells.
Let $G_j=\Phi_j \circ F_0$, where $F_0$ is
our triangle map.
For each $j=3,4,8,9$ we show that the
functions
$$G_j(x/2,y), \hskip 15 pt
  G_j(1-x/2,y), \hskip 15 pt
  G_j(1/2,y)$$
  are solidly positive dominant.
  But then $G_j>0$ on $(0,1)^2$.
  But then $F_j>0$ on $\Upsilon_o^o$, as claimed.
\endproof
  
\begin{lemma}
  If $p \in \Sigma_1^o$ then
      $VH_p$ has an edge $e_{38}$ which
    contains $p$ in its interior and
    is contained in $b_{38}$.
\end{lemma}

\startproof
The proof is similar to the proof in the
previous case, but somewhat more involved.
This result folllows from the claim that
\begin{equation}
\label{bisector3}
\Phi_j=|p-p_j|^2 - |p-p_2|^2, \hskip 30 pt j=1,4,6,9.
\end{equation}
is positive on $\Sigma_1^o$.

Using the same technique as in the previous case
we establish that each $\Phi_j$ is
positive on $\Upsilon_k^o$ for $k=1,2,3$.
It remains to deal with the two lines
$L_1=\Upsilon_1 \cap \Upsilon_2$ and
$L_2=\Upsilon_2 \cap \Upsilon_3$.
Under the maps $F_1$ and $F_2$,
the points $(x,1)$ correspond to points
on $L_1$ and $L_2$ respectively.
We set $y=1$ and observe that the resulting
functions $x \to \Phi_j \circ F_k(x,1)$
are strongly positive dominant for all
relevant indices.  This shows that
$F_j$ is also positive on the relative
interior of $L_1$ and $L_2$.  Now we have
covered all points of $\Sigma_0^o$.
\endproof

The existence of $e_{16}$ places strong
restrictions on what $VH_p$ can look like
when $p \in \Upsilon_0^o$.  Figure
6.1 shows $3$ of the $9$ possible combinatorial
types for $p \in \Upsilon_0^o$.
These $3$ cases will
give a strong suggestion as to what the
other $6$ possibilities are.  Note that
we are just giving a combinatorial
representation, and not a geometric one.
The white dot represents $p$.

\begin{center}
\resizebox{!}{1.8in}{\includegraphics{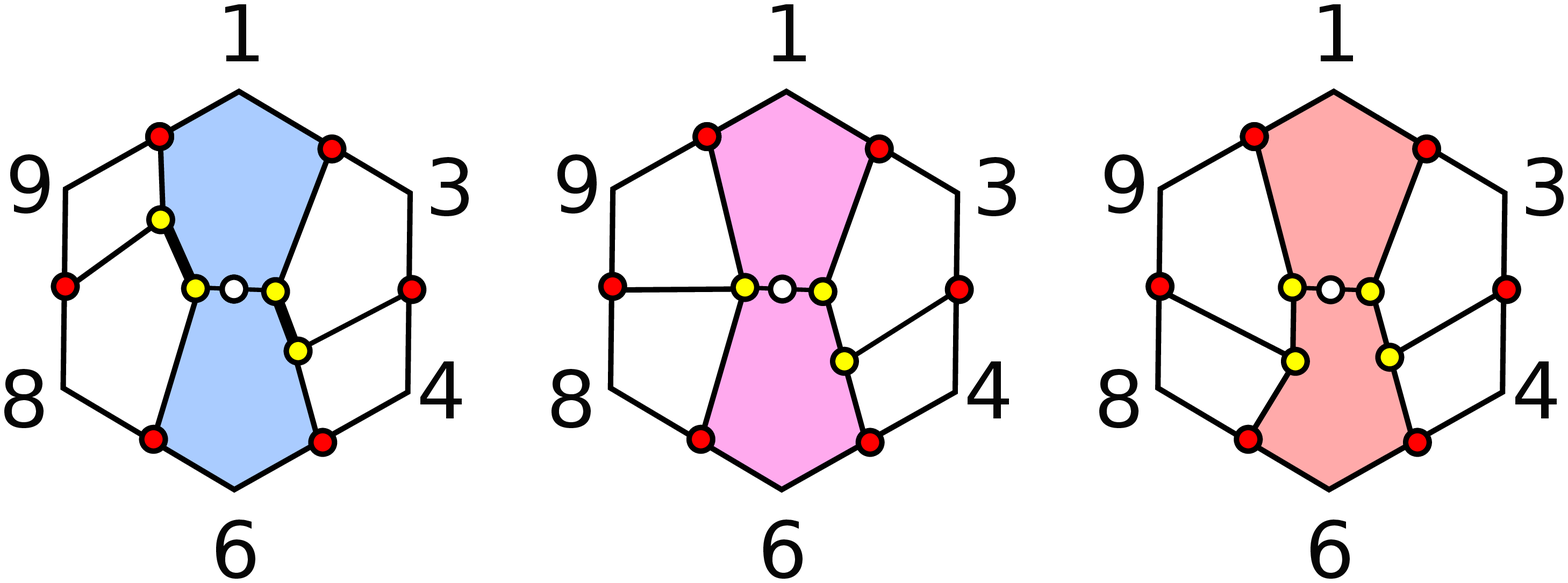}}
\newline
    {\bf Figure 6.1:\/} $3$ of $9$ combinatorial types
for $p \in \Upsilon_0$.
\end{center}

The picture for $p \in \Sigma_1^o$ is similar,
except that the indices $(1,3,4,6,8,9)$ are
replaced by the indices $(3,4,6,8,9,1)$.
That is, we just rotate the indices by one click.

Notice that the structure we have just established
places some restrictions on the vertices of
$VH_p$.  For instance, when $p \in \Upsilon_0^o$
the vertices must lie in the set:
$$
  \{(168;p),\  (189;p),\  (689;p),\  (169;p),\  (134;p),\  (146;p),\  (346;p),\  (136;p)\}.
$$
By continuity the same holds when $p \in \partial \Upsilon_0$.
One can view our result about the special edge as establishing
partial structural stability for $H_p$ in the sense of
\S \ref{voro}.  In the next section we take this further.

\subsection{Structural Stability}

Now we use our computational method to establish
some of the structural stability we mentioned
in \S \ref{voro}.  We first observe that
$VH_p$ is not structurally stable in the
small triangle $\Upsilon_2$.  This triangle
contains the interfaces between
the cities $C_{834}$ and $C_{836}$ for instance.
We have designed our partition of $\Sigma_1$ to
concentrate all the instability into this
one small triangle. 

\begin{lemma}
$VH_p$ is structurally stable for $p \in \Upsilon_0^o$.
\end{lemma}

\startproof
We consider the cases in turn.
We first check for some $p \in \Upsilon_0^o$ that
the graph for $VH_p$ is isomorphic to the
left one in Figure 6.1.  Since the
edge $e_{16}$ exists relative to all $p \in \Upsilon_0^o$
it suffices to show that the thick edges in
Figure 6.1 never shrink to points as $p$ varies
in $\Upsilon_0^o$.  These edges correspond to
the quadruples $\{1,6,8,9\}$ and $\{1,3,4,6\}$.
For each quadruple $(i,j,k,\ell)$ let
$P(i,j,k,l,x,y)$ denote the imaginary
part of the cross ratio of $p_i,p_j,p_k,p_{\ell}$
as a function of $p=x+iy$. This is just as
in Equation \ref{cr0}.   We let
$Q$ denote the numerator of $P \circ F_0$,
where $F_0$ is the triangle map to $\Upsilon_0$.
For both choices of $Q$ we check that
the two functions $Q(x,y/2)$ and $Q(x,1-y/2)$
are solidly positive dominant.  We also
check that $Q(x,1/2)$ is solidly positive
dominant.  This proves that $Q>0$ on $(0,1)^2$.
Hence $P>0$ on $\Upsilon_0^o$.

\begin{lemma}
$VH_p$ is structurally stable for $p \in \Upsilon_1^o$
and for $p \in \Upsilon_3^o$.
\end{lemma}

\startproof
The proof is similar to what we did in the
previous case.  This time (in both cases) the
two edges of interest to us correspond
to the quadruples $\{1,3,8,9\}$ and $\{3,4,6,8\}$
and the triangle maps are $F_1$ and $F_3$.
Let $Q$ be the polynomial that arises in
each of the $4$ cases, as in the previous
lemma.  We check that the $4$ functions
$$Q(x/2,y/2),\ Q(x/2,1-y/2),\ Q(1-x/2,y/2),\ Q(1-x/2,1-y/2)$$
are solidly positive dominant. This shows that
$Q>0$ on $(0,1)^2$ except perhaps on the
segment $\{1/2\} \times (0,1)$ and
$(0,1) \times \{1/2\}$. We then show that the $4$ 
single variable polynomials
$$Q(1/2,y/2),\ Q(1/2,1-y/2),\ Q(x/2,1/2),\ Q(1-x/2,1/22)$$
are solidly positive dominant, and we check explicitly
that $Q(1/2,1/2)>0$.  This shows that $Q>0$ on $(0,1)^2$.
The rest of the proof is as in the previous case.
\endproof

\subsection{Confining the Vertices}
\label{confine}

Now we prove Statement 1 of the Voronoi Structure Lemma.
In view of the structural results we have proved above,
we can say the following about the vertices of $VH_p$.
\begin{enumerate}
\item When $p \in \Upsilon_0$ the vertices are
$(163;p),\ (168;p),\ (346;p),\ (189;p).$

\item When $p \in \Upsilon_1$ the vertices are
$(831;p),\ (836;p),\ (189;p),\ (346;p).$

\item When $p \in \Upsilon_3$ the vertices are
$(834;p),\  (839;p),\ (139;p),\ (468;p).$

\item When $p \in \Upsilon_2$ the vertices are
amongst the $8$ total listed for
$\Upsilon_1$ and $\Upsilon_3$.
\end{enumerate}
Note that these vertices might not be distinct for
points in the boundaries of these various triangles.
Some of the triple points can coalesce.
In \S \ref{triple} we explain how we compute
these points in all cases.

As in \S \ref{triple} we introduce the map
$L: \C \to \R^3$ given by
\begin{equation}
L(x+i y)=(x,y,1).
\end{equation}
If $3$ complex numbers
$a,b,c$ are collinear then
\begin{equation}
\det(a,b,c):=L(a) \cdot (L(b) \times (L(c))=0.
\end{equation}

We check, for one point $p_0,p_1,p_2,p_3$ respectively in each of
$\Upsilon_0,\Upsilon_1,\Upsilon_2,\Upsilon_3$ that
all the corresponding vertices listed above lie in $\Delta$.
We just have to see that
this situation cannot change as we vary $p$ around each triangle.
Let $\ell \in \{0,1,2,3\}$ be any index.
Let $t_p$ be any triple point above associated
to $\Upsilon_{\ell}$.  Let $v,w$ be any two
vertices of the triangle $\Delta$. 
We consider the function
\begin{equation}
P(p)=\det(t_p,v,w).
\end{equation}
It suffices to prove that $\epsilon P \geq 0$ on $\Upsilon_{\ell}$
for some choice of sign $\epsilon \in \{-1,1\}$.

For this purpose we consider the function
$Q = P \circ F_{\ell}$.  We show that
the two functions $\epsilon Q(x/2,y)$ and $\epsilon Q(1-x/2,y)$ are
positive dominant for one of the two choices of sign -- the
same choice in each case.  This proves that
$\epsilon Q \geq 0$ on $[0,1]^2$ and hence
$\epsilon P \geq 0$ on $\Upsilon_{\ell}$.
There are $60=3 \times 4 + 3 \times 4 + 3\times 8 + 3\times 4$
functions in total, and we make the check in each case.
This completes the proof.
\newline
\newline
{\bf Remark:\/}
We might have taken more effort, as in the previous
chapter, to pick our signs in advance so that
$\epsilon=1$ in all cases, but we do not really
need to bother with this. Our function checks
that either $Q$ is positive dominant or
$Q$ is negative dominant by
taking the maxima and minima of all the coefficient
sums that arise in the definition, and then it
checks that the max and the min do not have
opposite signs.  This suffices.

\subsection{Comparing the Distances}

In this section we prove Statement 2 of the
Voronoi Structure Lemma.  Our proof here is
almost exactly the same as what we did for the
Comparison Lemma.
First of all, we check for some choice of $p$
in each triangle $\Upsilon_0,\Upsilon_1,\Upsilon_2,\Upsilon_3$ that
\begin{equation}
\label{compare}
|(ijk;p)-p_i|<|(ijk;p)-p_{\ell}|, \hskip 30 pt \ell \in \{0,2,5,7\}.
\end{equation}
Here $p_{\ell}$ is one of the vertices of the decagon $D_p$ which
is not a vertex of the hexagon $H_p$ and $(ijk;p)$ is one of the
triple points above associated to the relevant triangle.
We mean to say that we check all possibilities for all triangles.
This makes for $80=4 \times 4 + 4 \times 4 + 8 \times 4 + 4 \times 4$ checks.

Equation \ref{compare} is
perhaps stronger than Statement 2 of the Voronoi Structure Theorem
because perhaps some of the distances on the right hand side do
not correspond to geodesic segments in the dodecahedron connecting
the $2$ relevant points.  Also,
the point $(ijk;p)$ might not actually be a vertex of $VH_p$.
None of this bothers us. We will establish
Equation \ref{compare} for all relevant indices and for all
points in the interiors of our triangles.  By continuity,
we still have a weak inequality even for boundary points.
Hence, Equation \ref{compare} holds for all $p$ in the relevant
triangle provided that $(<)$ is replaced by $(\leq)$.  
This still implies Statement 2 for all points in the triangle.

Here are the calculation details.
As $p$ varies around one of the triangles, Equation \ref{compare} fails
only if the $4$ points $(p_i,p_j,p_k,p_{\ell})$ become cocircular.
As in the previous chapter, we let 
$P$ be the imaginary part of the cross ratio of these points and we
let $Q$ be the numerator of $P \circ F$, where $F$ is the relevant
triangle map. For each of the $80$ choices of $Q$ we check that
the functions $\epsilon Q(x/2,y)$ and $\epsilon Q(1-x/2,y)$ and $\epsilon Q(1/2,y)$ are
solidly positive dominant for one of the two choices
$\epsilon \in \{-1,1\}$. This proves that
$\epsilon Q>0$ on $(0,1)^2$.   Since $Q$ never vanishes on $(0,1)^2$, we
never lose the inequality in Equation \ref{compare}.
This completes the proof.

\subsection{Proof of the Selection Lemma}
\label{selection}

In this section we deduce the Selection Lemma from the
Voronoi Structure Lemma.
Our argument
refers to Figure 6.2, which shows $VH_p$
for some randomly chosen $p \in \Delta$
superimposed over the plan $P$ from Figure 3.1.
The reader can see much better pictures of this
using my program.
The right side of Figure 6.2 is a close-up of
the left side.

\begin{center}
\resizebox{!}{2.6in}{\includegraphics{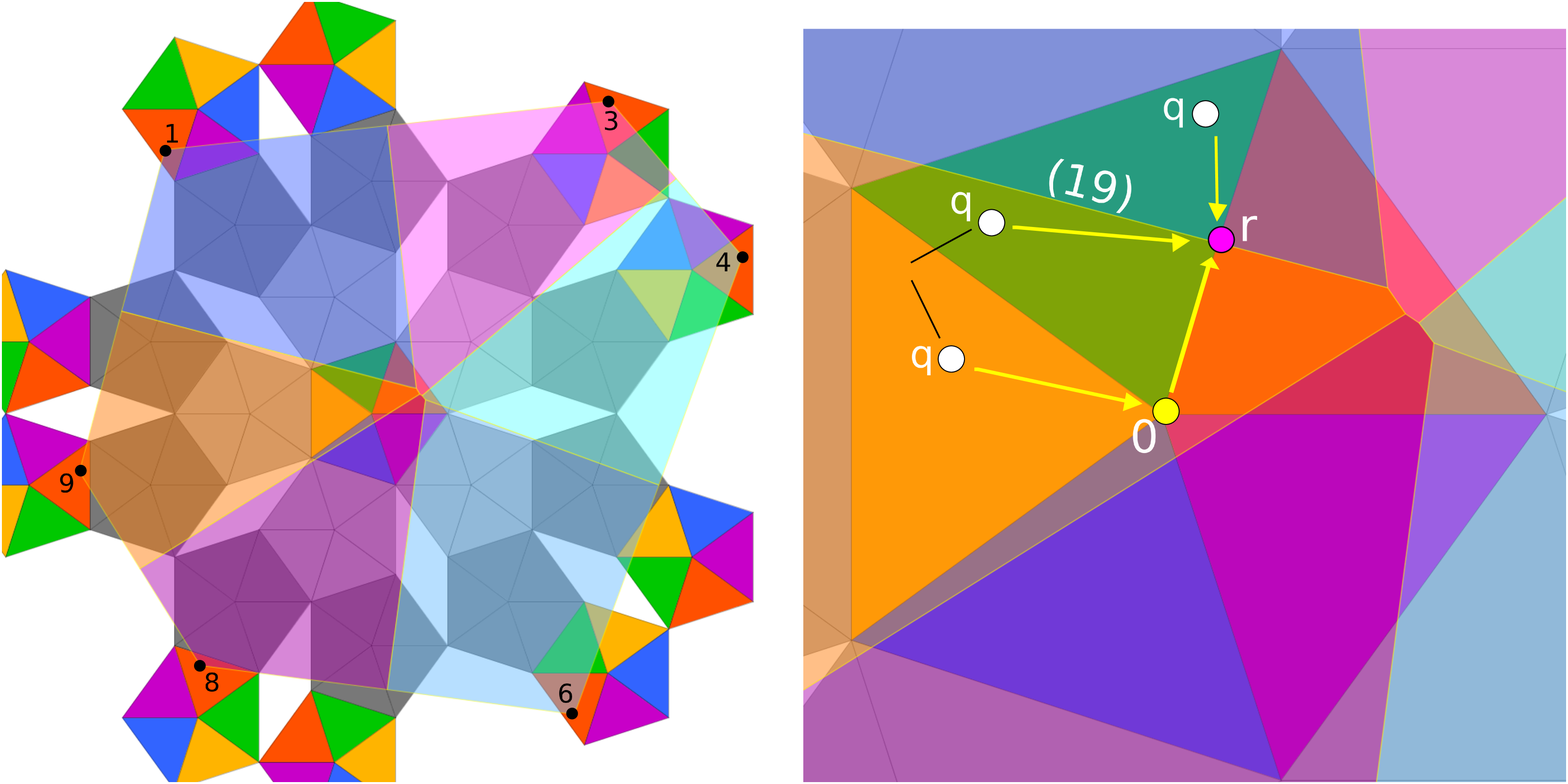}}
\newline
    {\bf Figure 6.2:\/} The set $VH_p$ superimposed over the dodecahedral plan $P$.
\end{center}

We parametrize the bisector $b_{19}$ so that we start at
$\partial H_p$, to the left of $\Delta$, and move
rightwards.  By the Voronoi Structure
Lemma, $b_{19}$ eventually intersects $\Delta$,
and therefore first hits $\Delta$ somewhere in
the left edge.  This intersection point is $r$ in Figure 6.2.

Let $\omega^k \Delta$ be the triangle obtained
from $\Delta$ by multiplying to $\omega^k$.
Here $\omega=\exp(2 \pi/5)$.  Let
$q \in \widehat {\cal G\/}_p$.  We want to
rule out the possibility that
$q \in \omega^k\Delta-\Delta$ for $k=1,2,3,4$.
By symmetry, it suffices to consider
$\omega \Delta$ and $\omega^2 \Delta$.
These are respectively the green and orange central triangles
in Figures 3.1 and 6.2.
Let $\Delta_j$ denote the red triangle in Figure 6.2
containing the point $p_j$.  (This point is just
denoted by a $j$ in the figure.)
\newline
\newline
    {\bf Case 1:\/}
Every point in $\omega \Delta$ can be joined to
every point of $\Delta_1$ and to every point
of $\Delta_9$ by line segments which remain inside
the dodecahedron plan $P$.  Such segments therefore
correspond to straightforward geodesic segments
joining $q$ to $A(p)$ in the dodecahedron.
Thus, if $q \in \omega \Delta$,
we have
\begin{equation}
  \label{eqright}
  \widehat d_X(p,A(q)) = \widehat d_X(q,A(p)) \leq \min(|q-p_1|,|q-p_9|)=\phi_{19}(q)
\end{equation}
The first equality just comes from antipodal symmetry.
The last equality is our way of defining the function $\phi_{19}$.

Since $q \in \omega \Delta-\Delta$ we have
$\phi_{19}(q)<\phi_{19}(r)$.
We have the inequalities
\begin{equation}
  \label{chain0}
  \widehat d_X(p,A(q)) \leq  \phi_{19}(q)<\phi_{19}(r) =
\mu_{H_p}(r) = \widehat d_X(p,A(r)).
\end{equation}
This contradicts that $q \in \widehat {\cal F\/}_p$.
The first equality comes from the fact,
thanks to the Voronoi Structure Theorem, that
$b_{19}$ encounters no triple points until
it reaches $\Delta$. Thus $r \in C_1 \cup C_9$,
the union of these two Voronoi cells.

The last equality in Equation \ref{chain0} needs some more justification.
Since $r \in H_p$, every line segment joining
$r$ to a vertex of $H_p$ lies in the plan $P$
and thus is the developing image of a
straightforward geodesic segment which joins $r$ to $A(p)$.
Let $\gamma$ be the
shortest of these and
let $\ell(\gamma)$ denote its length.
By construction
$$\mu_{H_p}(r)=\ell(\gamma).$$
Now, any other straightforward geodesic
segment joining $r$ to $A(p)$ corresponds
to some line segment connecting $r$ to
one of the vertices of the decagon $Q_p$.
But Statement 2 of the Voronoi Structure
Theorem says that such segments are no
shorter than $\gamma$.  Hence
$\gamma$ is a shortest straightforward
geodesic connecting $r$ to $A(p)$.
Hence $$\widehat d_X(r,A(p))=\ell(\gamma)$$ as well.
By symmetry,
the same is true for $d_X(p,A(r))$.
\newline
\newline
{\bf Case 2:\/}
Now we treat the case $k=2$.
Every point of $\omega^2 \Delta$ can
be joined to every point of
$\Delta_9$ by a line segment
that remains inside the plan $P$.
Hence
\begin{equation}
  \widehat d_X(p,A(q)) \leq |q-p_9|=\phi_9(q).
\end{equation}
Since $q \not \in \Delta$ we have
$\phi_9(q)<\phi_9(0)$.   Here $0$ is
the yellow dot on the right side
of Figure 6.2.   If we drop a
perpendiclar from $p_9$, an
arbitrary point of $\Delta_9$, to
the line extending the left edge of
$\Delta$, the intersection point
lies below $\Delta$.  Hence
$\phi_9(0) \leq \phi_9(r)$.
But now we have
\begin{equation}
  \widehat d_X(p,A(q)) \leq  \phi_{9}(q)<\phi_{9}(0) \leq \phi_{9}(r)=^*
  \phi_{19}(r)=
\mu_{H_p}(r) = \widehat d_X(p,A(r)).
\end{equation}
The starred equality comes from the fact
that $r \in b_{19}$.  The last two inequalities
are the same as in Equation \ref{chain0}.
We get the same contradiction as in Case 1.

\newpage

\section{The Vertex Competition}
\label{struct2}

\subsection{The Competition Lemma}

In this section we prove the Competition Lemma.
Let us take stock of what we know so far.
As in the previous chapters we take
$$p \in \Upsilon_0 \cup \Upsilon_1 \cup  \Upsilon_2 \cup \Upsilon_3=
\Upsilon_0 \cup \Sigma_1.$$
What we know so far is that $\widehat {\cal G\/}_p$ is always an
essential vertex of the Voronoi decomposition $VH_p$.
The possibilities are listed in \S \ref{confine}.
Also, for each of these points $q=(ijk;p)$ we know that
the distance from $q$ to $A(p)$ is given by
$$\mu_{H_p}(q)=\min_{j \in \{1,3,4,6,8,9\}} |q-p_j|.$$
We are done with the decagon, so we set
$\mu_p=\mu_{H_p}$.

\begin{lemma}
  \label{state0}
  Let
  $p \in \Upsilon_0$.
  \begin{enumerate}
    \item The point $(189;p)$ belongs to
      $\widehat {\cal G\/}_p$ only if   $(189;p)=(168;p)$.
    \item The point $(346;p)$ belongs to
      $\widehat {\cal G\/}_p$ only if $(346;p)=(163;p)$.
  \end{enumerate}
\end{lemma}

\startproof
Consider the first statement.
We show
$\mu_p((168;p))> \mu_p((189;p))$
whenever these points are distinct.
Since $1$ is an index for both points, this is the same as showing that
$|p_1-(168;p)|>|p_1-(189;p)|$ whenever the two points are distinct.
The two points $(168;p)$ and $(189;p)$
both lie along the bisector $b_{18}$ associated to $p$.
Looking at Figure 3.1 or Figure 6.2 we see that
the line $\overline{p_1p_8}$ lies entirely to the
left of the central pentagon $\Pi$ and in particular
entirely to the left of $\Upsilon_0$.  As we travel
rightward along $b_{18}$ away from $\overline{p_1p_8}$ we increase
the distance to $p_1$ (and to $p_8$).  Our rightward travel
brings us first to $(189;p)$ and then to $(168;p)$.
Hence the latter point is farther from $p_1$.

Exactly the same argument, with $b_{36}$ replacing
$p_{18}$, and {\it right\/} replacing
{\it left\/}, establishes the second statement.
\endproof

\begin{lemma}
  \label{state1}
  Let
  $p \in \Sigma_1$.
  \begin{enumerate}
    \item The point $(139;p)$ belongs to
      $\widehat {\cal G\/}_p$ only if   $(139;p)=(839;p)$.
    \item The point $(468;p)$ belongs to
      $\widehat {\cal G\/}_p$ only if $(468;p)=(834;p)$.
    \item The point $(189;p)$ belongs to
      $\widehat {\cal G\/}_p$ only if   $(189;p)=(831;p)$.
    \item The point $(364;p)$ belongs to
      $\widehat {\cal G\/}_p$ only if $(364;p)=(836;p)$.
  \end{enumerate}
\end{lemma}

\startproof
This has the same kind of proof as Lemma \ref{state0}.
The key points are that lines
$\overline{p_3p_9}, \overline{p_4p_8}, \overline{p_1p_8},
\overline{p_3p_6}$ respectively lie above,
below, left of, and right of, $\Sigma_1$.
\endproof

We use these lemmas to prove the Competition Lemma.
From the Vertex Lemma, we know that when $p \in \Upsilon_0$, the set
$\widehat {\cal G\/}_p$ is contained in
the union of $4$ triple points listed in
\S \ref{confine} in connection with
$\Upsilon_0$.  Lemma \ref{state0} eliminates
two of these, leaving the two listed in
Statement 2 of the Competition Lemma.
Statement 3 follows from Statement 2 and from
symmetry.

Statement 4 follows from
Lemma \ref{state1} in the
same way that Statement 2 follows
from Lemma \ref{state0}.

Now we turn to Statement 1.
Our analysis above eliminates
all vertices of the Voronoi
decomposition $VH_p$ except
those incident to the edge
$e_{16}$ (in $\Sigma_0$) and
the edge $e_{38}$ (in $\Sigma_1$)
which we analyzed in \S \ref{edge}.
A direct calculation shows that when
we approach the $\partial \Sigma_0$
the edge $e_{16}$ shrinks to a point
and when we approach $\partial \Sigma_1$
the edge $e_{38}$ shrinks to a point.
But this combines with Statements 1-3
to show that $\widehat {\cal G\/}_p=\{p\}$ in
all cases.

Each of the triple points listed in the
Competition Lemma does actually arise as
a member of $\widehat {\cal G\/}_p$ for
suitable choices of $p$.  So, it only
remains to analyze how the placement
of $p$ inside the two states
$\Sigma_0$ and $\Sigma_1$ determines
the triple with the largest $\mu_p$-value.
After a preliminary section which makes
some definitions we need, we  treat the
two states in turn.

\subsection{Preliminaries}

To help us analyze the dependence
$\widehat {\cal G\/}_p$ on the point $p$,
we consider the function from Equation \ref{disc0}.
Here it is again.

\begin{equation}
  \label{disc}
  f_{ijk\ell}=|(ijk,p)-p_i|^2 - |(ij\ell,p) - p_i|^2,
\end{equation}

Suppose we are in some open set $U$ where
$f_{ijkl}>0$ throughout $U$.   This means that
$\widehat {\cal G\/}_p$ cannot contain
$(ij\ell;p)$ when $p \in U$.  If
$f_{ijkl}(p)=0$ for some $p \in U$ it means that
$\widehat {\cal G\/}_p$ either contains both
$(ijk;p)$ and $(ij\ell;p)$ or neither of them.
Finally, if $f_{ijkl}<0$ throughout $U$ it
means that $\widehat {\cal G\/}_p$ cannot contain
$(ijk;p)$ when $p \in U$.

Given Statement 1 of
the Competition Lemma, we already know what happens
on $\partial \Sigma_0$ and $\partial \Sigma_1$.
So, we will work with the open
regions $\Upsilon_0^o$ and $\Sigma_0^o$.
Restricting our attention to
$\Upsilon_0^o$ ignores one part of
the space, namely the hypotenuse
of $\Upsilon_0$.  For ease of
exposition, we do not specially
treat this edge.  The behavior of
$\widehat {\cal G\/}$ on this edge
is just the continuous extension of
the behavior in $\Upsilon_0^o$.

For the indices of interest to us, the
function $f_{ijk\ell}$ always factors into
smaller factors, one of which is a
cubic function $g_{ijk\ell}$.  The
cubic determines the interface between the
cities $C_{ijk}$ and $C_{ij\ell}$.
In all but one case, the other factors
are linear, and we know them from geometric
reasoning:  Our function
$f_{ijk\ell}$ vanishes on two of the edges on
the boundary of the domain of interest, and
the defining functions for the lines extending
these edges are the linear factors.
In one case $h_{ijk\ell}$ factors into two
cubics and the other cubic $h_{ijk\ell}$ turns
out to be nonzero on the domain.  I do not understand
this other cubic geometrically.

We also mention a shortcut. Sometimes
we do not need to compute
$f$ but rather can get $g$ by more
direct means.  We consider one
of the two examples where this applies.
Given the structure of $VH_p$, the points
$(831;p)$ and $(839;p)$ are not both
essential vertices unless they are
equal.  So, instead of using the
function $f_{8319}$ above we can
directly compute the function
$g_{8319}$ which is the imaginary
part of the cross ratio of
$(p_8,p_3,p_1,p_9)$.
Similar remarks apply for the
quadriple $(8,3,4,6)$.

The analysis in $\Sigma_0$, which just involves
$g_{1638}$, is pretty straightforward.
We concentrate here on $\Sigma_1$.
The analysis in $\Sigma_1$ is harder because
we do not have the same bilateral symmetry.
We say that a smooth function $g$
is {\it cleanly\/} related to $\Sigma_1$ if
\begin{enumerate}
\item Some directional derivative of $g$ is nonzero throughout $\Sigma_1$.
\item $g$ vanishes on exactly two points of $\partial \Sigma_1$.
\end{enumerate}
If $g$ is cleanly related to $\Sigma_1$ then the
$0$-set $\gamma$ of $g$ intersects $\Sigma_1$ in a single
smooth arc which has two points on the boundary.

For Property 1, we will try to use the direction
parallel to the segment foliation in $\Sigma_1$.
These segments are parallel to the line $y=-\tan(\pi/5) x$.
We call this the {\it preferred direction\/}.
If we can't get the preferred direction to work, we will
use the direction $(1,1)$.

We use one of two methods to check Property 2.
One approach, the {\it calculus approach\/}, is
to check that the directional derivatives of $g$ along the directions
parallel to the sides of $\Sigma_1$ do not vanish
in $[0,1]^2$.  In these cases, $g$ vanishes at
most once in each edge of $\partial \Sigma$, but
the level curve can only connect up two
of these points.  Hence $g$ vanishes only
twice on $\partial \Sigma_1$.
The other approach,
the {\it restriction method\/} is that we
check explicitly that
that $g$ is nonzero on $\partial \Sigma$
except at a pair of vertices
of $\partial \Sigma_0$.  In this case,
$\gamma$ connects two vertices of $\Sigma_1$.

\subsection{The Zeroth State}
\label{zeroth}

We consider the picture in $\Upsilon_0^o$.
In light of the Competition Lemma, we have the
following implications.
\begin{enumerate}
  \item $f_{1638}(p)>0$ implies that $\widehat G(p)=(163,p)$.
  \item $f_{1638}(p)<0$ implies that $\widehat G(p)=(168,p)$.
  \item $f_{1638}=0$ implies that $\widehat {\cal G\/}_p \subset \{(163,p),(168,p)\}$.
\end{enumerate}

We mention that the vertical boundary of $\Sigma_0$ lies in the solution to the equation
$x=\cos(2 \pi/5)$.  The segment foliation in $\Sigma_0$ (and in
$\Upsilon_0$ by restriction) is parallel to
the line $L$ whose equation is $y=x \tan(\pi/5)$.
Using the formulas from
\S \ref{formula} we compute in Mathematica that
\begin{equation}
  f_{1638}(x,y)=C \times (x-\cos(2 \pi/5)) \times y \times g_{1638},
\end{equation}
where $C$ is some constant of no interest to us.
In other words, $f_{1368}$ vanishes along the horizontal
and vertical sides of $\Upsilon_0$.  The nontrivial factor,
$g_{1638}$, is a cubic function which is positive on
the interior of the vertical edge of $\Upsilon_0$,
negative on the horizontal edge of $\Upsilon_0$.
We list the formula for $g_{1638}$ in
\S \ref{formula}.

The directional derivative of $g_{1638}$ in the
direction of $L$ never vanishes in $\Sigma_0$.
Numerically, we have
$$\psi=\partial_L(g_{1638})=\nabla g_{1638} \cdot (\cos(\pi/5),\sin(\pi/5)) \approx $$
$$-38.5375 - 4.10995 x + 6.65003 x^2 - 2.16073 y + 6.99226 x y + 
 6.65003 y^2.$$
One can see the non-vanishing from these numerics: $\Sigma_0 \subset [0,1]^2$ and
the constant term is much more negative
than any of the other terms.

From the structure just discussed,
we conclude that the algebraic curve 
$$\gamma_{1368}=\{g_{1638}=0\} \cap \Upsilon_0$$
separates $\Upsilon_0$ into the
cities $C_{163}$ and $C_{168}$ in
 the combinatorial pattern shown in
Figure 1.4, such that each segment in the
foliation of $\Sigma$ (which is parallel to $L$)
intersects the curve exactly once.
This gives the claimed decomposition of
$\Upsilon_0$ into cities.
When $p \in \gamma_{1638}$, we have the
third option listed above, and our two maps
are inverses of each other.

We are done with $\Upsilon_0$. Now let us consider
the picture in the bigger domain $\Sigma_0$.
Let $\rho$ denote the reflection in $L$.
By symmetry, we have similar results for
$p \in \rho(\Upsilon_0^o)$, with the indices
$4,9$ replacing the indices $3,8$.
Also by symmetry, we have
$(163,p)=(164,p)$ and $(168,p)=(169,p)$ for $p \in L$.
This gives us the desired result even when
$p \in L \cap \Sigma_0^o$.
Finally, by symmetry, the combinatorial pattern
of the cities in $\rho(\Upsilon_0)$ is just obtained from
the pattern in $\Upsilon_0$ and reflecting it
across $L$.
This is as in Figure 1.4.
Hence Theorem \ref{main} describes $\widehat {\cal G\/}$ for all
points in the state $\Sigma_0$.

\subsection{The First State}
\label{firststate}

We will give the analysis of the
picture in $\Sigma_1$ modulo some
technical lemmas which we prove in the
next section.
Let $L$ be the line $y=-\tan(\pi/5)x$. This is the
direction parallel to the line segment foliation of $\Sigma_1$.
Let $\cal L$ denote the segment foliation.

Reflection in the real axis has the following
effect on the cities
$C_{168} \leftrightarrow C_{831}$ and
$C_{163} \leftrightarrow C_{836}$.
For this reason,
\begin{equation}
f_{8316}(x,y)=f_{1638}(x,-y), \hskip 40 pt
g_{8316}(x,y)=f_{1638}(x,-y).
\end{equation}
What we mean is to say that $f_{8316}$ factors
just as $f_{1638}$ does, and the cubic
factor $g_{8316}$ satisfies the identity just given.
We use the cross ratio method to compute
$g_{8319}$, as we already discussed.

\begin{lemma}
  \label{intersect}
  The functions $g_{8316}$ and $g_{8319}$ are cleanly related to
  $\Sigma_1$ and their gradients are linearly independent at each point of
  $\Sigma_1$.
\end{lemma}

\startproof
See \S \ref{tech}.
\endproof

The level curve $\gamma_{8319}$
intersects $\partial \Sigma_0$ at the $(1,0)$
and at the vertex opposite $(1,0)$.  In particular,
the endpoints of
$\gamma_{8316}$ and $\gamma_{8319}$ are interlated
on $\partial \Sigma_1$.
This means that
$\gamma_{8316}$ and
$\gamma_{8319}$ intersect at least once.
Given the properties in the previous
lemma, these two curves intersect exactly once
in $\Sigma_1$.

Now we can say that
the subset where $g_{8316}>0$ and $g_{8319}>0$ is
a simply connected domain
bounded by $1$ line segment and $2$ cubic curves.  This
counts as a city.  We call this city
$C_{831}$.  We check that the sign of
$g_{8319}$ is correct in the sense that
$g_{8319}>0$ implies that
$(831;p)$ rather than
$(839;p)$ is the essential vertex of
$VH_p$.
From all this, we conclude that
\begin{itemize}
\item $G(p)=(831,p)$ for $p \in C_{831}^o$.
\item $\widehat {\cal G\/}_p$ does not contain
  $(831,p)$ when $p \in \Sigma_0-C_{831}$.
\item $(831,p)=(839,p)$ along $\gamma_{8319}$.
\end{itemize}
Figure 7.1 shows $C_{831}$ in red.

\begin{center}
\resizebox{!}{1.6in}{\includegraphics{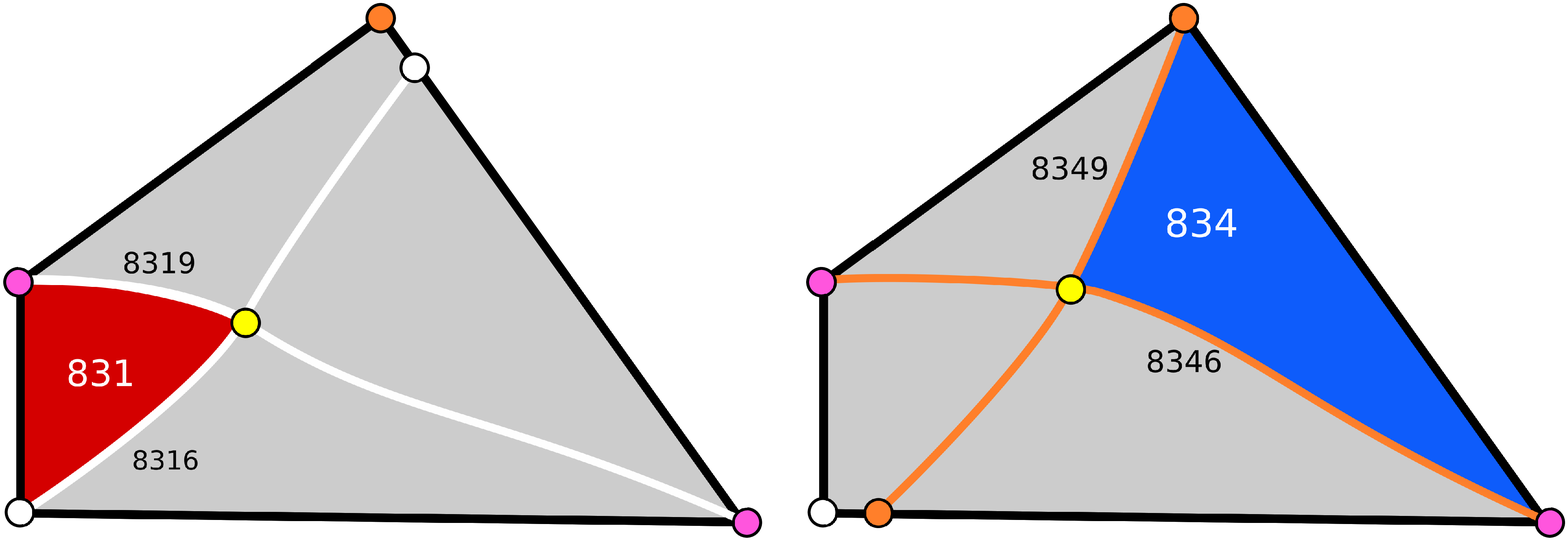}}
\newline
    {\bf Figure 7.1:\/} The cities $C_{831}$ and $C_{834}$ cut out by
        intersecting curves.
\end{center}

We define
$g_{8349}$ and $g_{8346}$ in the same way, respectively,
just as we defined
$g_{8316}$ and $g_{8319}$.

\begin{lemma}
  \label{intersect2}
  The functions $g_{8349}$ and $g_{8346}$ are cleanly related to
  $\Sigma_1$ and their gradients are linearly independent at each point of
  $\Sigma_1$.
\end{lemma}

\startproof
This has the same proof as Lemma \ref{intersect}.
\endproof

We define $C_{834}$ just as we defined
$C_{831}$.  From the structure above we conclude the
following:
\begin{itemize}
\item $G(p)=(834,p)$ for $p \in C_{834}^o$.
\item $\widehat {\cal G\/}_p$ does not contain
  $(834,p)$ when $p \in \Sigma_0-C_{834}$.
\item $(834,p)=(836,p)$ along $\gamma_{8346}$.
\end{itemize}
Figure 7.1 draws $C_{834}$ in blue.  We have
drawn these cities as disjoint.  We will justify
the disjointness at the end.

We have $\widehat {\cal G\/}_p \subset \{(836,p),(839,p)\}$ for
$p \in \Sigma_0^o-C_{831}-C_{834}$.   The function
$f_{8369}$ is the most complicated one we need to consider.
After a lot of {\it ad hoc\/} work on Mathematica,
I found that $f=gh$ where
$g_{8369}$ and $h_{8369}$ are both cubics.
The formula for $h_{8369}$ is given in
Equation \ref{h8369}, and the reader can
verify that indeed $f=gh$.  I don't have
a geometric understanding of $h_{8369}$.

\begin{lemma}
  $h_{8369}$ is positive on $\Sigma_1^o$ and
  cleanly related to $\Sigma_1$.
\end{lemma}

\startproof
See \S \ref{tech},
  \endproof
  
  Thus, $g_{8369}$ determines the competition between
  $(836;p)$ and $(839;p)$ in $\Sigma_1$.
  Let $\gamma_{8369}$ be the zero level curve in
  $\Sigma_1$ associated to $g_{8369}$.

  \begin{lemma}
    $\gamma_{8369}$ intersects each of
    $\gamma_{8316}, \gamma_{8319}, \gamma_{8346},
    \gamma_{8349}$ exactly once.
  \end{lemma}

  \startproof
  See \S \ref{tech}.
  \endproof

  Now let us analyze the locations of the intersection points.
The intersection point
$\gamma_{8369} \cap \gamma_{8319}$ must lie
on the curve $\gamma_{8319}$ because,
  as we have already mentioned, we have
  $(831,p)=(836,p)$ at this point.
  Hence $\gamma_{8369}$ intersects
  $C_{831}$ at its interior vertex.
  Likewise  $\gamma_{8369}$ intersects
  $C_{834}$ at its interior vertex.
  
  Consider the arc
  \begin{equation}
    \beta=\gamma_{8316} \cup \gamma_{8369} \cup \gamma_{8319}.
    \end{equation}
  This is the red-blue interface on the right side of
  Figure 1.4.

  \begin{lemma}
    The arc $\beta$ is embedded and joins the origin
    to the opposite vertex of $\Sigma_1$. Moreover, each segment of
    $\cal L$ meets $\beta$ exactly once.
  \end{lemma}

  \startproof
   Each of the three arcs of $\beta$ is embedded,
  and the middle arc intersects the other two
  exactly once.  To show that $\beta$ is embedded,
  it suffices to prove   that $\gamma_{8316}$ and $\gamma_{8319}$ are
  disjoint.  Recall that each segment of
  $\cal L$. meets each of these arcs at most once.
  We just have to show that no segment of
  $\cal L$ meets both arcs.

  The segments of $\cal L$ have a
  {\it transverse order\/}:  We order them
  according to the $y$-intercepts of the
  lines extending these segments.
    The two vertices of $\beta$ are
  $b_1=\gamma_{8316} \cap \gamma_{8369}$ and
  $b_2=\gamma_{3869} \cap \gamma_{8319}$.  
  As we move along $\gamma_{8369}$ away from
  the origin and transverse to the segments in $\cal L$,
   we encounter $b_1$ before $b_2$.
  Hence, all the segments of $\cal L$ which
  meet $\gamma_{8316}$ precede, in the
  transverse order, all the
  segments of $\cal L$ which meet
  $\gamma_{8319}$.  This establishes
  the claim that $\beta$ is embedded.

  We check that $g_{8316}$ vanishes at the origin and
  locally takes opposite signs on the edges of
  $\Sigma_1$ incident to the origin. Hence
  $\gamma_{8316}$ has one endpoint at the origin.
  Call the origin $v_0$.
  Similarly, $\gamma_{8319}$ has one endpoint
  at the vertex $v_1$ of $\Sigma_1$ opposite the origin.
  Hence $\beta$ is an embedded arc joining
  $v_0$ to $v_1$.
  
  From what we have said above, each segment of
  $\cal L$ meets $\beta$ at most once.  Since
  $\beta$ joins a pair of opposite vertices of
  $\Sigma_1$ and since every segment of $\cal L$
  meets both components of $\partial \Sigma_1-\{v_0,v_1\}$,
  we see that every segment of $\cal L$ meets
  $\beta$.  Hence every segment of $\cal L$
  meets $\beta$ exactly once.
  \endproof

  \begin{corollary}
    $C_{831}$ and $C_{834}$ are disjoint.
  \end{corollary}

  \startproof
  Each of these cities shares an arc with
  $\beta$ and hence lies in one component
  of $\Sigma_1-\beta$.  We check that the
  one city lies in one component and the
  other city lies in the other component.
  \endproof

  Our results above give us the combinatorial picture
  shown on the right side of Figure 1.4.   We have
  established all the claimed properties about
  the cities in $\Sigma_1$.
    This establishes
  that the dynamical system in Theorem \ref{main}
  describes $\widehat {\cal G\/}$ in the state
  $\Sigma_1$.

  \subsection{Technical Details}
  \label{tech}

  Now we prove the lemmas whose we deferred in the previous
  section.  For convenience we repeat the lemmas.

\begin{lemma}
  The functions $g_{8316}$ and $g_{8319}$ are cleanly related to
  $\Sigma_1$ and their gradients are linearly independent at each point of
  $\Sigma_1$.
\end{lemma}

\startproof
We check that $\partial_L g_{8316}$ does not
vanish in $[0,1]^2$, a region which
contains $\Sigma_1$.
This is Property 1 above. We use the
calculus method to check Property 2 for
$g_{8316}$.  The method works in an obvious
way, with the constant term dominating the others,
as it did for $g_{1638}$ above.

Now we consider $g_{8319}$.
Let $\psi= \partial_{(1,1)} g_{8319}$.
Let $F_1,F_2,F_3$ be the triangle
maps defined in the previous chapter.
We check that $\psi \circ F_j$ is
strongly positive dominant for
$j=1,2,3$.  This checks that
$\psi$ does not vanish in $\Sigma_1$.  This is
Property 1.
We use the evaluation method to check Property 2
for $g_{8319}$.

Let $$g=\det(\nabla g_{8316},\nabla g_{8319}).$$  This
function vanishes if and only if
two gradients are parallel.
We check that the $3$ functions
$g \circ F_j$ are strongly positive dominant.
This does it.
\endproof

\begin{lemma}
  $h_{8369}$ is positive on $\Sigma_1^0$ and
 cleanly related to
    $\Sigma_1$.
\end{lemma}

\startproof
Let $F_1,F_2,F_3$ be the same triangle
maps above.  For each polynomial
$Q_j=h_{8369} \circ F_j$ and some sign choice
$\epsilon_j$ we check that
the $8$ functions
$$\epsilon_jQ_j(x/2,y/2), \hskip 4 pt
  \epsilon_jQ_j(1-x/2,y/2), \hskip 4 pt
  \epsilon_jQ_j(x/2,1-y/2), \hskip 4 pt
  \epsilon_jQ_j(1-x/2,1-y/2),$$
  $$
  \epsilon_jQ_j(1/2,y/2), \hskip 4 pt
  \epsilon_jQ_j(1/2,1-y/2), \hskip 4 pt
  \epsilon_jQ_j(x/2,1/2), \hskip 4 pt
  \epsilon_jQ_j(1-x/2,1/2)$$
  are solidly positive dominant.
  We also check that $\epsilon_jQ_j(1/2,1/2)>0$
  in all $3$ cases.  This shows
  that $h_{8369}$ is positve on $\Sigma_1^o$.
  
  We check the preferred directional derivative
  $\partial_L g$ and we use the calculus
  method for the other $4$ directions.
  In all cases, the constant term is
  much larger than the other terms and the
  result is numerically obvious.
  \endproof

  \begin{lemma}
    $\gamma_{8369}$ intersects each of 
    $\gamma_{8316}$ and
    $\gamma_{8319}$ once.
  \end{lemma}

  \startproof
  The proof is the same in both cases.
  We give the proof for $\gamma_{8316}$.
  We check that the endpoints are interlaced,
    so we just have to verify the gradient condition.  The
  gradient condition does not quite work, so we have
  to scramble.

  Consider the function $g_{8319}$.  We note that
  $\Sigma_1 \subset [0,1] \times [0,1/2]$. We check
  that $\epsilon g_{8369}(1-x/4,y/2)$ is strongly positive dominant
  for some sign choice $\epsilon$.  This
  means that $g_{8369} \not =0$ on $[3/4,1] \times [0,1/2]$.
  Hence $$\gamma_{8369} \subset [0,3/4] \times [0,1/2].$$
  We set
  $\psi={\rm det\/}(\nabla g_{8369},\nabla g_{8316})$
  and check that $\psi(3 x/4,y/2)$ is positive dominant.
  This shows that the two gradients are linearly
  independent in the rectangle $[0,3/4] \times [0,1/2]$.
  Hence $\gamma_{8369}$ and $\gamma_{8316}$ intersect
  exactly one.
  \endproof
  
  \begin{lemma}
    $\gamma_{8369}$ intersects
    $\gamma_{8346}$ and
    $\gamma_{8349}$ once each.
  \end{lemma}

  \startproof
  The proof is the same in both cases. We give
  the proof for $\gamma_{8346}$.
  The endpoints in each case are interlaced and so we
  just have to check the linear independence of the gradients.
  We observe that  $\Sigma_1 \in R=[3/10,1] \times [0,1/2]$
  so we just have to check the linear independence
  in this smaller rectangle. Define
  $\psi={\rm det\/}(\nabla g_{8369},\nabla g_{8346})$.
    check that $\epsilon \psi(1-7 x/10,y/2)$ is strongly positive dominant
  for one of the sign choices.  This shows that
  $\psi$ does not vanish on $R$.
  \endproof

  \newpage

\section{Formulas}

\subsection{Transplant Codes}
\label{code}

This section concerns chains whose sequence ends in $11$.
These correspond to geodesic segments which connect a point
in $\Pi$ with a point in $A(\Pi)$.   Given a chain $C$ there
is a $6$ digit {\it transplant code\/}
$(c_0,...,c_5)$ with the following property.
Given $p \in \Pi$ the point $A(p) \in A(\Pi)$ develops
out to the point
\begin{equation}
  \label{transcode}
 \langle C,p \rangle= \sum_{k=0}^4 c_k \exp(2 \pi i k/5) + \exp(\pi i c_5/5) \overline p.
\end{equation}
Here, as usual, we identity $\Pi$ with the pentagon in $\C$
whose vertices are the $5$th roots of unity.
We call $\langle C,p \rangle$ the
{\it transplant\/} of $p$ with respect to $C$.

In our proof of the Correspondence Lemma, we considered
$6 + 60$ chains.  The first $6$ chains let us define
the vertices $p_1,p_3,p_4,p_6,p_8,p_9$ of
the hexagon $H_p$.  The other $60$ chains are the
competing chains which we eliminate.
The point $p_j$ is give by $\langle c_j,p \rangle$, where
\begin{itemize}
  \item $c_1: (0,3,3,1,0,7)$.
  \item $c_3: (3,3,1,0,0,3)$.
  \item $c_4: (3,3,0,0,1,1)$.
  \item $c_6: (3,0,0,1,3,7)$.
  \item $c_8: (0,0,1,3,3,3)$.
  \item $c_9: (0,0,3,3,1,1)$.
\end{itemize}

 Here we $6$ of the sequences corresponding
to the $60$ competing chains and the
corresponding transplant codes.
\begin{itemize}
\item $0,2,1,9,11 \to 2,4,3,0,1,4$
\item $2,10,9,11 \to 1,4,3,2,0,6$
\item $0,3,2,9,11 \to 1,3,4,2,0,6$
\item $0,3,10,9,11 \to 0,2,4,3,1,8$
\item $0,3,2,10,9,11 \to 2,2,5,3,0,7$
\item $0,4,3,10,9,11 \to 0,2,3,5,2,9$.
\end{itemize}
We can deduce the remaining transplant codes
by symmetry.
Given a chain $C$ we define $C^{\#}$ to be the mirror image of $C$.
We define $\omega C$ to be the chain whose developing image
is obtained from that of $C$ by multiplying the whole picture
by $\exp(2 \pi i/5)$.  For instance if the sequence associated
to $C$ is $0,2,9,11$ then the sequence associated to
$\omega C$ is $0,3,10,11$. Given any chain $C$ we have
the $10$  chains $\omega^k C$ and $\omega^k C^{\#}$ for
$k=0,1,2,3,4$. We call these new chains the
{\it dihedral images\/} of $C$.  Here are the rules for figuring out the
transplant codes for the dihedral images.  Assume that $C$ has
transplant code $c_0,...,c_6$ as above. Then...

\begin{enumerate}
\item $C^{\#}$ has transplant code $c_2,c_1,c_0,c_4,c_3,8-c_5$.
\item $\omega C$ has transplant code $c_4,c_0,c_1,c_2,c_3,c_5+4$.
\end{enumerate}

\subsection{Triangle Maps}
\label{trianglemap}

In the proof of the
Comparison Lemma and the Voronoi Structure Lemma,
we relied on certain polynomial maps
from $[0,1]^2$ to certain triangles $\Upsilon$
We call these the {\it triangle maps\/}.
Here are the domains
\begin{enumerate}
\item $\Upsilon_0=\Upsilon_{1638}$.  This contains the cities $C_{834}$ and $C_{839}$.
\item $\Upsilon_{8349}$. This contains the cities $C_{831}$ and $C_{836}$.
\item $\Upsilon_{8316}$. This contains the cities $C_{163}$ and $C_{168}$.
\item $\Upsilon_1$.
\item $\Upsilon_2$.
\item $\Upsilon_3$.
\end{enumerate}
Together, $\Upsilon_1, \Upsilon_2,\Upsilon_3$ partition the state $\Sigma_0$.

Let $\Upsilon$ be one of the triangles of interest to us.
We want to construct a surjective polynomial map
$F: [0,1]^2 \to \Upsilon$.  In all case let
$T_0$ be the triangle with vertices
$(0,0)$, $(1,0)$ and $(1,1)$.  
We write $F=f_1 \circ f_2$, where
\begin{itemize}
\item $f_1$ is an affine map from the triangle $T_0$ to $\Upsilon$.
\item $f_2(x,y)=(x,xy)$.
\end{itemize}

To define $F$ in each case, we just need to write down the
affine map $f_1$ we use in each case.  Here are the $6$ maps,
listed in the same order as the corresponding triangles.
The various versions of $f_1$ send $(x,y)$ to ...

\begin{equation}
\bigg(\frac{x}{\sqrt{5}+1},\frac{\sqrt{5-2 \sqrt{5}} y}{\sqrt{5}+1}\bigg).
\end{equation}

\begin{equation}
\bigg(\frac{2 \sqrt{5} x-\sqrt{5} y+2}{2 \sqrt{5}+2},
\frac{\sqrt{10-2 \sqrt{5}} \left(-2 x+\left(\sqrt{5}+2\right) y+2\right)}{4 \left(\sqrt{5}+3\right)}\bigg).
\end{equation}

\begin{equation}
\bigg(x+\frac{1}{8} \left(\sqrt{5}-7\right) y,\frac{1}{4} \sqrt{\frac{1}{2} \left(5-\sqrt{5}\right)} y\bigg),
\end{equation}

\begin{equation}
\left(\frac{2 \sqrt{5} x-\sqrt{5} y+2}{2 \sqrt{5}+2},\frac{\sqrt{10-2 \sqrt{5}} \left(-2 x+\left(\sqrt{5}+2\right) y+2\right)}{4 \left(\sqrt{5}+3\right)}\right)
\end{equation}

\begin{equation}
\left(\frac{-4 \sqrt{5} x+\sqrt{5} y-y+4 \sqrt{5}+4}{4 \sqrt{5}+4},\frac{\sqrt{10-2 \sqrt{5}} \left(4 x+\left(\sqrt{5}-1\right) y\right)}{8 \left(\sqrt{5}+3\right)}\right)
\end{equation}

\begin{equation}
\left(\frac{1}{8} \left(\left(\sqrt{5}-5\right) x-2 y+8\right),\frac{1}{16} \sqrt{10-2 \sqrt{5}} \left(\sqrt{5} x+x-\sqrt{5} y+y\right)\right)
\end{equation}

\subsection{Triple Points}
\label{triple}

We will need to compute various triple points
$(ijk;p)$ where $$i,j,k \in \{1,3,4,6,8,9\}.$$
The point $(ijk;p)$ is equidistant from
$p_i,p_j,p_k$.  There is a rational expression
which computes this, but we prefer to use
an alternate approach.

The points
\begin{equation}
p_{ij} = (p_i+p_j)/2, \hskip 30 pt
q_{ij}=p_{ij}+ i (p_i-p_j)
\end{equation}
are both points on the bisector
$b_{ij}$.  We likewise define
$p_{jk}$ and $q_{jk}$.

We define maps $L: \C \to \R^3$ and
$P: \R^3 \to \C$ as follows:
\begin{equation}
L(x+iy)=(x,y,1), \hskip 30 pt
P(x,y,z)=(x/z,y/z).
\end{equation}
The map $PL$ is the identity.  These
maps are familiar from projective geometry.

We have
\begin{equation}
(ijk;p)=P\bigg((L(p_{12}) \times L(q_{12})) \times (L(p_{23}) \times L(q_{23}))\bigg).
\end{equation}
Here $(\times)$ denotes the vector cross product.

\subsection{Formulas for the City Boundaries}
\label{formula}

  We give formulas for the curves considered above.
  Every formula for a city edge can be obtained from
  the ones below by pre-composing these formulas with
  a dihedral symmetry of the pentagonn $\Pi$.
  We have the relation
  $g_{1683}(x,y)=g_{8316}(x,-y)$, so we won't give the formula explicitly
  for $g_{1683}$. This leaves us with the quadruples above which begin
  with $83$.

As we mentioned in the introduction, these functions all have the form
  \begin{equation}
    \label{FORM2}
    \sum_{i+j \leq 3} \bigg(s_{ij} \sqrt{a_{ij}+b_{ij}\sqrt 5}\bigg) x^iy^j=0, \hskip 30 pt
    s_{ij} \in \{-1,0,1\}, \hskip 10 pt
    a_{ij},b_{ij} \in \Z.
  \end{equation}
  We will supply the matrices $\{s_{ij}\}$ and $\{a_{ij}\}$ and $\{b_{ij}\}$ in all
  the relevant cases.
  Since it is easy to mix up matrices with their transposes,
    let me say explicitly that the top horizontal row corresponds to the
    monomials $1,x,x^2,x^3$.  With that said, here is the data for
    $g_{8316}$.

{\tiny
\begin{equation}
  \left[\begin{matrix} + & - & - & +\cr + & + & - & 0 \cr -&+&0&0\cr - & 0 & 0 & 0 \end{matrix}\right]
  \hskip 10 pt
  \left[ \begin{matrix} 50 & 585 & 225 & 60 \cr 283 & 12 & 8 & 0 \cr 25 & 40 &0 &0 \cr 8&0&0&0 \end{matrix}\right]
  \hskip 10 pt
 \left[  \begin{matrix} 20 & 171 & -99 & -16 \cr 105 & -4 & 0 & 0 \cr 0 & 40 &0 &0 \cr 8&0&0&0 \end{matrix}\right]
\end{equation}
\/}
Here are the matrices for $g_{8349}$:
{\tiny
\begin{equation}
  \left[\begin{matrix} + & - & - & +\cr + & - & + & 0 \cr -&+&0&0\cr + & 0 & 0 & 0 \end{matrix}\right]
  \hskip 10 pt
  \left[ \begin{matrix}   940 & 9480 & 235 & 60 \cr 7780 & 100 &20 &0 \cr 75 &60 &0& \cr 20 &0&0&0 \end{matrix}\right]
  \hskip 10 pt
 \left[  \begin{matrix} 420 & 4400 & 105 & 20 \cr 3476& 44 & -4 & 0 \cr 25 & 20 &0 &0 \cr -4&0&0&0 \end{matrix}\right]
\end{equation}
\/}
Here are the matrices for $g_{8319}$:
{\tiny
\begin{equation}
  \left[\begin{matrix} 0 & - & + & -\cr + & + & - & 0 \cr -& - &0&0\cr - & 0 & 0 & 0 \end{matrix}\right]
  \hskip 10 pt
  \left[ \begin{matrix}   0 & 85 & 130 & 5 \cr 29 & 206 &9 &0 \cr 20 &5 &0& \cr 9 &0&0&0 \end{matrix}\right]
  \hskip 10 pt
 \left[  \begin{matrix} 0 & 38 & 58 &2  \cr 12 & 90 & 4 & 0 \cr 8 & 2 &0 &0 \cr 4&0&0&0 \end{matrix}\right]
\end{equation}
\/}
Here are the matrices for $g_{8346}$:
{\tiny
\begin{equation}
  \left[\begin{matrix} 0 & - & + & +\cr + & + & + & 0 \cr - & + &0&0\cr + & 0 & 0 & 0 \end{matrix}\right]
  \hskip 10 pt
  \left[ \begin{matrix}   0 & 126075 & 58835 & 16810 \cr 109265 & 336200 &16810 &0 \cr 294175 &5 &0& \cr 16810 &0&0&0 \end{matrix}\right]
  \hskip 10 pt
 \left[  \begin{matrix} 0 & 42045 & 25215 &0  \cr 31939 & 26896 & 6724 & 126075 \cr 0 & 2 &0 &0 \cr 6724&0&0&0 \end{matrix}\right]
\end{equation}
\/}
Here are the matrices for $g_{8369}$:
{\tiny
\begin{equation}
  \left[\begin{matrix} + & - & - & +\cr + & + & - & 0 \cr -&+&0&0\cr - & 0 & 0 & 0 \end{matrix}\right]
  \hskip 10 pt
  \left[ \begin{matrix} 4700 & 41160 & 210 & 175 \cr 14180 & 1220 & 5 & 0 \cr 2010 & 175 &0 &0 \cr 5&0&0&0 \end{matrix}\right]
  \hskip 10 pt
 \left[  \begin{matrix} 2100 & 18400 & 80 & 75 \cr 6316 & 524 & 1 & 0 \cr 880 & 75 &0 &0 \cr 1&0&0&0 \end{matrix}\right]
\end{equation}
\/}
Here are the matrices for $h_{8369}$:
{\tiny
  \begin{equation}
    \label{h8369}
    \left[\begin{matrix}  0 & - & + & + \cr + & - & + & 0 \cr  - & + & 0 & 0 \cr + &0 &0 &0 \end{matrix}\right] \hskip 10 pt
    \left[\begin{matrix} 0& 196800& 103040& 15040  \cr 1056320& 2044160& 28480& 0 \cr  515200& 15040& 0& 0 \cr 28480& 0& 0& 0 \end{matrix} \right] \hskip 10 pt
   \left[\begin{matrix}  0& 88000& 46080& 6720 \cr 472384& 914176& 12736& 0 \cr 230400& 6720& 0& 0 \cr 12736& 0& 0& 0 \end{matrix} \right]
  \end{equation}
  \/}

\subsection{Formulas for the Special Points}

The quadruple point in $\Sigma_0^o$ is
$(\cos(\pi/5)t,\sin(\pi/5)t)$, where  $t=0.25016...$ is a root of the following cubic.
    $$(5 + 3 \sqrt 5) + (-24 - 10 \sqrt 5) t + (-5 + \sqrt 5) t^2 + 4t^3.$$
In the case of the triple points in $\Sigma_1^o$,
I don't know how to prove that the formulas I got from
Mathematica are correct, but I list them anyway.
In the equations below, the list $(a_0,...,a_{10})$ stands for the
polynomial $$a_0+a_1t + ... + a_{10}t^{10}.$$
The two triple points in $\Sigma_1^o$ are the
$$(\cos(2 \pi/5) u_1,\sin(2\pi/5) v_1), \hskip 30 pt
(\cos(2 \pi/5) u_2,\sin(2\pi/5) v_2),$$
where
$u_1,v_1,u_2,v_2$ respectively are the roots of the following
polynomials.
\newline

{\tiny
$$(316255, -1021235 , 1187259 , 628411 , -2861623 , 
  3126530 , -1726141 , 440390 , -15077 , -8998 , 
  604).$$
  
$$(-495, 9045 , -59511, 
  170103 , -171269 , -112328 , 339489 , -267720 , 
  108905 , -25870 , 3020)$$

  $$(-1044164, 
  4232724 , -10713465 , 20137044 , -23128795 , 
  14627289 , -5047850 , 960889 , -66285 , -10636 ,
  724 )$$
  
  $$(-3820, 14590 , 3825 , -149495 , 131854 ,
   97712 , -165546 , -51200 , 15 , -10750 , 
   3620 10)$$
   \/}

To specify the roots exactly it is enough to note that
$$u_1=1.4799...., \hskip 20 pt
v_2=.21542..., \hskip 20 pt u_2=1.4984..., \hskip 20 pt  v_2=.23169...$$

These $4$ degree $10$ polynomials are irreducible, and Sage tells us that
their Galois groups are all degree $2$ extensions of $S_5 \times S_5$ where $S_5$ is the symmetric
group on $5$ symbols. Hence the coordinates for these triple points are
not solvable numbers.    I also found the polynomials for $x_1,y_1,x_2,y_2$.
The formulas for $x_1,x_2$ are similar to the ones for
$u_1,u_2$. The formulas for $y_1,y_2$ are degree $20$ even polynomials with
enormous integer coefficients.

\newpage

\section{References}

\noindent
[{\bf R1\/}] J. Rouyer, {\it Antipodes sur le t\'etra\`edra r\'egulier\/}, J. Geom. {\bf 77\/} (2003), no. 4, pp. 152-170.
\newline
\newline
[{\bf R2\/}] J. Rouyer, {\it On antipodes on a convex polyhedron\/}, Adv. Geom. {\bf 5\/} (2005), no. 4, pp. 497-507.
\newline
\newline
[{\bf R3\/}] J. Rouyer, {\it On antipodes on a convex polyhedron II\/}, Adv. Geom. {\bf 10\/} (2010), no. 3, pp. 403-417.
\newline
\newline
[{\bf S1\/}] R. E. Schwartz, {\it The Farthest Point Map on the Regular Octahedron\/}, J. Experimental Mathematics, 2021 (to appear).
\newline
\newline
[{\bf S2\/}] R. E. Schwartz, {\it The Projective Heat Map\/},  A.M.S. Research Monograph (2016)
\newline
\newline
[{\bf V1\/}] C. Vılcu, {\it On two conjectures of Steinhaus\/}, Geom. Dedicata {\bf 79\/} (2000), no. 3, pp. 267-275.
\newline
\newline
[{\bf V2\/}] C. Vılcu, {\it Properties of the farthest point mapping on convex surfaces\/}, Rev. Roum. Math. Pures Appl. {\bf 51\/} (2006),
no. 1, pp. 125-134. 
\newline
\newline
[{\bf VZ\/}] C. Vılcu and T. Zamfirescu, {\it Multiple farthest points on Alexandrov surfaces\/}, Adv. Geom. {\bf 7\/} (2007), no. 1, pp. 83-100.
\newline
\newline
[{\bf W\/}] Z. Wang, {\it Farthest Point Map on a Centrally Symmetric Convex Polyhedron\/}, Geometriae Dedicata {\bf 204\/} (2020),
pp. 73-97.
\newline
\newline
[{\bf Wo\/}] S. Wolfram, {\it Mathematica\/} (2020) wolfram.com/mathematica.
\newline
\newline
[{\bf Z\/}] T. Zamfirescu, {\it Extreme points of the distance function on a convex surface\/}, Trans. Amer. Math. Soc. {\bf 350\/} (1998),
no. 4, pp. 1395-1406.

\end{document}